%% file: WENOZD.tex
\documentclass[english,a4paper,12pt,fleqn,nonatbib,preprint]{elsarticle}
\usepackage{geometry}
\usepackage{amsmath,amsthm,amsfonts,amssymb}
\usepackage{mathtools}
\usepackage{yhmath}
\usepackage{array}
\usepackage{float}
\usepackage{units}
\usepackage{multirow}
\usepackage{booktabs}
\usepackage{graphicx}
\usepackage{subcaption}
\usepackage{xcolor}
\usepackage{soul}
\usepackage{csquotes}
\usepackage[english]{babel}
\makeatletter
\let\c@author\relax
\makeatother
\usepackage[maxbibnames=99]{biblatex}
\usepackage{hyperref}

\def\compilapesado{1}

\newcommand{\graficopesado}[2]{
  \centering
  \fbox{
    \begin{minipage}[c][0.33\textheight][c]{0.5\textwidth}
      \centering{Overleaf não me deixa compilar \\ com esta figura}
    \end{minipage}
  }
}

\if\compilapesado1
   \renewcommand{\graficopesado}[2]{\includegraphics[#1]{#2}}
\fi



\sethlcolor{yellow}

\graphicspath{{Images/arXiv/}}
\input{commands.tex}

\begin{document}

\newpageafter{abstract}

\begin{frontmatter}
\title{New Weighting Strategies for WENO Schemes}
\author[1,4]{Daniel Barreto}
\ead{danielw@dme.ufrj.br}

\author[2,4]{Rafael B. de~R. Borges}
\ead{rafael.borges@uenf.br}

\author[3,4]{Bruno Costa\texorpdfstring{\corref{cor1}}{}}
\ead{bcosta@im.ufrj.br}

\author[1,4]{Silvaneo dos~Santos}
\ead{silvaneo@dme.ufrj.br}

\cortext[cor1]{Corresponding author}

\affiliation[1]{
   organization={Departamento de Métodos Estatísticos, Universidade Federal do Rio de Janeiro},
   city={\\Rio de Janeiro/RJ},
   country={Brazil}
}
\affiliation[2]{
   organization={Laboratório de Ciências Matemáticas, Universidade Estadual do Norte Fluminense Darcy Ribeiro},
   city={\\Campos dos Goytacazes/RJ},
   country={Brazil}
}
\affiliation[3]{
   organization={Departamento de Matemática Aplicada, Universidade Federal do Rio de Janeiro},
   city={\\Rio de Janeiro/RJ},
   country={Brazil}
}
\affiliation[4]{
   organization={Laboratório de Matemática Aplicada, Universidade Federal do Rio de Janeiro},
   city={\\Rio de Janeiro/RJ},
   country={Brazil}
}

\begin{abstract}
In this article, we propose a modified convex combination of the polynomial reconstructions of odd-order WENO schemes to maintain the central substencil prevalence over the lateral ones in all parts of the solution. New ''centered'' versions of the classical WENO-Z and its less dissipative counterpart, WENO-Z+, are defined through very  simple modifications of the classical nonlinear weights and show significantly superior numerical properties; for instance, a well-known dispersion error for long-term runs is fixed, along with decreased dissipation and better shock-capturing abilities. Moreover, the proposed centered version of WENO-Z+ has no ad-hoc parameters and no dependence on the powers of the grid size. All the new schemes are thoroughly analyzed concerning convergence at critical points, adding to the discussion on the relevance of such convergence to the numerical simulation of typical hyperbolic conservation laws problems. Nonlinear spectral analysis confirms the enhancement achieved by the new schemes over the standard ones.
\end{abstract}
\end{frontmatter}


\tableofcontents\newpage


\input{1-Introduction}
\input{2-Review-WENO}
\input{2.1-Convergence-WENO}
\input{2.2-Dispersion-Error}
\input{3-The-WENO-C}

\input{3.1-Convergence-Analysis-WENO-C}
\input{4-WENOZD-and-ZDplus}
\input{4.1-WENOZDplus}

\input{4.2-ADR-Analysis}
\input{4.3-Weight-Distribution}
\input{5-Numerical-Experiments}
\input{6-Conclusions}

\newpage

\printbibliography

\newpage

\appendix
\input{A-Taylor-Series}
\input{B-Proof-Conditions}
\input{C-Convergence-WENO-ZDplus}

\end{document}

%% file: commands.tex
\newcommand{\dx}{{\Delta x}}

\newcommand{\fhat}{{\hat f}}
\newcommand{\half}[1][1]{\frac{#1}{2}}
\newcommand{\cp}{\mathrm{cp}}
\newcommand{\ncp}{n_\cp}
\newcommand{\eps}{\varepsilon}

\def\wmin*{\textit{min}}

\DeclareMathOperator{\Ord}{O}

\newcommand{\percentdone}[1]{ [{#1}\%]}  \renewcommand{\percentdone}[1]{}

\newcommand{\opcional}[1]{}
\newcommand{\opcionalfigura}[1]{}

\newtheorem{condition}{Condition}
\theoremstyle{definition}
\newtheorem{definition}{Definition}
\theoremstyle{remark}
\newtheorem*{remark}{Remark}

%% file: 1-Introduction.tex
\section{Introduction}

WENO schemes have been proposed to overcome the Gibbs-like oscillations when dealing with solutions of shock-turbulence interactions containing large gradients and discontinuities in the numerical simulation of hyperbolic conservation laws. The main idea is to
use a nonlinear convex combination of lower-order approximation polynomials that adapts either to a higher-order approximation in smooth regions of the solution, or to a lower-order spatial discretization that avoids interpolation across discontinuities and provides the necessary numerical dissipation for shock capturing. See \cite{shu97,shu09} and references therein for a detailed introduction.

The need for stable high-order numerical methods without introducing extra nonphysical oscillations has inspired the developments of total variation diminishing (TVD) schemes \cite{harten83}, the essentially nonoscillatory (ENO) schemes \cite{harten87b}, the weighted essentially nonoscillatory (WENO) schemes \cite{liu94,jiang96,balsara00,henrick05,borges08,castro11}, and so on. 
In other words, WENO schemes belong to the search of a method that both correctly represents fine and smooth structures and capture shock discontinuities in the numerical simulation of inviscid compressible flows. This all started back in 1983, with the TVD and ENO schemes, passing through the first WENO scheme by Liu, Osher, and Chan in 1994 \cite{liu94}, and getting to the smoothness indicators of the WENO scheme proposed by Jiang and Shu, hereafter dubbed WENO-JS, in 1996 \cite{jiang96}. 

These WENO schemes used so-called smoothness indicators, traditionally named as $\beta_{k}$, which satisfy \(\beta_{k} = \Ord(1)\) if the function has a jump discontinuity in the substencil \(S_k\), and \(\beta_{k} = \Ord(\dx^q)\) (i.e., \(\beta_{k}\) is small) otherwise. 
The weights of the convex combination decrease with $\beta_{k}$, also decreasing the importance of substencils containing discontinuities and/or high gradients, and at smooth regions of the solution, the weights must mimic ideal weights generating a higher order approximation.

Henrick et al. \cite{henrick05} pointed out a deficiency of WENO-JS at critical points 
which they were able to correct by using a mapping function (WENO-M). Nevertheless, it ended up being an expensive fix that Borges et al. \cite{borges08,castro11}, by introducing a global smoothness indicator, $\tau$, were able to work around with a lower computational cost (WENO-Z) \cite{zhao14}. Many other WENO-Z-like variants were proposed aiming higher order schemes as well as less dissipative schemes to better capture the finer features and shocks of the solutions, see for instance \cite{yamaleev09b, hu10, ha13, vanlith17, wang19, hong20}. In particular, Acker et al. \cite{acker16} proposed the addition of an extra anti-dissipative term to improve the representation of curvature features of the numerical solution (WENO-Z+). However, WENO-Z+ shows overamplification of some features of the solution that depends on an ad hoc grid size parameter \cite{luo21a,luo21b}.

\begin{figure}[H]
    \centering
    \includegraphics[width=\textwidth]{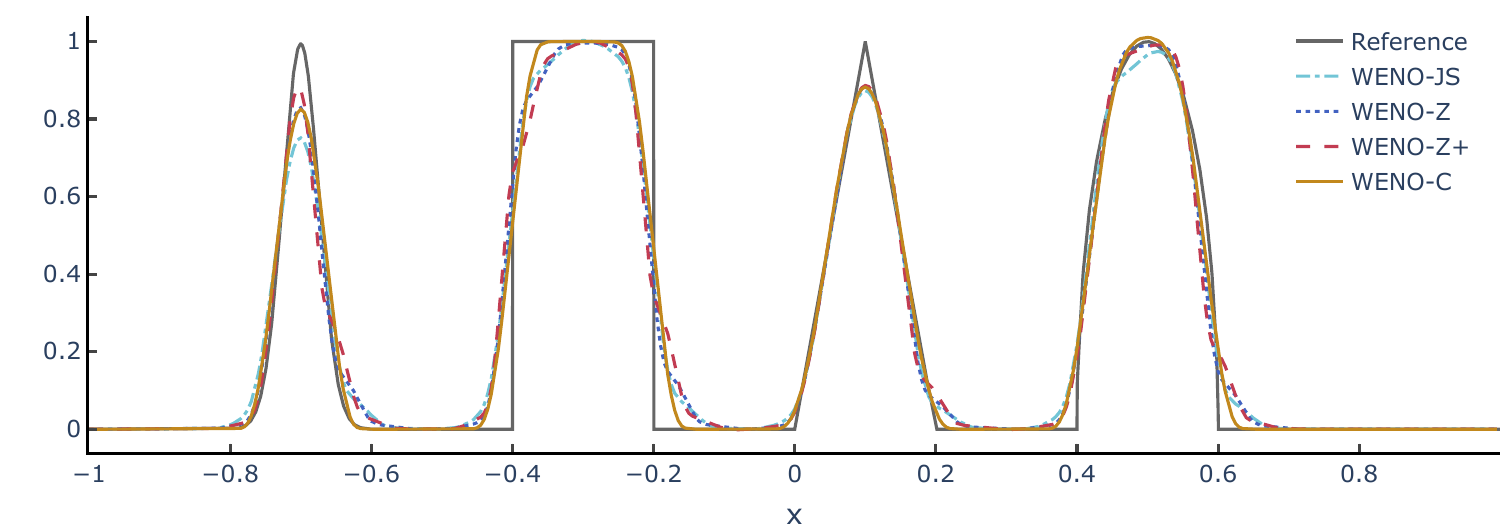} \\
    \includegraphics[width=\textwidth]{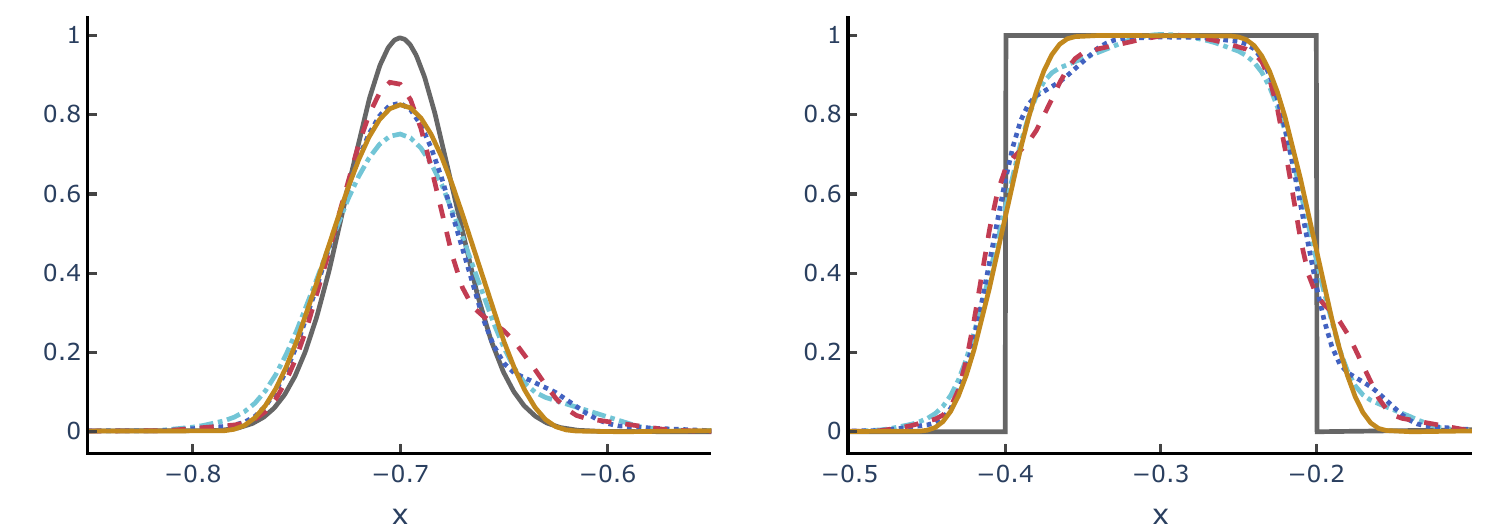}
    \caption{Numerical solutions of the Gaussian-Square-Triangle-Ellipse (GSTE) advection test at \(t=100\) with $N=400$ and CFL = $0.45$.}
    \label{fig:intro 1a}
\end{figure}

The lack of convergence at critical points is not the only obstacle faced by most of the WENO schemes; a dispersion error for long-term runs is also present and it was noticed by some researchers \cite{hong20,don23,li23}.  One of the main novelties presented in this article is the building of new WENO schemes that stay much closer to the structure of the central upwind scheme to avoid dispersion errors, which we credit to the weakening of the information of the central substencil in the WENO convex combination. For that,  we devised new weighting strategies, where the central substencil keeps its relevance in relation to the lateral ones in a way that still maintains the nonoscillatory property.
This new centered configuration, dubbed WENO-C, also improved on dissipative behavior  when compared with the classical WENO-JS, WENO-M, and WENO-Z and, as shown in Figure  \ref{fig:intro 1a}, it fixed the dispersion error existing in these schemes. 

Nevertheless, WENO-C also needs adjustment when dealing with critical points of the solution, and the idea is to use a convergence accelerator term that acts on the presence of these points. Among the many choices to do this, the most appropriate form  was found in relation to a parallel investigation of searching for a better anti-dissipative term for WENO-Z+, one that prevents the overamplification referred to above. The new scheme, consisting of the coupling of WENO-C and  this accelerating term,  was named  WENO-ZC, and it recovered optimal convergence on critical points, maintaining all the previous enhancements of WENO-C with respect to WENO-Z. The quest for an improved WENO-Z+ was also successful through the same centered structure proposed in this article; however, it no longer suffers from the overamplification of curvature features as the original WENO-Z+. Moreover, the anti-dissipative term of the new scheme, WENO-ZC+, involves no extra parameter nor any dependence on any power of the grid size, and enhanced dissipation and shock-capturing capabilities are easily seen in the classical numerical experiments.

\begin{figure}[ht]
    \centering
    \begin{subfigure}[b]{0.5\textwidth}
         \centering
         \includegraphics[width=\textwidth]{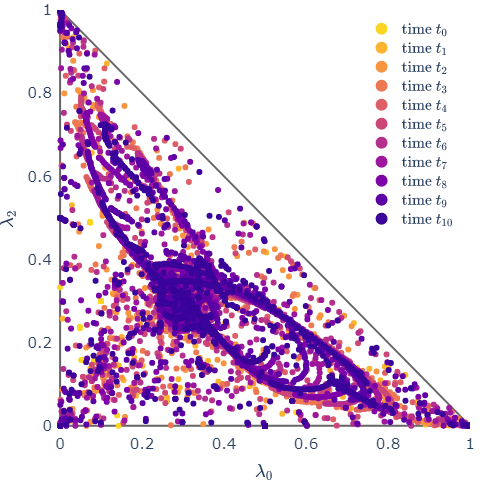}
         \caption{WENO-JS}
     \end{subfigure}\hfill
    \begin{subfigure}[b]{0.5\textwidth}
         \centering
         \includegraphics[width=\textwidth]{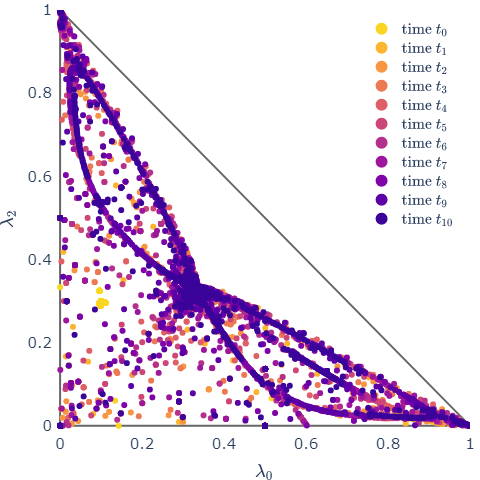}
         \caption{WENO-Z}
     \end{subfigure}\hfill
     \caption{Distribution map of the weights for the Titarev--Toro problem at different times \(t_i = i/2\), with $N=1000$ and CFL = $0.5$.}
     \label{fig:intro 1b}
\end{figure}

Observing the dynamics of the nonlinear components of the WENO weights also generates further understanding. For instance, Figure \ref{fig:intro 1b} shows the temporal evolution of the nonlinear components of the unnormalized weights, divided by the ideal weights, of the lateral substencils, $\lambda_0$ and $\lambda_2$. The size of $\lambda_1$, the one corresponding to the central substencil is deduced through the normalization $\lambda_0 + \lambda_1 +\lambda_2=1$. In this way, it is easily seen  in Figure \ref{fig:intro 1b} that WENO-Z puts a floor to the size of $\lambda_1$, whereas WENO-JS shows convex combinations where  $\lambda_0$ and $\lambda_2$  are  both  large  and their sum, $\lambda_0 + \lambda_2$, is very close to $1$. In this article, we also present such  comparisons to obtain new insights into the inner workings of the new WENO schemes in  relation  to the classical ones.

The convergence of the new WENO schemes at critical points, as well as the proximity of its nonlinear weights to the ideal weights, is thoroughly discussed throughout the article concerning the practical aspects obtained on the numerical tests. For this matter, the WENO-D method \cite{wang19} is brought into the comparison, for it is a scheme that achieves optimal order of convergence  at critical points of any order. Nevertheless, its numerical properties are similar to the ones of WENO-Z on typical hyperbolic problems, and this may shed some light on the practical importance of convergence at high-order critical points when dealing with numerical solutions containing shocks and discontinuities.

The remainder of this article is divided as follows: In Section \ref{sec:review-weno}, we review the classical WENO-JS, WENO-M, and WENO-Z schemes, as well as the more recent WENO-Z+ and WENO-D schemes, and discuss two issues these schemes may have: the loss of convergence near critical points and the long-term cumulative dispersion errors in the solutions. In Section \ref{sec:weno-c}, we present the fundamental modification of the classical WENO schemes that fixes the long-term dispersion error when allowing a combination where the central substencil has bigger relevance than the lateral ones. However, this new scheme has an issue near critical points, and in Section \ref{sec:weno-zc}, we improve the centered combination of substencils to achieve the same order of convergence of WENO-Z at critical points of the solution, leading to the novel WENO-ZC and WENO-ZC+ schemes. We also present numerical results and an ADR analysis showing that the resulting schemes have superior spectral  properties in relation to the original WENO-Z and WENO-Z+.
In Section \ref{sec:numerical-experiments}, the enhanced numerical results obtained with the 2D Euler Equations show that the centered WENO schemes are a viable path of investigation for future improvements. Conclusions are drawn in Section \ref{sec:conclusions}.

%% file: 2-Review-WENO.tex
\section{Review of Classical WENO Schemes} \label{sec:review-weno}

In this section, we present a brief summary of WENO schemes and, more specifically, how they are employed for solving conservation laws.  In particular, we make a brief review of the following classical WENO schemes---the WENO-JS \cite{jiang96,balsara00,gerolymos09}, the mapped WENO (WENO-M) \cite{henrick05} and the WENO-Z \cite{borges08,castro11}---as well as the more recent WENO-Z+ \cite{acker16} and WENO-D \cite{wang19} schemes. The present explanation is admittedly terse; the references \cite{shu97, shu09}, as well as the others found in the text, are recommended for a more comprehensive and detailed discussion.

WENO schemes are applied to hyperbolic conservation laws:
\begin{equation}
   \frac{\partial u}{\partial t} + \nabla\cdot F(u) = 0
   \label{eq:Hyperbolic Law}
\end{equation}
where $u$ and $F$ are the vectors representing the conservative variables and the flux functions. Spatial discretization is achieved by an uniform grid $x_{i}=i\Delta x$, $i=0,\,\ldots,\,N$, where\newline $\Delta x=\frac{1}{N-1}$ is the constant grid size. The $x_{i}$ are also called cell centers, with cells defined as $I_{i}=\left[x_{i-\frac{1}{2}},x_{i+\frac{1}{2}}\right]$, where $x_{i\pm\frac{1}{2}}=\left(i\pm\frac{1}{2}\right)\Delta x$ are the cell boundaries.

The semi-discretized form of \eqref{eq:Hyperbolic Law} yields the system of ordinary differential equations:
\[
   \frac{du_{i}(t)}{dt} = -\frac{\partial f}{\partial x}\bigg|_{x=x_{i}}, \quad i=0\ldots N,
\]
where $u_{i}(t)$ is a numerical approximation of $u$ at $\left(x_{i},t\right)$. Once the derivative $\frac{\partial f}{\partial x}|_{x=x_{i}}$ has been evaluated, a time discretization method such as the third-order TVD Runge--Kutta advance the solution in time. The spatial derivative above is exactly computed by a conservative finite difference formula at the cell boundaries
\begin{equation}
   \frac{du_{i}(t)}{dt}=\frac{1}{\Delta x}\left(h_{i+\frac{1}{2}}-h_{i-\frac{1}{2}}\right)\label{eq:conservative}
\end{equation}
where the numerical flux function $h(x)$ is implicitly defined as
\begin{equation}
   f(x)= \frac{1}{\dx} \int_{x-\frac{\dx}{2}}^{x+\frac{\dx}{2}} h(\xi) \: d\xi.
   \label{eq:numerical-flux}
\end{equation}

\begin{figure}
   \centering
   \includegraphics[scale=0.8]{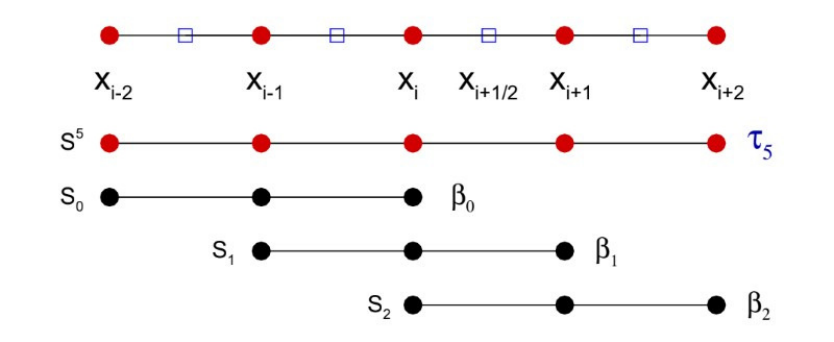}
   \caption{WENO uniform grid and global stencil, $S^{5}$, and its 3 substencils $S_{k}$, $k=0,1,2$, along with their corresponding smoothness indicators $\tau_{5}$ and $\beta_{k}$.}
   \label{fig:WENO-uniform-grid}
\end{figure}

In the classical fifth-order WENO scheme, a global 5-points stencil, $S^{5}$, is subdivided into three 3-points substencils $\left\{ S_{0},S_{1},S_{2}\right\}$, as shown in Fig. \ref{fig:WENO-uniform-grid}. The values $h_{i\pm\frac{1}{2}}$ are computed by fifth-order polynomial interpolations on the grid values of $f$, which are the results of the convex combinations below:
\begin{equation}
   \hat{f}_{i\pm\frac{1}{2}}=h_{i\pm\frac{1}{2}}+\Ord(\Delta x^{5})=\sum_{k=0}^{2}\omega_{k}^{\pm}\hat{f^{k}}(x_{i\pm\frac{1}{2}}).
   \label{eq:convex combination}
\end{equation}
 Here, each $\hat{f^{k}}(x)$ is the second-order Lagrangian polynomial computed with the values of $f$ in each $S_{k}$, satisfying
 \begin{align}
    \fhat^0_{i+\frac{1}{2}} &= \frac{1}{6} (2f_{i-2} - 7f_{i-1} + 11f_{i}), \nonumber \\
    \fhat^1_{i+\frac{1}{2}} &= \frac{1}{6} (-f_{i-1} + 5f_{i} + 2f_{i+1}), \nonumber \\
    \fhat^2_{i+\frac{1}{2}} &= \frac{1}{6} (2f_{i} + 5f_{i+1} - f_{i+2}), \label{eq:fhatk}
 \end{align}
 and $\omega_{k}$ are nonlinear weights which depend on the smoothness of the numerical solution in $S_{k}$. The main idea is to decrease the importance in \eqref{eq:convex combination} of the components residing on a stencil where the solution is not smooth, thus avoiding interpolations through high gradients and their ensuing oscillations.

\subsection{The WENO-JS Scheme}
The first successful WENO method to be proposed, by Jiang and Shu \cite{jiang96}, hereafter denominated WENO-JS, used the following \textit{smoothness indicators} to quantify the regularity of the $k$-th polynomial approximation $\hat{f^{k}}(x)$ at stencil $S_{k}$:
\begin{equation}
   \beta_{k}=\sum_{l=1}^{2}\Delta x^{2l-1}\int_{x_{i-\frac{1}{2}}}^{x_{i+\frac{1}{2}}}\left(\frac{d^{l}}{dx^{l}}\hat{f^{k}}(x)\right)^{2}dx,
   \quad k = 0,1,2.
   \label{eq:betaksintegral}
\end{equation}

The unnormalized and normalized WENO-JS weights were respectively defined as:
\begin{equation}
   \alpha_{k}^{JS} = \frac{d_{k}}{\left(\beta_{k}+\eps\right)^{p}} \quad \text{and} \quad 
   \omega_{k}^{JS}=\frac{\alpha_{k}^{JS}}{\sum_{j=0}^{2}\alpha_{j}^{JS}},
   \label{eq:WENO-JS Weights}
\end{equation}
where $d_{k}$, the \textit{ideal weights}, namely
\begin{equation} \label{eq:ideal-weights}
   (d_0, d_1, d_2) = \left(\frac{1}{10}, \frac{6}{10}, \frac{3}{10}\right),
\end{equation}
are those providing in \eqref{eq:convex combination} the fifth-order Lagrangian interpolation of the original central-upwind scheme:
\[
   \sum_{k=0}^2 d_k \fhat^k(x_{i+\half}) = \fhat(x_{i+\half}) = h(x_{i+\half}) + \Ord(\dx^5).
\]
The parameters $\eps$ and $p$ were respectively designed to avoid division by zero and to increase the difference of scales of distinct weights at nonsmooth parts of the solution. The general idea is that, at smooth parts of the solution, the $\beta_{k}$ are small and with the same magnitude, thus the $\omega_{k}$ are close to the ideal weights \(d_k\), reproducing the fifth-order upwind scheme. However, if one of the $S_{k}$ contains high gradients, then $\beta_{k}=\Ord(1)$ and the corresponding $\omega_{k}$ is small relative to the other weights, decreasing the importance of the component $\hat{f^{k}}$ in the final convex combination \eqref{eq:convex combination}, avoiding interpolation across discontinuities and the ensuing numerical oscillations.

\subsection{The WENO-M Scheme}
Still, in smooth regions, where $\omega_{k}$ would ideally approach $d_{k}$, flat points of the solution degrade the order of this convergence. In other words, too much smoothness, surprisingly, is also a problem for WENO schemes to generate the ideal convex combination of the substencils. In \cite{henrick05}, a mapping $g_k(\omega)$ was proposed to fix this issue, yielding the mapped WENO weights:
\begin{equation}
   \alpha_{k}^{M} = g_k(\omega^{JS}_k) 
   \quad \text{and} \quad 
   \omega_{k}^{M} = \frac{\alpha_{k}^{M}}{\sum_{j=0}^{2}\alpha_{j}^{M}}.
   \label{eq:WENO-M Weights}
\end{equation}
The resulting WENO scheme, here dubbed WENO-M, was able to fix the critical points convergence issue of WENO-JS, although at a high computational cost. Additional convergence issues were soon discovered involving $\eps$ and $p$, and many further works relate on these \cite{arandiga11, don13}, which we shall address later. 

\subsection{The WENO-Z Scheme}
Borges et al. \cite{borges08} proposed a different path to enhance the WENO-JS scheme through the utilization of higher-order smoothness information contained in the five-point global stencil. The WENO-Z weights are given by
\begin{equation}
   \alpha_{k}^{Z} = d_{k}\Bigg[1+\left(\frac{\tau}{\beta_{k}+\eps}\right)^{p}\Bigg]
   \quad \text{and} \quad 
   \omega_{k}^{Z} = \frac{\alpha_{k}^{Z}}{\sum_{j=0}^{2}\alpha_{j}^{Z}},
   \label{eq:WENO-Z Weights}
\end{equation}
where \(\tau\), the global smoothness indicator, measures the smoothness of the solution in the global five-points stencil. In the original work, it was defined as
\[
    \tau = \big|\beta_2 - \beta_0\big|.
\]
The Z-weights are easily seen as a modification of the ideal weights $d_{k}$ through the addition of the nonlinear component,namely, \(\left(\dfrac{\tau}{\beta_k + \eps}\right)^p\). The Z-weights accelerate the convergence to the ideal weights at all regions of the solution, not only at the smooth ones, also decreasing the smearing of shocks and discontinuities.

\subsection{The WENO-Z+ Scheme}
The main idea of the WENO-Z+ is to better capture parts of the solution with curvature features \cite{acker16}. Its weights are a modification of the WENO-Z ones:
\[
   \alpha_{k}^{Z+} = d_{k}\Bigg[1 + \left(\frac{\tau}{\beta_{k}+\eps}\right)^p + \eta \frac{\beta_{k}}{\tau+\eps}\Bigg]
   \quad\text{and}\quad
   \omega_{k}^{Z+} = \frac{\alpha_{k}^{Z+}}{\sum_{j=0}^{2}\alpha_{j}^{Z+}}.
\]
The extra term $\dfrac{\beta_{k}}{\tau + \eps}$ increases the weights for those substencils with high gradients and, at the same time, \(\eta\) limits this growth when discontinuities are present. It also works as an anti-dissipative term and the grid-size-dependent parameter \(\eta\) is used to decrease overamplification. This was the main idea of the WENO-Z+, allowing the substencils with curvature to be better represented in the final convex combination. $\eta = \Delta x^{2/3}$ was defined as a compromise between stability and resolution power for the standard tests suite; however, this value was found in an ad hoc manner and no general result was obtained. Moreover, this makes WENO-Z+ not a self-similar scheme and, therefore, sensible to changes in scaling (see \cite{shu97,don23,li23} for more detailed discussions on the subject of self-similarity).

The numerical experiments with sinusoidal solutions, like Shu--Osher and Titarev--Toro, showed the occurrence of an overamplification of some features, depending on the chosen \(\eta\) \cite{acker16}. Later works have further explored the overamplification tendencies of WENO-Z+ \cite{luo21a,luo21b} and in this article we will propose some changes on the anti-dissipative term in order to prevent that.

\subsection{The WENO-D Scheme}

The WENO-D scheme of Wang, Wang, and Don \cite{wang19} is also a modification of the WENO-Z scheme, defined as
\begin{equation}
   \alpha_{k}^{D} = d_{k}\Bigg[1+\Phi\!\left(\frac{\tau}{\beta_{k}+\eps}\right)^{p}\Bigg]
   \quad \text{and} \quad 
   \omega_{k}^{D} = \frac{\alpha_{k}^{D}}{\sum_{j=0}^{2}\alpha_{j}^{D}},
   \label{eq:WENO-D Weights}
\end{equation}
where
\begin{equation*}
   \Phi = \min\{1, \phi\}, \qquad
   \phi = \sqrt{|\beta_0 - 2\beta_1 + \beta_2|}.
\end{equation*}
Remarkably, and in contrast with all previous WENO schemes, the WENO-D scheme achieves optimal order $5$ even in the presence of very flat critical points due to the use of the term $\Phi$, which works as a convergence accelerator near critical points \cite{wang19}. Along the article, we will use this optimality of WENO-D to analyze the numerical results of the new proposed schemes, particularly when referring to the issue of convergence at critical points.

%% file: 2.1-Convergence-WENO.tex
\subsection{The Classical Analysis of Convergence of WENO Schemes} \label{sec:convergence-weno}

First noticed by Hendrick et al. \cite{henrick05}, slower convergence at critical points is a collateral trait of WENO schemes. WENO schemes avoid numerical oscillations by detecting relative differences, not absolute differences, on gradients of the distinct substencils numerical solutions. Thus, on smooth parts of the solution, specially flat ones, gradients of the solution on distinct stencils may pertain to relatively different scales, although they all have absolute small scales. This situation causes the WENO combination to break away from the ideal central-upwind configuration, decreasing the numerical convergence, particularly on critical points of high order.

In this section, we state the general framework to study the convergence of WENO schemes in the presence of critical points that will permeate the investigation of the convergence of the new WENO schemes proposed in Sections \ref{sec:weno-c} and \ref{sec:weno-zc}. We start by recalling two sets of conditions on the normalized weights \(\omega_k\), deduced in \cite[Eqs. (25) and (29)]{henrick05}, which assure convergence with optimal order:

\begin{condition} \label{cond:sufficient}
   Suppose \(f\) is smooth.  If the weights \(\omega_k^{\pm}\) satisfy
   \begin{align}
      \sum_{k=0}^{2}(\omega_{k}^{\pm}-d_{k}) &= \Ord(\Delta x^{6}), \label{eq:normalizationcond1} \\ 
      \omega_{k}^{\pm}-d_{k} &= \Ord(\Delta x^{3}), \quad k = 0, 1, 2, \label{eq:condition-dx3}
   \end{align}
   then the associated WENO scheme is fifth-order accurate in space.
\end{condition}

\begin{condition} \label{cond:necessary}
   Suppose \(f\) is smooth. A WENO scheme is fifth-order accurate in space if, and only if, its weights \(\omega_k^{\pm}\) satisfy
   \begin{align}
      \sum_{k=0}^{2}(\omega_{k}^{\pm}-d_{k}) 
         &= \Ord(\Delta x^{6}), \label{eq:normalizationcond} \\
      \omega_{k}^{\pm}-d_{k}
         &= \Ord(\Delta x^{2}), \qquad k = 0, 1, 2, \label{eq:idealweightscond} \\[6pt]
      3(\omega_{0}^{+}-\omega_{0}^{-}) &- (\omega_{1}^{+}-\omega_{1}^{-}) + (\omega_{2}^{+}-\omega_{2}^{-}) 
         = \Ord(\Delta x^{3}). \label{eq:Aks}
   \end{align}
\end{condition} 
\noindent Here, \(\omega_{k}^{\pm}\) denotes respectively the weights of \(\fhat^k(x_{i\pm\half})\) in the convex combination.

In the absence of critical points, it has been shown that both WENO-M and WENO-Z satisfy Condition \ref{cond:sufficient} \cite{henrick05,borges08,don13}, and WENO-JS satisfies Condition \ref{cond:necessary} \cite{don13}. The proofs of Conditions \eqref{cond:sufficient} and \eqref{cond:necessary} give more insights about the inner workings of WENO schemes and we present them in \ref{sec:proof-conditions}. We will be referring to these conditions throughout the following sections when analyzing the new WENO schemes; nevertheless, there are some remarks we deem necessary to clarify the following discussion on the convergence of WENO schemes in the presence of flat regions of the solution.

\begin{remark}
It is easily seen that Eqs. \eqref{eq:normalizationcond1} and \eqref{eq:normalizationcond} are automatically satisfied due to normalization. Also, requirement \eqref{eq:condition-dx3} is stronger than  \eqref{eq:idealweightscond} and exempts those WENO schemes satisfying it from satisfying requirement \eqref{eq:Aks}. 
\end{remark}

\begin{remark}
It is also easily verifiable that
\[
   \alpha_k = d_k + \Ord(\dx^q)
   \quad \text{and} \quad
   \omega_k = \frac{\alpha_k}{\sum_{j=0}^2 \alpha_j}
   \quad \implies \quad 
   \omega_k = d_k + \Ord(\dx^{\tilde{q}}), \quad \tilde{q} \geq q.
\]
So, if either condition \eqref{cond:sufficient} or \eqref{cond:necessary} is satisfied for \(\alpha_k\), then it is immediately satisfied for \(\omega_k\). Importantly, it was shown in \cite[Lemma 6]{don13} that if the leading-order term of \(\alpha_k - d_k\) is independent of \(k\), then
\[
   \alpha_k = d_k + \Ord(\dx^q)
   \quad \text{and} \quad
   \omega_k = \frac{\alpha_k}{\sum_{j=0}^2 \alpha_j}
   \quad \implies \quad 
   \omega_k = d_k + \Ord(\dx^{\tilde{q}}), \quad \tilde{q} \geq q+1.
\]
That is, in this case the order of \(\omega_k - d_k\) is at least one unit greater than the order of \(\alpha_k - d_k\).
\end{remark}

In the presence of critical points, the above WENO schemes lose orders of accuracy depending on the values of the parameters \(\eps\) and \(p\) \cite{don13}, as we shall see in the discussion below.
 First, we recall the definition of the order of a critical point:
\begin{definition} \label{def:ncp}
    \(x_i\) is a critical point of order \(\ncp\) of \(f\) if, and only if,
    \begin{equation*}
       f'(x_i) = f''(x_i) = \ldots = f^{(\ncp)}(x_i) = 0, \qquad
       f^{(\ncp + 1)}(x_i) \neq 0.
    \end{equation*}
\end{definition}

It has been shown that WENO-JS lose accuracy at simple critical points (i.e., those of order \(\ncp = 1\)) if the parameter \(\eps\) is too small \cite{henrick05}, achieving order 3 only. In fact, it was shown that, in the presence of critical points, WENO-JS converges with nominal order 5 if, and only if, \(\eps = \Omega(\dx^2)\) (that is, \(\eps \geq C \dx^2\) as \(\dx \to 0\) for some constant \(C > 0\)) \cite{don13}. So, the parameter \(\eps\), which originally had the single role of avoiding a division by zero, can be used for improving the accuracy properties of WENO-JS. However, a large $\eps$ may harm the gradient detection ability of $\beta_{k}$, causing the loss of the nonoscillatory property and subsequent spurious oscillations \cite{don13}.

The WENO-Z scheme has a different fix for its loss of convergence at critical points. For simplicity, consider \(\eps = 0\), at simple critical points, the Taylor series of \(\beta_k\) \eqref{eq:betataylor} and \(\tau\) \eqref{eq:tautaylor} yield
\begin{align}
   \alpha_k^{Z(+)} 
      &= d_{k} + d_{k}\left(\frac{\tau}{\beta_{k}}\right)^{p}
      = d_{k} + d_{k}\left(\frac{\dfrac{13}{3} |f_i''f_i'''| \dx^5 + \Ord(\dx^7)}{\dfrac{13}{12} (f_i'')^2 \dx^4 + \Ord(\dx^5)}\right)^{p}
      = d_{k} + \Ord(\Delta x^{p}), \nonumber \\
   \alpha_k^{Z(-)} 
      &= d_{k} + d_{k}\left(\frac{\tau}{\beta_{k}}\right)^{p}
      = d_{k} + d_{k}\left(\frac{\dfrac{13}{3} |f_i''f_i'''| \dx^5 + \Ord(\dx^7)}{\dfrac{25}{12} (f_i'')^2 \dx^4 + \Ord(\dx^5)}\right)^{p}
      = d_{k} + \Ord(\Delta x^{p}). \label{eq:alpha-z-ncp-1}
\end{align}
For \(p=1\), Conditions \ref{cond:sufficient} and \ref{cond:necessary} are not satisfied; however, for \(p=2\), it has been shown \cite[Lemma 6]{don13} that
\[
   \omega_k^{(\pm)} = d_k + \Ord(\dx^3)
\]
even though we only have \(\alpha_k^{(\pm)} = d_k + \Ord(\dx^2)\). As previously mentioned, this is because the leading-order term in \eqref{eq:alpha-z-ncp-1} does not change with \(k\).

The discussion above shows that either $\eps$, in the case of WENO-JS, or $p$, in the case of WENO-Z, may be used to accelerate convergence at critical points of the solution. 
Nevertheless, for small \(\eps\) and regardless of \(p\), the orders of both WENO-JS and WENO-Z decay to three in the presence of a critical point of order \(\ncp = 2\). On the other hand, WENO-D uses a different convergence accelerator at critical points through the term $\Phi$ in \ref{eq:WENO-D Weights}, which detects through a clever conditional statement if the solution is rough or smooth, strongly decreasing the relevance of the nonlinear components of the weights in the later case.

%% file: 2.2-Dispersion-Error.tex
\subsection{Dispersion Error of the Classical WENO Schemes} \label{sec:dispersion-error}

The loss of accuracy in the presence of critical points is not the only issue WENO schemes may present. For long-time simulations, numerical solutions computed with WENO schemes may suffer from significant  dispersion error \cite{fleischmann19, hong20}.
Below, we show the numerical results of the Gaussian-Square-Triangle-Ellipse (GSTE) test, consisting of the linear advection equation $u_{t}+u_{x}=0$, under periodic boundary conditions on the domain $x\in\left[-1,1\right)$. The initial condition is represented by these four shapes, which are advected to the right:
\begin{gather}
   u(x,0) =
   \begin{cases}
      \frac{1}{6}[G(x,\beta,z-\delta)+4G(x,\beta,z)+G(x,\beta,z+\delta)], & x\in[-0.8,-0.6],\\
      1, & x\in[-0.4,-0.2],\\
      1-10(x-0.1), & x\in[0,0.2],\\
      \frac{1}{6}[F(x,\alpha,a-\delta)+4F(x,\alpha,a)+F(x,\alpha,a+\delta)], & x\in[0.4,0.6],\\
      0, & \text{otherwise}.
   \end{cases} \label{eq:linear advection} \\
   G(x,\beta,z) = e^{-\beta(x-z)^{2}}, \qquad 
   F(x,\alpha,a) = \sqrt{\max\{1-\alpha^{2}(x-a)^{2},\, 0\}},\nonumber
\end{gather}
with $z=-0.7$, $\delta=0.005$, $\beta=\dfrac{\log2}{36\delta^{2}}$, $a=0.5$, and $\alpha=10$.

\begin{figure}[htb]
    \centering
    \includegraphics[width=\textwidth]{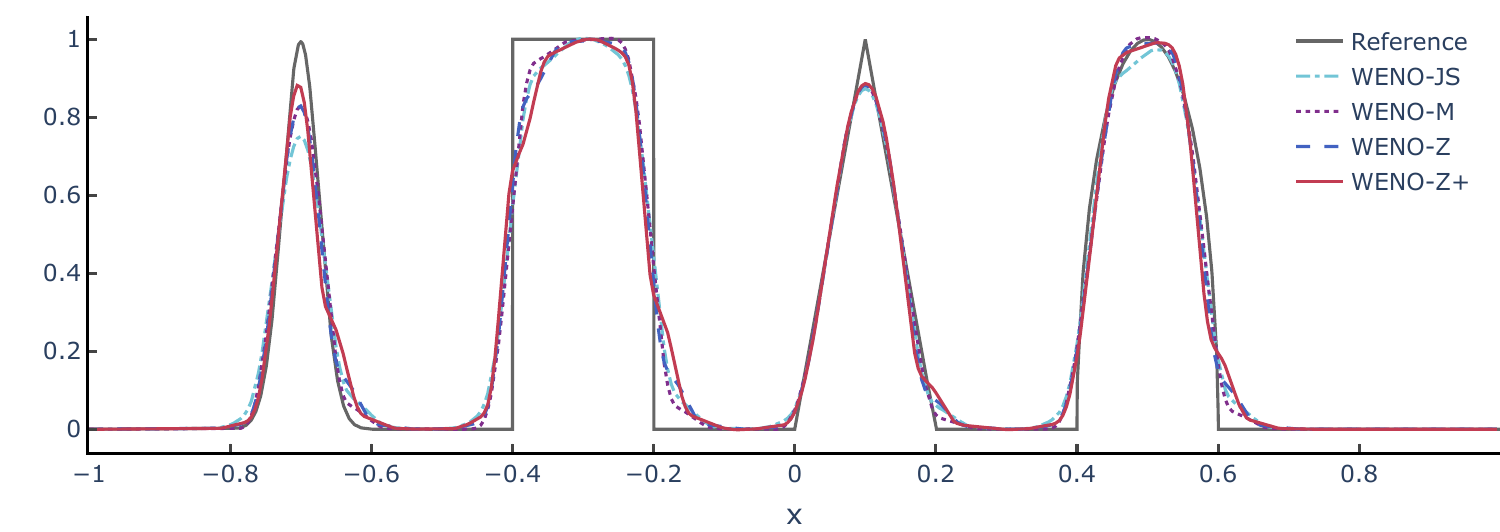} \\
    \centering
    \includegraphics[width=\textwidth]{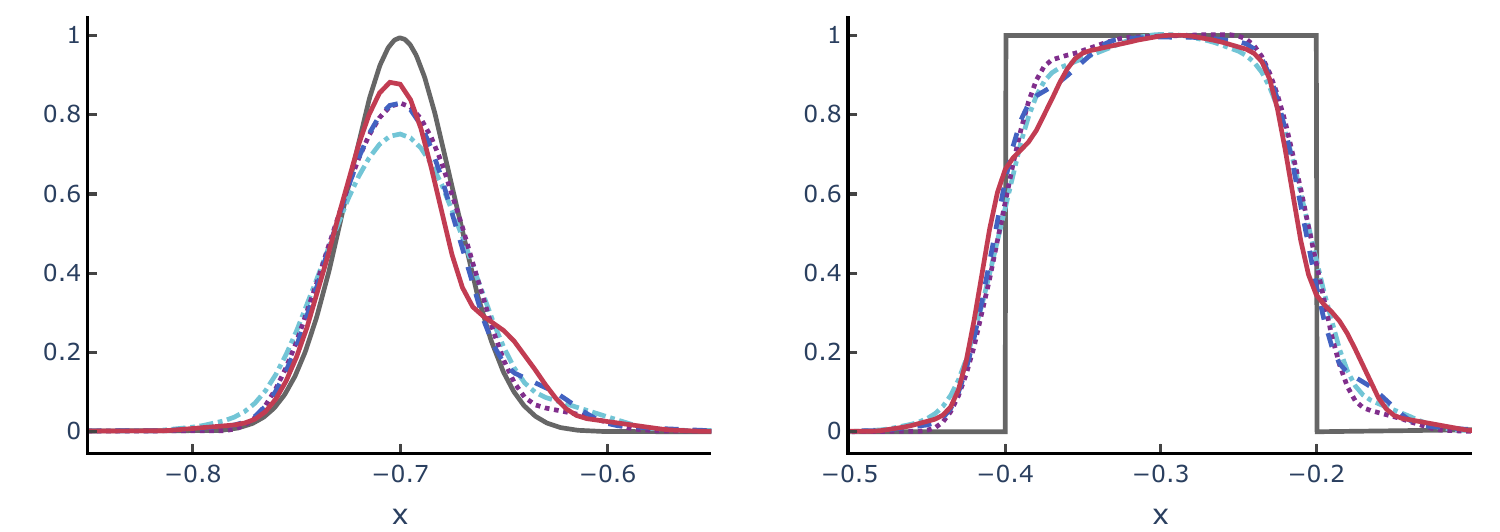}
    \caption{Numerical solutions of the Gaussian-Square-Triangle-Ellipse (GSTE) advection test at \(t=100\) with $N=400$ and CFL = $0.45$.}
    \label{fig:2a-1}
\end{figure}

It is well known that the solutions of WENO-JS, WENO-M, WENO-Z, and WENO-Z+ differ by many aspects, dissipation being the most clearly visible one. Yet, once the final time is increased to $t=100$, all numerical solutions suffer from dispersion error, as it is shown in Figure \ref{fig:intro 1a}, in the Introduction, and, below, in Figure \ref{fig:2a-1}.

\begin{figure}[htb]
     \centering
     \begin{subfigure}[b]{0.45\textwidth}
         \centering
         \includegraphics[width=\textwidth]{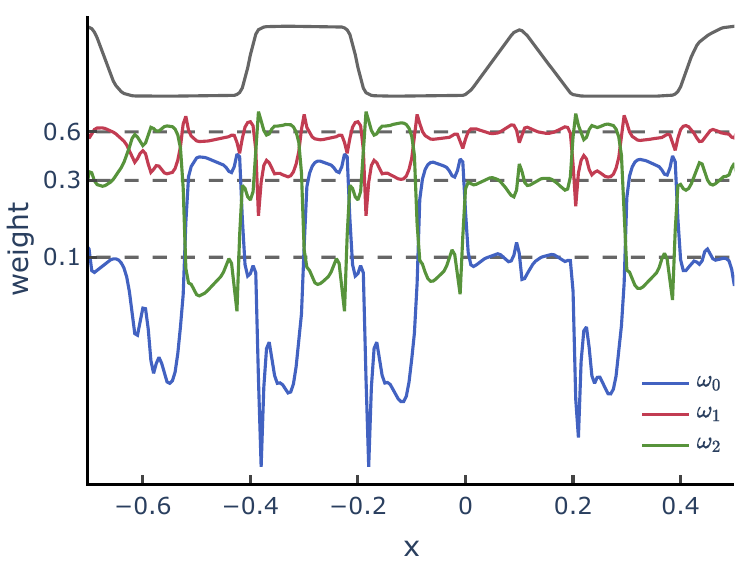}
         \caption{WENO-JS}
     \end{subfigure}
     \hfill
     \begin{subfigure}[b]{0.45\textwidth}
         \centering
         \includegraphics[width=\textwidth]{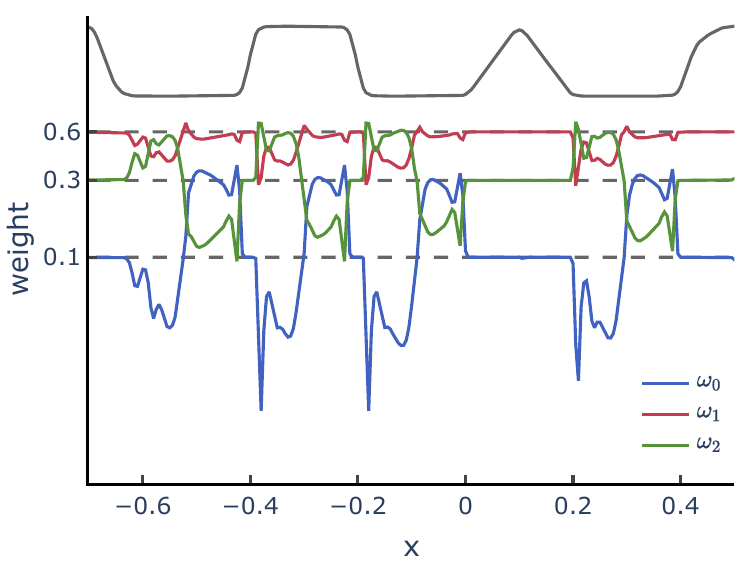}
         \caption{WENO-Z}
     \end{subfigure}
     \\
     \begin{subfigure}[b]{0.45\textwidth}
         \centering
         \includegraphics[width=\textwidth]{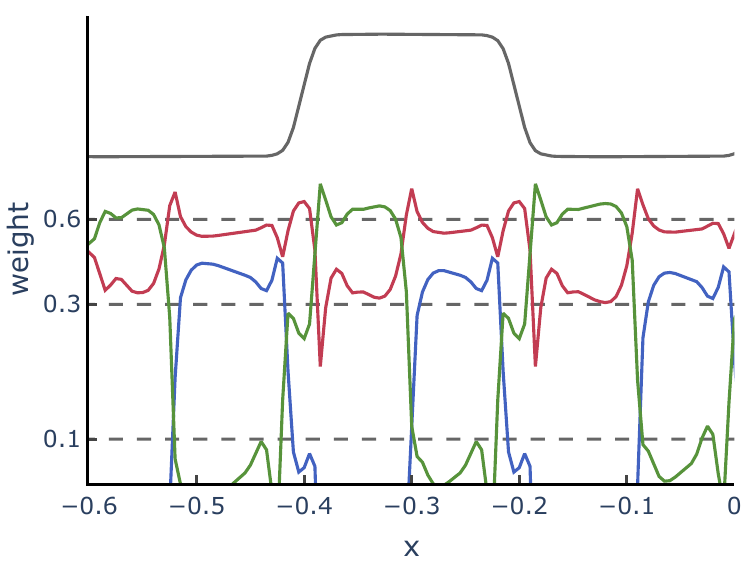}
         \caption{WENO-JS, zoom at the square wave}
     \end{subfigure}
     \hfill
     \begin{subfigure}[b]{0.45\textwidth}
         \centering
         \includegraphics[width=\textwidth]{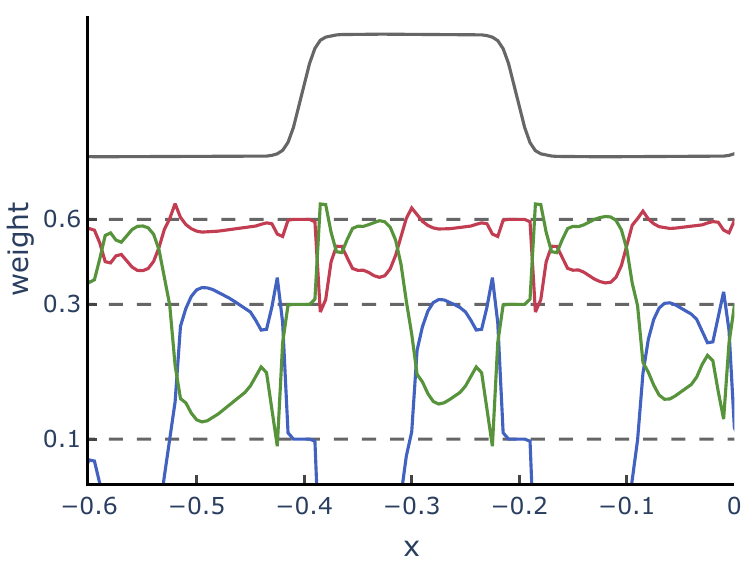}
         \caption{WENO-Z, zoom at the square wave}
     \end{subfigure}
    \caption{Weights \(\omega_k\) of WENO-JS and WENO-Z schemes for the GSTE numerical solution generated by WENO-Z) at \(t=2\) with $N=400$ and CFL = $0.45$.}
    \label{fig:GTSEweightsJSandZ}
\end{figure}

Due to their adaptive combination of  substencils, the numerical properties of WENO schemes may fit anywhere between those of the smoothest three-point substencil to those of the optimal five-point global stencil at every time step. For instance, in Figure \ref{fig:GTSEweightsJSandZ}, we show the weights response of WENO-JS and WENO-Z when facing the instantaneous solution of \eqref{eq:linear advection} generated by WENO-Z up to a final time \(t=2\). Some differences between the weights are clearly observed. 
As noted in \cite{borges08}, the WENO-Z weights of the discontinuous substencils are larger than the corresponding ones of WENO-JS. In general, the minimum value of each weight close to discontinuities are larger for WENO-Z than for WENO-JS, and to this fact has been attributed the better dissipative performance of WENO-Z over WENO-JS. 

Nevertheless, a very important aspect to the main discussion of this article is also present in Figure \ref{fig:GTSEweightsJSandZ}, when, for both schemes, the lateral weight \(\omega_{2}\) is larger than the central weight \(\omega_{1}\), not only at points of strong gradients, but also at smooth points of the solution, like \(x=-0.55\), \(-0.35\), and \(0.25\). This is due to the small (not visible to scale) fluctuations of the numerical solution, that forces the lateralization of the final convex combination, making it far from the original fifth-order central upwind one. Thus, it is natural that numerical dispersion may be present, even in odd-order WENO simulations.

It is also important to point out that on these smooth, although (essentially non)oscillatory, regions, the central weight of WENO-Z, $\omega_{1}^{Z}$, is overall much closer to the ideal value $d_{1}=0.6$ than $\omega_{1}^{JS}$, showing that the convex combination for WENO-Z is more centered than the one of WENO-JS. This  hints that when facing an unavoidable suboptimal convergence, since the problems WENO schemes deal with have shocks and discontinuities, one could relax the ties to the original fifth-order central-upwind scheme in favor of enhancing desirable numerical characteristics like small dispersion and dissipation errors, favoring a more centered combination of substencils.

%% file: 3-The-WENO-C.tex
\section{Centered WENO Schemes} \label{sec:weno-c}
In the first part of this section, we develop the main idea of the article: to augment the influence of the central substencil $S_{1}$ to decrease the departure from a centered upwind scheme due to an eventual overvaluation of the lateral substencils $S_{0}$ and $S_{2}$.  Later, we investigate the convergence of the centered schemes to be defined, showing that they achieve order $5$ at smooth parts of the solution. Nevertheless, we also show that the order decreases at critical points, and in the next section we provide a correction to recover order $5$ at such points.

\subsection{The WENO-C Scheme} 
We name this prototypical scheme WENO-C, and its weights are defined as the following, and incredibly simple, modification of the WENO-Z ones:
\begin{align*}
   \alpha_{k}^{C} & =d_{k}\Bigg[1+c_{k}\left(\frac{\tau}{\beta_{k}+\eps}\right)^{p}\Bigg]
   \quad \text{and} \quad 
   \omega_{k}^{C} = \frac{\alpha_{k}^{C}}{\sum_{j=0}^{2}\alpha_{j}^{C}}.
\end{align*}

\begin{figure}[H]
     \centering
     \begin{subfigure}[b]{0.45\textwidth}
         \centering
         \includegraphics[width=\textwidth]{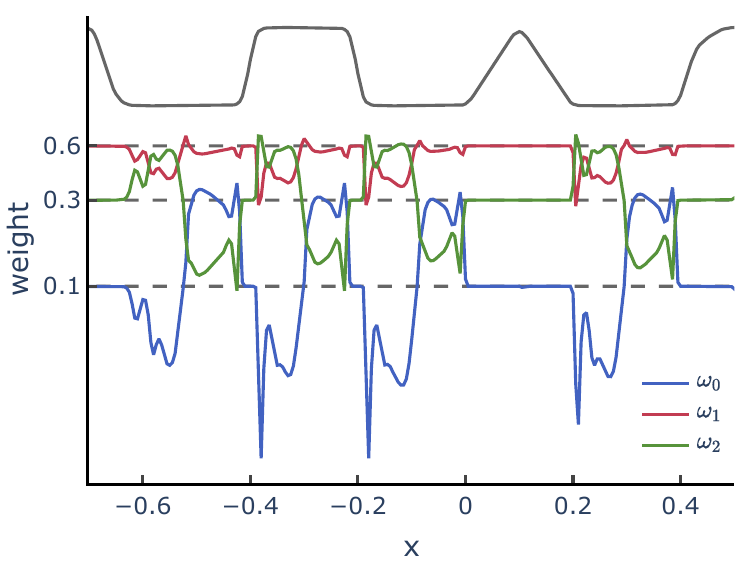}
         \caption{WENO-Z}
     \end{subfigure}
     \hfill
     \begin{subfigure}[b]{0.45\textwidth}
         \centering
         \includegraphics[width=\textwidth]{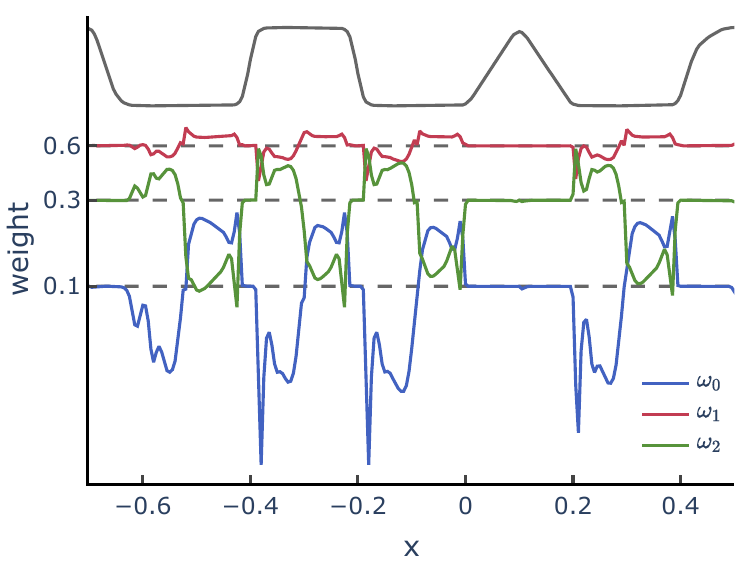}
         \caption{WENO-C}
     \end{subfigure}
     \\
     \begin{subfigure}[b]{0.45\textwidth}
         \centering
         \includegraphics[width=\textwidth]{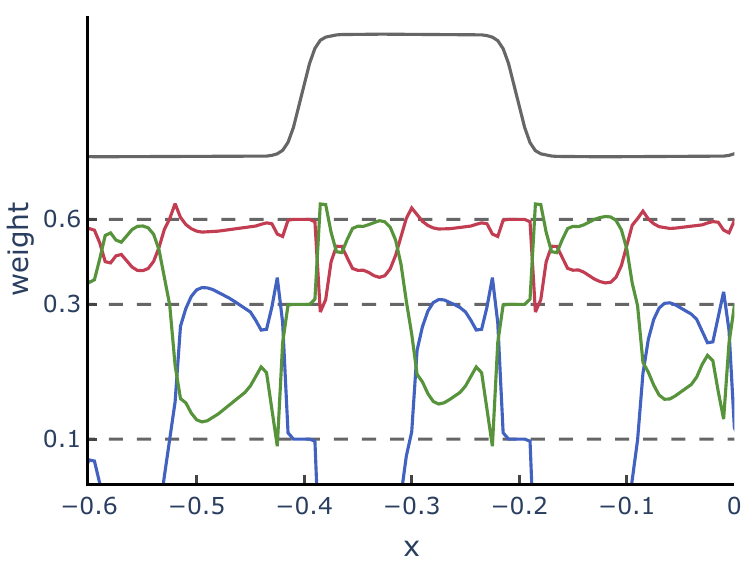}
         \caption{WENO-Z, zoom at the square wave}
     \end{subfigure}
     \hfill
     \begin{subfigure}[b]{0.45\textwidth}
         \centering
         \includegraphics[width=\textwidth]{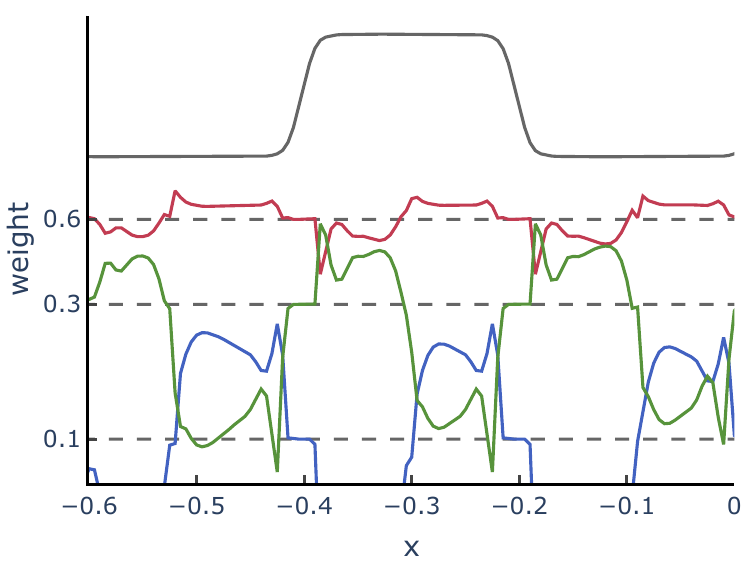}
         \caption{WENO-C, zoom at the square wave}
     \end{subfigure}
    \caption{Weights \(\omega_k\) of WENO-Z and WENO-C schemes for the GSTE advection problem at \(t=2\) with $N=400$ and CFL = $0.45$.}
    \label{fig:2a}
\end{figure}

We choose the coefficients $c_{k}$ to maintain the same relative balance of the nonlinear components of the weights in relation to the linear ones as in the WENO-Z weights. Thus, we require that 
\[
    \frac{c_{0}+c_{1}+c_{2}}{3}=1.
\]
It is also fundamental to perceive that, in the case of  nonsmoothness, information on the central WENO substencil $S_1$ is doubly penalized since $\beta_0$ and $\beta_2$, both measure roughness in complementary parts of $S_1$, as can be seen in Fig. \ref{fig:WENO-uniform-grid}. Thus, to correct this underestimation of the relative importance of the weight $\omega_{1}$ in relation to $\omega_{0}$ and $\omega_{2}$ we  impose that:
\[
   c_{1}=2c_{0}=2c_{2}
\]
 yielding
\begin{equation} \label{eq:ck-weno-c}
   (c_0, c_1, c_2) = \left(\frac{3}{4}, \frac{3}{2}, \frac{3}{4}\right).
\end{equation}
In the next section, we show there is sufficient leeway in the WENO analytical framework for such modifications, allowing for the strengthening of the role of the central substencil without losing numerical convergence. Before that, in Figure \ref{fig:2a}, we  perform the same comparison as of the last section, between the weights of WENO-Z and WENO-C when dealing with the numerical solution of \eqref{eq:linear advection} generated by WENO-Z up to \(t=2\). The WENO-C weights related to the central substencil are always larger than those of WENO-Z, i.e., $\omega_{1}^{C} > \omega_{1}^{Z}$ . Also, the dips in $\omega_{1}^{C}$ are clearly more contained than those of the $\omega_{1}^{Z}$ and, at smooth parts of the solution, $\omega_{1}^{C}$ is never smaller than the ideal weight $d_{1}$, nor it is  smaller than \(\omega_{2}^{C}\). The 
weights distribution of WENO-C tends to be more centered than the corresponding one of WENO-Z.
The result of this reinforcement of the central weight $\omega_{1}$ can be seen in Figure \ref{fig:intro 1a} through the fixing of the long-term advection error. It is remarkable that this  correction is not a feature of WENO-Z alone, for even WENO-JS can be ``centered'' if we simply multiply each $\alpha_{k}^{JS}$ by the same $c_{k}$ as above, obtaining:
\[
   \alpha_{k}^{JSC} 
      = c_{k} \alpha_{k}^{JS} = c_{k}\frac{d_{k}}{(\beta_{k} + \eps)^{p}},
\]

Again, the long-term advection distortions are eliminated,  confirming that when the solution contains discontinuities and/or high gradient, simply increasing the importance of the central Lagrangian interpolation $\hat{f^{1}}$ improves on the numerical properties of the convex combination \eqref{eq:convex combination}. Thus, it is not difficult to conclude that putting a floor to the central weight is the main factor to significantly decrease the dispersion error of the fifth-order WENO schemes.

For one more piece of evidence, Table \ref{tab:GSTE, weights Relative-error} shows the relative error \[e_{k}=\sum_{i=0}^N\dfrac{\left|\omega_{k}-d_{k}\right|}{d_{k}}\dx\] of the several WENO weights with respect to the ideal weights. The idea of this test is to show the behavior of WENO schemes in a typical situation, i.e., the case with a smoothed discontinuous numerical solution. For this, we obtain the solution of the GSTE problem by WENO-Z at \(t=2\) with \(N=400\) points and then we evaluate the weights of each WENO scheme in the whole domain of the solution, confirming that the centering imposed by the coefficients $c_{k}$ decreases the individual and global errors. 
\begin{table}[htb]
   \centering
   \begin{tabular}{lcccc}
   \toprule
   & $e_{0}$ & $e_{1}$ & $e_{2}$ & $e_{0} + e_{1} + e_{2}$ \\
   \midrule
   WENO-JS  & 2.29721 & 0.38900 & 1.25938 & 3.94561 \\
   WENO-Z   & 1.52174 & 0.25985 & 0.85565 & 2.63724 \\
   WENO-C   & 1.03149 & 0.17170 & 0.67824 & 1.88145 \\
   \bottomrule
   \end{tabular}
   \caption{Relative $L^1$  error of the weights computed for the WENO-Z solution of the GSTE at time $t=2$, $N=400$ and CFL = $0.45$.}
   \label{tab:GSTE, weights Relative-error}
\end{table}

\begin{figure}[htb]
    \centering
    \includegraphics[width=\textwidth]{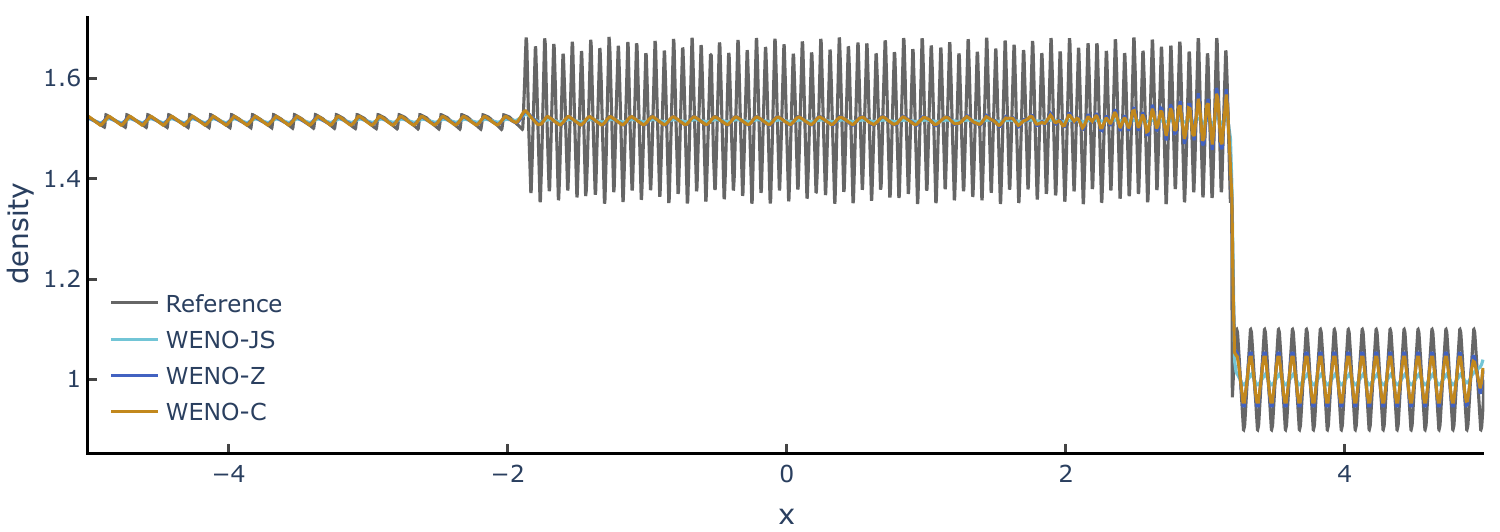} \\
    \centering
    \includegraphics[width=\textwidth]{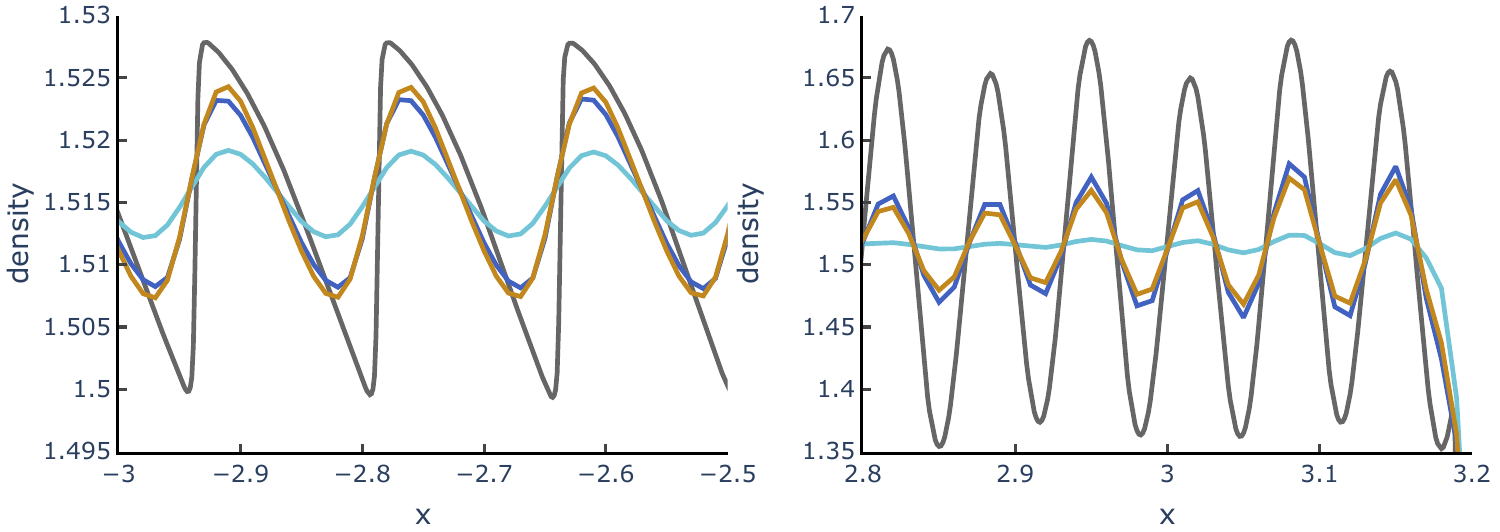}
    \caption{Numerical solutions of the Titarev--Toro shock-density wave problem \cite{titarev04} at \(t=5\) with $N=1000$ and CFL = $0.5$.}
    \label{fig:2d}
\end{figure}

Nevertheless, doubling the relative importance of the central weight is not enough to correct the scale jumps of the $\beta_k$ that occur in the vicinity of critical points, and WENO-C still has the same convergence deficiencies at these points as the other WENO schemes.
Let us exemplify this by comparing the results obtained by the WENO-C method when applied to the Titarev--Toro shock-density wave problem \cite{titarev04}:
\begin{equation}
   (\rho,u,p) = 
   \begin{cases}
      \left(1.515695,\, 0.523346,\, 1.805000\right), & x < -4.5, \\
      \left(1+\sin(20 \pi x)/10,\, 0,\, 1\right), & x \geq -4.5.
   \end{cases}
   \label{eq:Titarev-Toro}
\end{equation}
The Titarev--Toro problem is a variation of the classic Shu--Osher shock-density wave problem \cite{shu89}. It is traditionally used as an ideal one-dimensional flow configuration for testing numerical methods that must simultaneously capture shocks and avoid the damping of smooth turbulent structures. Its initial conditions consist of a perfect normal shock with imposed density fluctuations downstream, and as time progresses, these fluctuations interact with and become altered by the shock wave. Being more centered,
WENO-C better captures the shocklets of the solution, as seen in Figure \ref{fig:2d}. On the other hand, some amplitude is lost at critical points of the solution, requiring the investigation of the convergence of WENO-C at these points. In the next section, we thoroughly discuss this issue, arriving at a fix that results in a centered WENO scheme showing enhanced numerical properties when compared to WENO-Z.  

%% file: 3.1-Convergence-Analysis-WENO-C.tex
\subsection{Critical Points Convergence Analysis for WENO-C} \label{sec:convergence-weno-c}

We now analyze the convergence properties of WENO-C and show that it converges with optimal order at regular points of the solution but, indeed, loses convergence in the presence of critical points.
First, if \(f(x)\) does not contain a critical point, then Eqs. \eqref{eq:betataylor}--\eqref{eq:tautaylor} 
give:
\begin{align*}
   \alpha_{k}^{C(\pm)}
      &= d_{k}\left[1+c_{k}\left(\frac{\left|\left(\dfrac{13}{3}f_{i}''f_{i}''' - f_i'f_i^{(4)}\right)\Delta x^{5}+\Ord(\Delta x^{6})\right|}{\left(f_{i}'\right)^{2}\Delta x^{2}+\Ord(\Delta x^{3})}\right)^{p}\right]\\
      &= d_{k} + \Ord(\dx^{3p}).
\end{align*}
This implies \(\omega_{k}^{C(\pm)} = d_{k} + \Ord(\dx^{3p})\) and, therefore, by Condition \ref{cond:sufficient}, WENO-C has optimal order 5 if \(p \geq 1\).

Now, consider \(x_i\) a simple critical point of order \(\ncp=1\) of \(f\). In this case, Eqs. \eqref{eq:betataylor}--\eqref{eq:tautaylor} now give:
\begin{align*}
   \alpha_{k}^{C(+)}
      &= d_{k}\left[1+c_{k}\left(\frac{\left|\dfrac{13}{3}f_{i}''f_{i}'''\Delta x^{5}+\Ord(\Delta x^{7})\right|}{\dfrac{13}{12}\left(f_{i}''\right)^{2}\Delta x^{4}+\Ord(\Delta x^{5})}\right)^{p}\right]\\
      &= d_{k}\left[1+4^{p}c_{k}\left(\left|\frac{f_{i}'''}{f_{i}''}\right|\Delta x+\Ord(\Delta x^{2})\right)^{p}\right]
\end{align*}
and
\begin{align*}
   \alpha_{k}^{C(-)}
      &= d_{k}\left[1+c_{k}\left(\dfrac{\left|-\dfrac{13}{3}f_{i}''f_{i}'''\Delta x^{5}+\Ord(\Delta x^{7})\right|}{\dfrac{25}{12}\left(f_{i}''\right)^{2}\Delta x^{4}+\Ord(\Delta x^{5})}\right)^{p}\right] \\
      &= d_{k}\left[1+\left(\frac{52}{25}\right)^{p}c_{k}\left(\left|\frac{f_{i}'''}{f_{i}''}\right|\Delta x+\Ord(\Delta x^{2})\right)^{p}\right]
\end{align*}

resulting in
\[
   \alpha_{k}^{C(\pm)} = d_k + \Ord(\dx^p).
\]
 
We let $p=2$  and focus on the requirements of Condition \ref{cond:necessary}. For that, let us consider \[\overline{c}=\sum_{k=0}^{2}c_{k}d_{k}=\dfrac{c_{0}+6c_{1}+3c_{2}}{10}\] and check the order of convergence of the weights of WENO-C to the ideal weights $d_k$:
\begin{align*}
   \omega_{k}^{C(+)}
      &= \frac{d_{k}\left[1+16c_{k}\left(\dfrac{f_{i}'''}{f_{i}''}\right)^{2}\Delta x^{2}+\Ord(\Delta x^{3})\right]}{d_{k}\left[1+16\overline{c}\left(\dfrac{f_{i}'''}{f_{i}''}\right)^{2}\Delta x^{2}+\Ord(\Delta x^{3})\right]}\\
      &= d_{k}+16d_{k}\left(c_{k}-\overline{c}\right)\left(\frac{f_{i}'''}{f_{i}''}\right)^{2}\Delta x^{2}+\Ord(\Delta x^{3}),
\end{align*}
and analogously we have
\[
   \omega_{k}^{C(-)} = d_{k}+\frac{2704}{625}d_{k}\left(c_{k}-\overline{c}\right)\left(\frac{f_{i}'''}{f_{i}''}\right)^{2}\Delta x^{2}+\Ord(\Delta x^{3}).
\]
so Eq. \eqref{eq:idealweightscond} is satisfied for the WENO-C weights, but not necessarily Eq. \eqref{eq:Aks}. We could fix it by taking $c_0 = c_1 = c_2 = \overline{c}$; nevertheless, this leads to the uninteresting case where the relative importance of the stencils would not be changed. Alternatively, plugging the expansions of $\omega_{k}^C$ above into \eqref{eq:Aks}, we obtain the sufficient condition:
\begin{gather*}
   48d_{0}(c_{0}-\overline{c}) - \frac{8112}{625}d_{0}(c_{0}-\overline{c}) - 16d_{1}(c_{1}-\overline{c}) + {} \\
   \quad {} + \frac{2704}{625}d_{1}(c_{1}-\overline{c}) + 16d_{2}(c_{2}-\overline{c}) - \frac{2704}{625}d_{2}(c_{2}-\overline{c})=0
\end{gather*}
and this amounts to the simple relation
\begin{equation}
   4c_{0} - 11c_{1} + 7c_{2} = 0. \label{eq:relationcks}
\end{equation}
For instance, if we choose $\vec{c} = (c_0, c_1, c_2) =\left(\dfrac{7}{4},\dfrac{14}{11},1\right)$, we recover fifth-order convergence at critical points. Still, this one is definitely not a ``centered'' convex combination (since \(c_1 < c_0\)) and, as expected, the long-term advection error as discussed in Section \ref{sec:dispersion-error} is not corrected. 

Another idea is to choose from the vectors $\vec{c}$ satisfying \eqref{eq:relationcks} the closest in the least squares sense to \((3/4, 3/2, 3/4)\) (see Eq. \eqref{eq:ck-weno-c}). The solution is $\vec{c} \approx (0.463,0.508,0.530)$, which again does not ``center'' the convex combination and does not generate good nondispersive results on the long-term advection GSTE problem.

In sum, it is not possible to make WENO-C with \(p=2\) to achieve optimal order at critical points and be centered at the same time. In the next section, we will keep the initial chosen values of $c_k$, maintaining its centered structure, and search for a term that, close to critical points, is small enough to accelerate the convergence to zero of the nonlinear parts of the WENO-C weights.

%% file: 4-WENOZD-and-ZDplus.tex
\section{The WENO-ZC and WENO-ZC+ Schemes} \label{sec:weno-zc}

The search for an appropriate convergence accelerator for WENO-C was deeply influenced by a parallel investigation related to WENO-Z+, which we describe in the following. The main goal was to improve on the anti-dissipative term 
\[
   \xi_{k}^{Z+} = \lambda \frac{\beta_k}{\tau + \eps}
\]
which is known to overamplify curvature features of the numerical solution, as can be seen in Figure \ref{fig:4c} (see also \cite{luo21a} and \cite{luo21b}). The idea was to define another anti-dissipative term, one that was limited above, and for that, we adopted the term
\( \dfrac{\beta_k}{\tau + \overline{\beta} + \eps} \), with $\overline{\beta}=\frac{1}{3}\sum_{j=0}^{2}\beta_{j}$.
Now there is no grid-size dependency, and the addition of $\overline{\beta}$ in the denominator excludes the overamplification caused by a small $\tau$:
\[
   \frac{\beta_{k}}{\tau+\overline{\beta}+\eps}
      = \dfrac{\beta_{k}}{\tau+\frac{1}{3}\sum_{j=0}^{2}\beta_{j}+\eps}
      = \dfrac{3\beta_{k}}{3\tau+\beta_{0}+\beta_{1}+\beta_{2}+3\eps}
      \leq \dfrac{3\beta_{k}}{\beta_{k}}
      = 3.
\]

It is worth noting that, at first glance, it would be difficult to evaluate how the new anti-dissipative term would interact with the other components of the unnormalized weights. However, due to the relation 
\begin{equation*}
     \frac{1}{3}\sum_{k=0}^{2} \frac{\beta_{k}}{\tau+\overline{\beta}} 
    = \frac{1}{\tau + \overline{\beta}} \left(\frac{1}{3}\sum_{k=0}^{2}\beta_{k}\right)
    = \frac{\overline{\beta}}{\tau + \overline{\beta}} 
    = \frac{\overline{\beta} + \tau - \tau}{\tau + \overline{\beta}} 
    = 1 - \frac{\tau}{\tau + \overline{\beta}} 
    = 1 - \zeta,
\end{equation*}
we know $\zeta$ is small whenever the influence of the new anti-dissipative terms is large, and vice-versa. The relation above is the original motivation to choose $\zeta$ as the convergence accelerator to the nonlinear components of the WENO-C weights.
Thus, the new WENO schemes, respectively named as WENO-ZC and WENO-ZC+, are respectively defined as
\begin{equation}
   \alpha_{k}^{ZC} = d_{k}\Bigg[1+c_{k}\left(\frac{\tau}{\beta_{k}+\eps}\right)^{p}\left(\frac{\tau}{\tau+\overline{\beta}+ \eps}\right)^{p}\Bigg],
   \quad
   \omega_{k}^{ZC} = \dfrac{\alpha_{k}^{ZC}}{\sum_{j=0}^2 \alpha_{j}^{ZC}}, \label{eq:weno-zc-weights}
\end{equation}
and
\begin{align}
   \alpha_{k}^{ZC+}
      &= d_{k} \Bigg[1 + c_{k}\left(\frac{\tau}{\beta_{k}+\eps}\right)^{p} \left(\frac{\tau}{\tau+\overline{\beta}+\varepsilon}\right)^{p} + \frac{\beta_{k}}{\tau+\overline{\beta}+ \eps}\Bigg], \nonumber \\
   \omega_{k}^{ZC+} &= \frac{\alpha_{k}^{ZC+}}{\sum_{j=0}^{2}\alpha_{j}^{ZC+}},\label{eq:weno-zcplus-weights}
    \xi_{k}^{Z+} = \frac{\beta_{k}}{\tau+\overline{\beta}+ \eps}.
\end{align}
We will first analyze the characteristics of the WENO-ZC method, showing that it represents an enhancement of WENO-Z. After that, we will pass to investigate the numerical properties of WENO-ZC+.

\begin{figure}[htb]
    \centering
    \includegraphics[width=\textwidth]{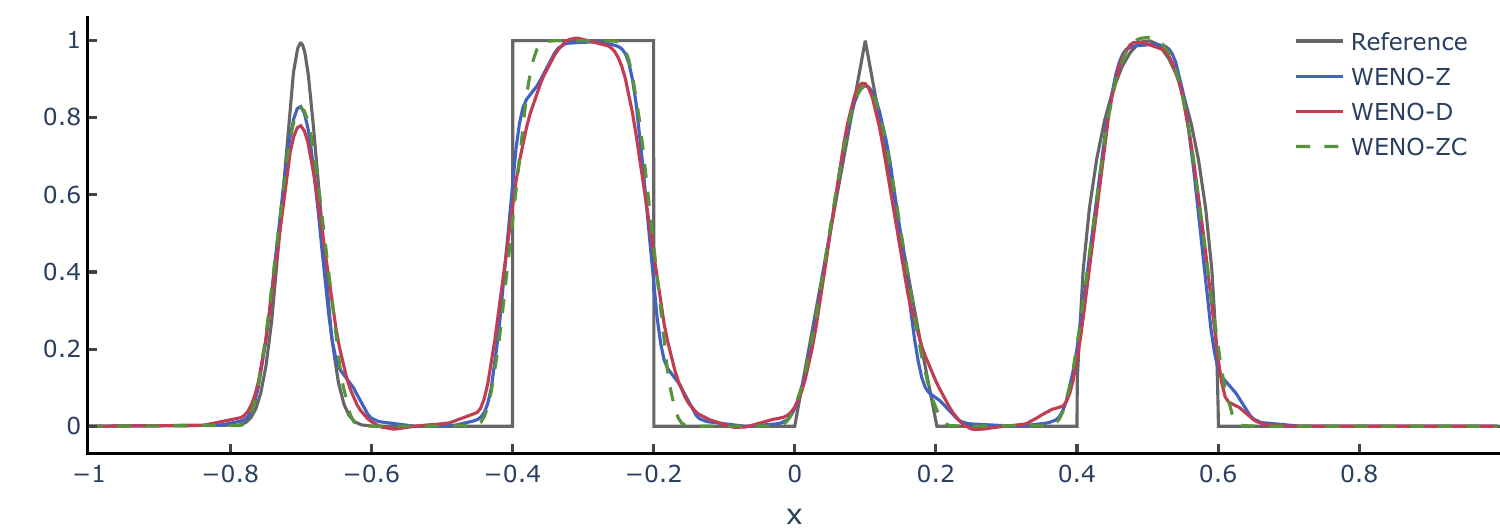} \\
    \includegraphics[width=\textwidth]{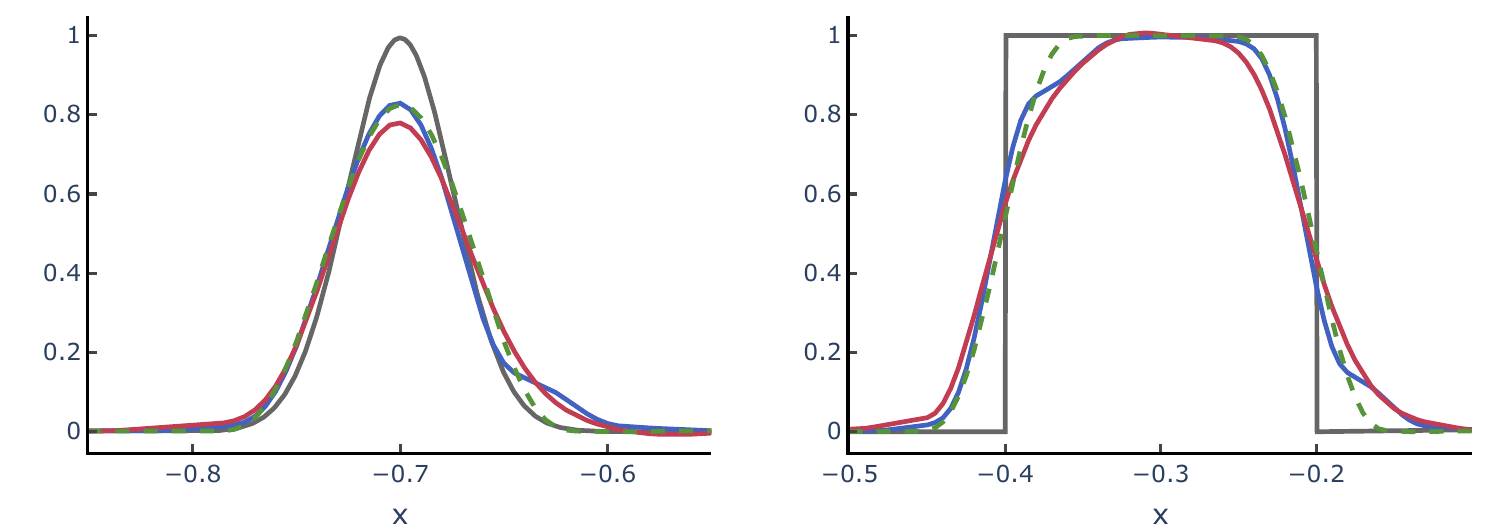}
    \caption{Numerical solutions of the Gaussian-Square-Triangle-Ellipse (GSTE) advection test at \(t=100\) with $N=400$ and CFL = $0.45$.}
    \label{fig:2b}
\end{figure}

\begin{figure}[htb]
    \centering
    \includegraphics[width=\textwidth]{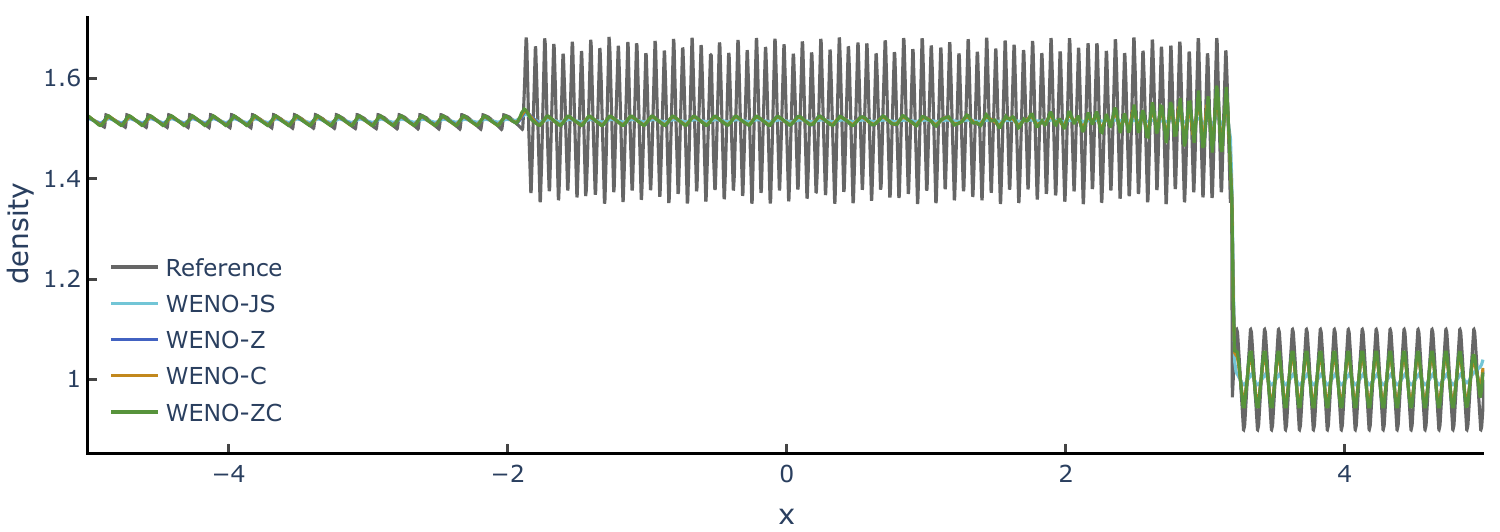} \\
    \includegraphics[width=\textwidth]{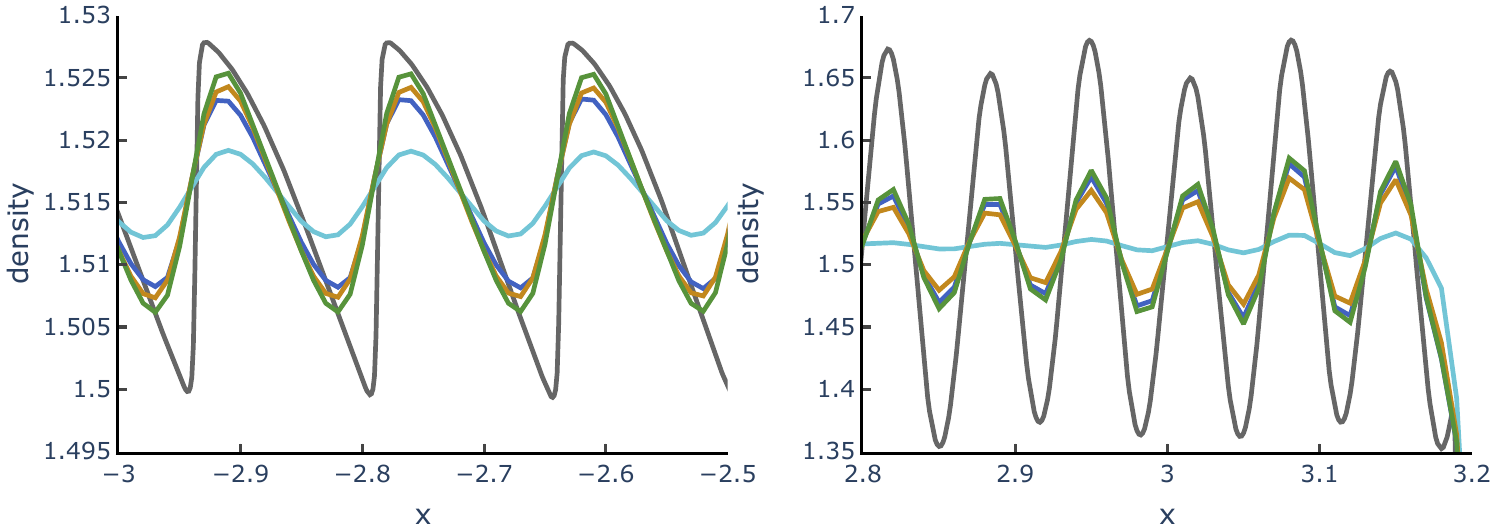}
    \caption{Numerical solutions of the Titarev--Toro shock-density wave problem \cite{titarev04} at \(t=5\) with $N=1000$ and CFL = $0.5$.}
    \label{fig:3e}
\end{figure}

We start by showing that WENO-ZC maintains the good dispersive characteristics of WENO-C in the GTSE problem, see Figure \ref{fig:2b}; which is also a good place to show that WENO-D suffers from the same long-term dispersion error as WENO-Z.
Figure \ref{fig:3e} shows the numerical results of WENO-ZC in the Titarev--Toro problem. One can see that not only the convergence of the amplitude at critical points has been recovered but also that WENO-ZC shows less dissipation than WENO-Z in all regions of the numerical solution, particularly at the shocklets.

\subsection{Critical Points Convergence Analysis of WENO-ZC} \label{sec:convergence-analysis-weno-zc}
To understand how the convergence-accelerating term works, let us analyze its Taylor series on critical points of orders $0$, $1$, and $2$. First, in the absence of critical points, and neglecting \(\eps\), Eqs. \eqref{eq:betataylor}--\eqref{eq:tautaylor} (see \ref{sec:appendix-taylor}) give
\begin{align*}
   \left(\dfrac{\tau^{\pm}}{\tau^{\pm} + \overline{\beta}^{\pm}}\right)^p
      &= \left(\dfrac{\left|\dfrac{13}{3}f''_{i}f'''_{i} - f_{i}'f_{i}^{(4)}\right|\Delta x^{5} + \Ord(\Delta x^6)}{\left|\dfrac{13}{3}f''_{i}f'''_{i} - f_{i}'f_{i}^{(4)}\right|\Delta x^{5} + \big(f_{i}'\big)^{2}\Delta x^{2} + \Ord(\dx^3)}\right)^p \\
      &= \left(\left|\frac{13}{3}\frac{f''_{i}f'''_{i}}{\big(f_{i}'\big)^{2}} - \frac{f_{i}'f_{i}^{(4)}}{\big(f_{i}'\big)^{2}}\right| \dx^3 + \Ord(\dx^4)\right)^p
      = \Ord(\dx^{3p}).
\end{align*}
This, coupled with the fact that \(\left(\dfrac{\tau}{\beta_{k}+\eps}\right)^{p} = \Ord(\dx^{3p})\) \cite{don13}, gives
\[
   \alpha_k^{ZC\pm} = d_k + \Ord(\dx^{6p})
   \quad \therefore \quad
   \omega_k^{ZC\pm} = d_k + \Ord(\dx^{6p}),
\]
and, therefore, the ZC weights satisfy Condition \ref{cond:sufficient} for \(p \geq 1/2\). Note that this result is entirely independent of the values of the coefficients \(c_k\).

Now, if \(x_i\) is a critical point of order \(\ncp = 1\), we have
\begin{align*}
   \left(\dfrac{\tau^{+}}{\tau^{+} + \overline{\beta}^{+}}\right)^p
      &= \left(\dfrac{\left|\dfrac{13}{3}f''_{i}f'''_{i}\right|\Delta x^{5} + \Ord(\Delta x^6)}{\left|\dfrac{13}{3}f''_{i}f'''_{i}\right|\Delta x^{5} + \dfrac{13}{12}\big(f_{i}''\big)^{2}\Delta x^{4} + \Ord(\dx^5)}\right)^p \\
      &= \left(4 \left|\frac{f'''_{i}}{f_{i}''}\right| \dx + \Ord(\dx^2)\right)^p
      = \Ord(\dx^{p}), \\
   \left(\dfrac{\tau^{-}}{\tau^{-} + \overline{\beta}^{-}}\right)^p
      &= \left(\dfrac{\left|\dfrac{13}{3}f''_{i}f'''_{i}\right|\Delta x^{5} + \Ord(\Delta x^6)}{\left|\dfrac{13}{3}f''_{i}f'''_{i}\right|\Delta x^{5} + \dfrac{25}{12}\big(f_{i}''\big)^{2}\Delta x^{4} + \Ord(\dx^5)}\right)^p \\
      &= \left(\dfrac{52}{25} \left|\frac{f'''_{i}}{f_{i}''}\right| \dx + \Ord(\dx^2)\right)^p
      = \Ord(\dx^{p}).
\end{align*}
This, together with the fact that \(\left(\dfrac{\tau}{\beta_{k}+\eps}\right)^{p} = \Ord(\dx^{p})\) \cite{don13}, yields
\[
   \alpha_k^{ZC\pm} = d_k + \Ord(\dx^{2p})
   \quad \therefore \quad
   \omega_k^{ZC\pm} = d_k + \Ord(\dx^{2p}),
\]
and, therefore, the ZC weights satisfy Condition \ref{cond:sufficient} for \(p \geq 3/2\)---a result that, again, does not depend on the \(c_k\).

For critical points of order \(\ncp=2\), the situation changes. We now have
\begin{align*}
   \left(\dfrac{\tau^{+}}{\tau^{+} + \overline{\beta}^{+}}\right)^p
      &= \left(\dfrac{\left|\dfrac{103}{36}f'''_{i}f^{(4)}_{i}\right|\Delta x^{7} + \Ord(\Delta x^9)}{\left|\dfrac{103}{36}f'''_{i}f^{(4)}_{i}\right|\Delta x^{7} + \dfrac{29}{36}\big(f_{i}'''\big)^{2}\Delta x^{6} + \Ord(\dx^7)}\right)^p \\
      &= \left(\dfrac{103}{29} \left|\frac{f^{(4)}_{i}}{f_{i}'''}\right| \dx + \Ord(\dx^2)\right)^p
      = \Ord(\dx^{p}), \\
   \left(\dfrac{\tau^{-}}{\tau^{-} + \overline{\beta}^{-}}\right)^p
      &= \left(\dfrac{\left|\dfrac{13}{3}\big(f'''_{i}\big)^2\right|\Delta x^{6} + \Ord(\Delta x^7)}{\left|\dfrac{13}{3}\big(f'''_{i}\big)^2\right|\Delta x^{6} + \dfrac{71}{36}\big(f_{i}'''\big)^{2}\Delta x^{6} + \Ord(\dx^7)}\right)^p \\
      &= \left(\dfrac{156}{227} + \Ord(\dx)\right)^p
      = \Ord(1).
\end{align*}
Since \(\left(\dfrac{\tau^-}{\beta^-_{k}+\eps}\right)^{p} = \Ord(1)\) \cite{don13}, this gives
\[
   \alpha_k^{ZC-} = d_k + \Ord(1)
   \quad \therefore \quad
   \omega_k^{ZC-} = d_k + \Ord(1),
\]
and Condition \ref{cond:necessary} cannot be satisfied for any value of $p$.

The comparison of the WENO discrete derivative with the actual derivative \(f'\) of a given function \(f(x)\) at the grid points \(x_i\), \(i = 0, \ldots, N\), with increasing grid sizes \(N\), is shown in Table \ref{tab:weno-ZC-ncp0} and confirm the theoretical results above.  We use
\begin{equation} \label{eq:f0}
    f_0(x) = \exp\!\left(x - \frac{\sin(\pi x)}{2\pi}\right),
    \quad
    x \in [-1, 1].
\end{equation}
for a function with no critical points, and 
\begin{equation} \label{eq:f1}
    f_1(x) = \sin\!\left(\pi x - \frac{\sin(\pi x)}{\pi}\right),
    \quad
    x \in [-1, 1].
\end{equation}
for a function with critical points of order \(\ncp=1\), with the important property that the third derivative does not vanish at the same location (see \cite{henrick05}).
Finally, for a function with a critical point of order \(\ncp = 2\) (at \(x=1/2\)), we choose
\begin{equation} \label{eq:f2}
    f_2(x) = \sin\!\left(\pi x + \cos(\pi x) + \sin(\pi x) + \frac{\cos^2(\pi x)}{2} + \cos^3(\pi x)\right),
    \quad
    x \in [-1, 1].
\end{equation}

\begin{table}[H]
\centering
\begin{tabular}{lllllll}
\toprule
 & \multicolumn{2}{l}{\(f_0(x)\) (no critical points)} & \multicolumn{2}{l}{\(f_1(x)\) (\(\ncp = 1\))} & \multicolumn{2}{l}{\(f_2(x)\) (\(\ncp = 2\))} \\
\cmidrule(lr){2-3} \cmidrule(lr){4-5} \cmidrule(lr){6-7}
\(1/\Delta x\) & \(L^1\) error & \(L^1\) order & \(L^1\) error & \(L^1\) order & \(L^1\) error & \(L^1\) order \\
\midrule
\(25\) & \(2.76205 \times 10^{-5}\) & --- & \(8.31844 \times 10^{-4}\) & --- & \(2.53246 \times 10^{-1}\) & --- \\
\(50\) & \(8.83108 \times 10^{-7}\) & \(4.96701\) & \(2.70148 \times 10^{-5}\) & \(4.94449\) & \(1.23091 \times 10^{-2}\) & \(4.36274\) \\
\(100\) & \(2.76013 \times 10^{-8}\) & \(4.99978\) & \(7.99497 \times 10^{-7}\) & \(5.07851\) & \(1.01371 \times 10^{-3}\) & \(3.60200\) \\
\(200\) & \(8.60551 \times 10^{-10}\) & \(5.00333\) & \(2.41364 \times 10^{-8}\) & \(5.04981\) & \(9.54303 \times 10^{-5}\) & \(3.40906\) \\
\(400\) & \(2.68545 \times 10^{-11}\) & \(5.00202\) & \(7.47436 \times 10^{-10}\) & \(5.01312\) & \(1.10972 \times 10^{-5}\) & \(3.10425\) \\
\(800\) & \(8.43036 \times 10^{-13}\) & \(4.99343\) & \(2.33412 \times 10^{-11}\) & \(5.00100\) & \(1.36255 \times 10^{-6}\) & \(3.02581\) \\
\bottomrule
\end{tabular}
\caption{Numerical results of the accuracy test for the WENO-ZC scheme. The expressions for the test functions \(f_0(x)\), \(f_1(x)\), and \(f_2(x)\) are given in Eqs. \eqref{eq:f0}--\eqref{eq:f2}.}
\label{tab:weno-ZC-ncp0}
\end{table}

The accelerating factor \(\left(\dfrac{\tau}{\tau + \overline{\beta}}\right)^p\) fixes WENO-C's issue of loss of accuracy at critical points of order \(\ncp=1\), but not for \(\ncp=2\). Nevertheless, for second-order critical points, both WENO-M and WENO-Z are equally not able to achieve such convergence, even though their good performance at problems with shocks is well established.

\begin{figure}[H]
     \centering
     \begin{subfigure}[b]{0.45\textwidth}
         \centering
         \includegraphics[width=\textwidth]{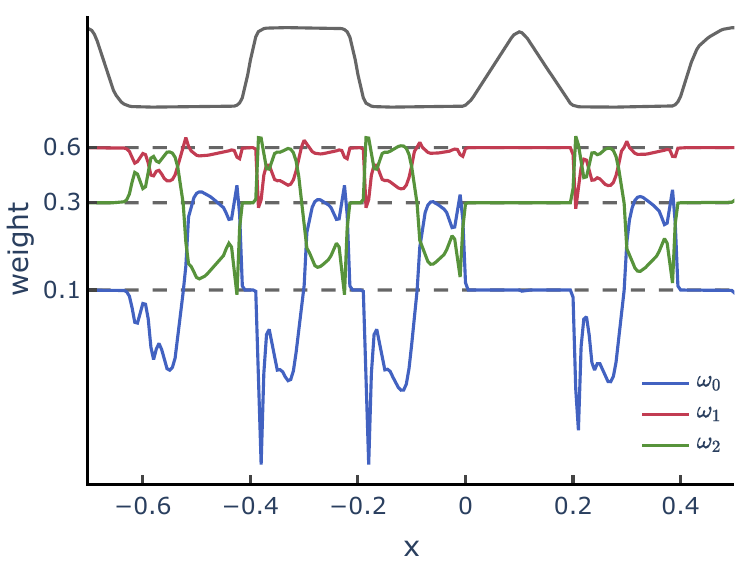}
         \caption{WENO-Z}
     \end{subfigure}
     \hfill
     \begin{subfigure}[b]{0.45\textwidth}
         \centering
         \includegraphics[width=\textwidth]{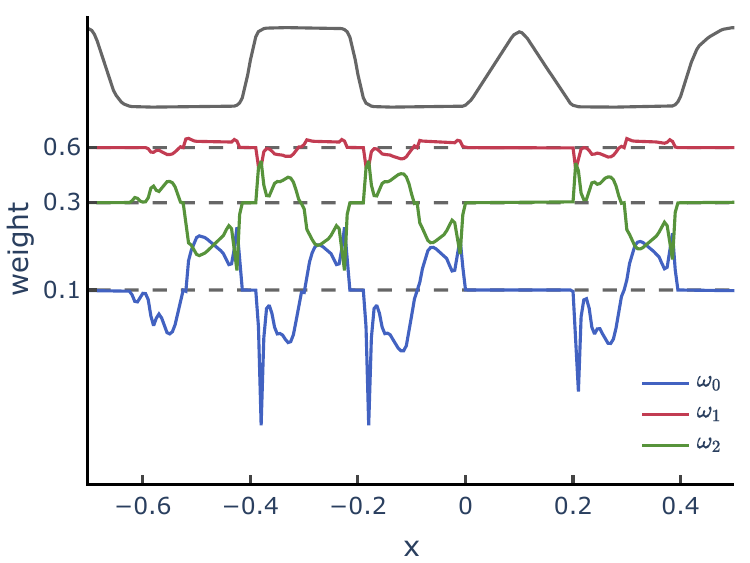}
         \caption{WENO-ZC}
     \end{subfigure}
     \\
     \begin{subfigure}[b]{0.45\textwidth}
         \centering
         \includegraphics[width=\textwidth]{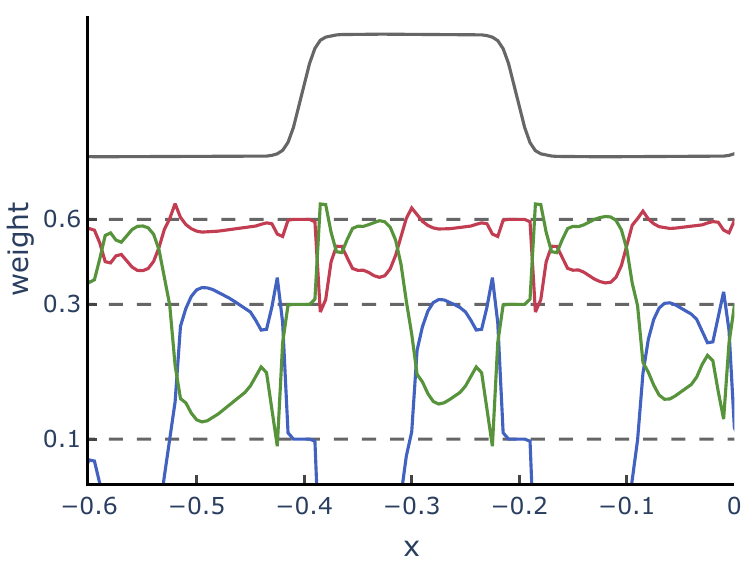}
         \caption{WENO-Z, zoom}
     \end{subfigure}
     \hfill
     \begin{subfigure}[b]{0.45\textwidth}
         \centering
         \includegraphics[width=\textwidth]{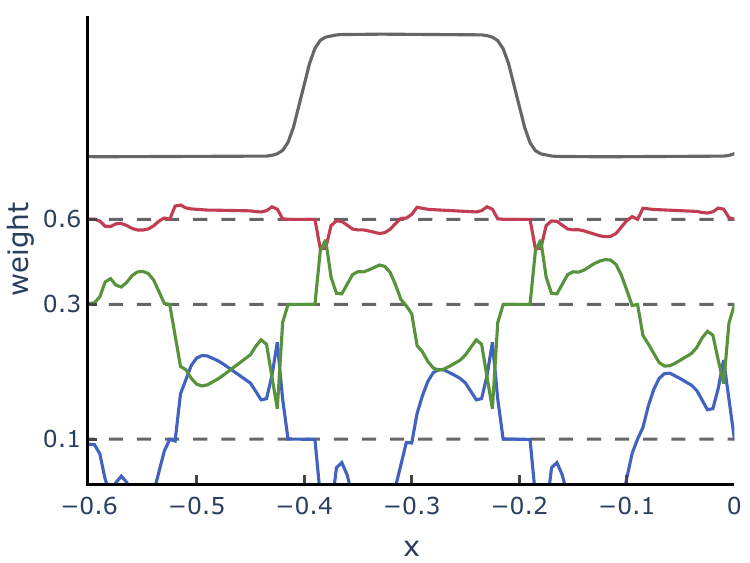}
         \caption{WENO-ZC, zoom}
     \end{subfigure}
    \caption{Weights \(\omega_k\) of WENO-Z and WENO-ZC schemes for the GSTE advection problem at \(t=2\) with $N=400$ and CFL = $0.45$.}
    \label{fig:3g}
\end{figure}

Figure \ref{fig:3g} shows the normalized weights of WENO-Z and WENO-ZC for the GSTE advection problem at \(t=2\). The central weight \(\omega_1^{ZC}\) is never smaller than \(\omega_2^{ZC}\) at smooth parts of the solution and is much closer to the ideal weight $d_1$ than in the WENO-C case. This indicates that the essentially nonoscillatory property is now being enforced in a stronger way.

%% file: 4.1-WENOZDplus.tex
\subsection{The WENO-ZC+ Scheme} \label{sec:weno-zc-plus}

The WENO-ZC+ scheme is 
intended to be an analogous extension to WENO-ZC as WENO-Z+ was to WENO-Z. We now pass to analyze its numerical properties related to the new anti-dissipative term, $ \xi_{k}^{ZC+}$. The first thing to notice is that in Eq. \eqref{eq:weno-zcplus-weights}, the $c_k$ do not multiply $ \xi_{k}^{ZC+}$ because doing so the resulting scheme would not be of order $5$. Thus, the changing of the balance between the linear and nonlinear parts of the weights decreased the relevance of the central substencil, demanding the increasing of the $c_{k}$ to account for the addition of $ \xi_{k}^{ZC+}$ as
\[
    \frac{c_{0}+c_{1}+c_{2}}{3}=\frac32
\]
to recover the nondispersive results of WENO-ZC. This yielded the following new values for the coefficients \(c_k\):
\[
    (c_0, c_1, c_2) = \left(\frac{9}{8}, \frac{9}{4}, \frac{9}{8}\right).
\]

\begin{figure}
    \centering
    \includegraphics[width=\textwidth]{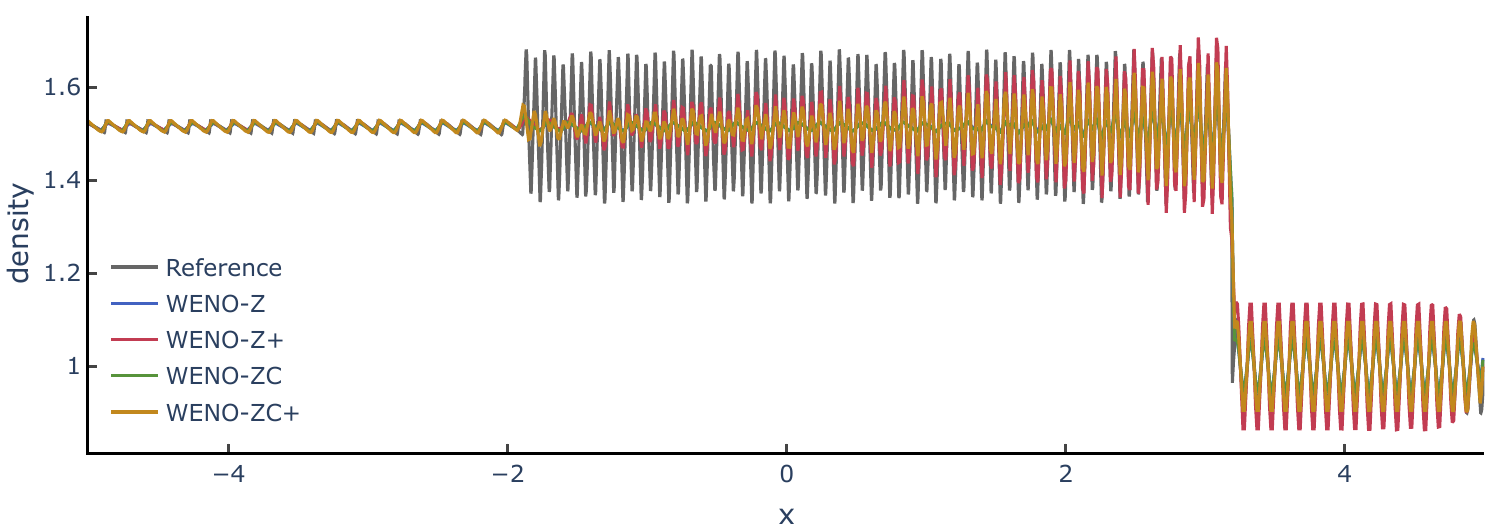} \\
    \includegraphics[width=\textwidth]{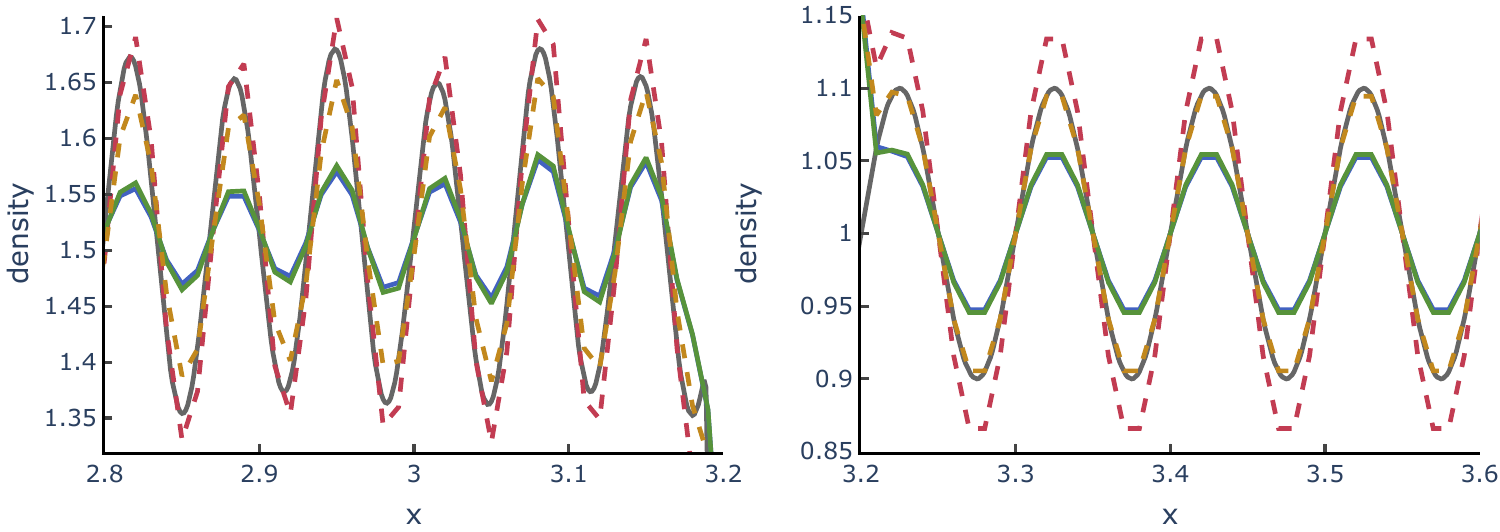}
    \caption{Numerical solutions of the Titarev--Toro shock-density wave problem \cite{titarev04} at \(t=5\) with $N=1000$ and CFL = $0.5$.}
    \label{fig:4c}
\end{figure}

The numerical experiments with the Shu--Osher and Titarev--Toro 1D problems show that WENO-ZC+ mitigates the well-known overamplification issue of WENO-Z+ \cite{luo21a,luo21b} as it is displayed in Figure \ref{fig:4c}, albeit showing a lesser dissipation when compared to WENO-ZC, in the same way WENO-Z+ enhanced WENO-Z.

\begin{figure}[H]
     \centering
     \begin{subfigure}[b]{0.49\textwidth}
         \centering
         \includegraphics[width=\textwidth]{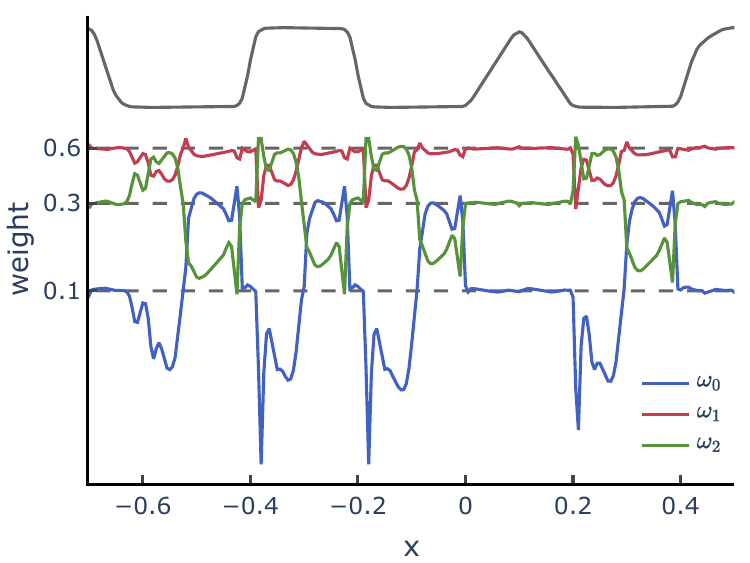}
         \caption{WENO-Z+}
     \end{subfigure}
     \hfill
     \begin{subfigure}[b]{0.49\textwidth}
         \centering
         \includegraphics[width=\textwidth]{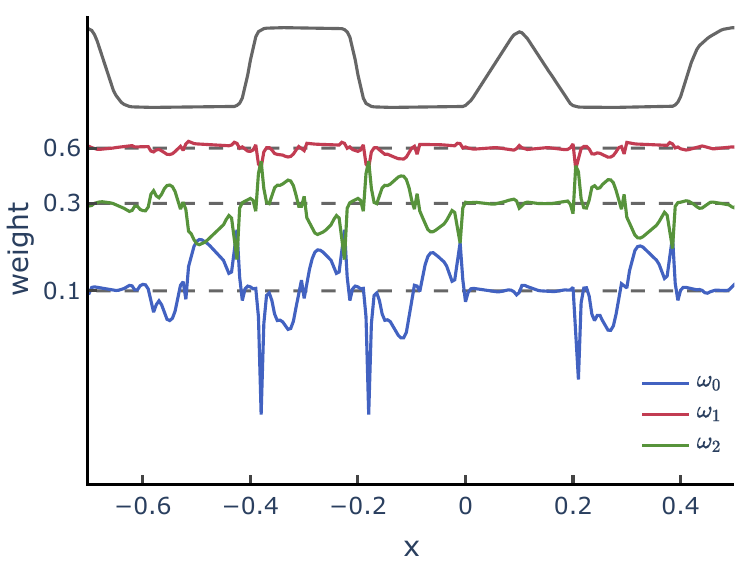}
         \caption{WENO-ZC+}
     \end{subfigure}
     \\
     \begin{subfigure}[b]{0.49\textwidth}
         \centering
         \includegraphics[width=\textwidth]{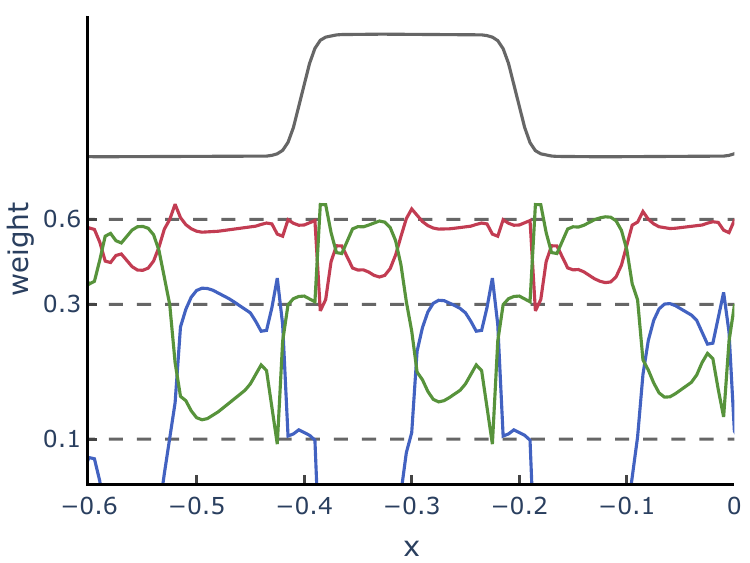}
         \caption{WENO-Z+, zoom}
     \end{subfigure}
     \hfill
     \begin{subfigure}[b]{0.49\textwidth}
         \centering
         \includegraphics[width=\textwidth]{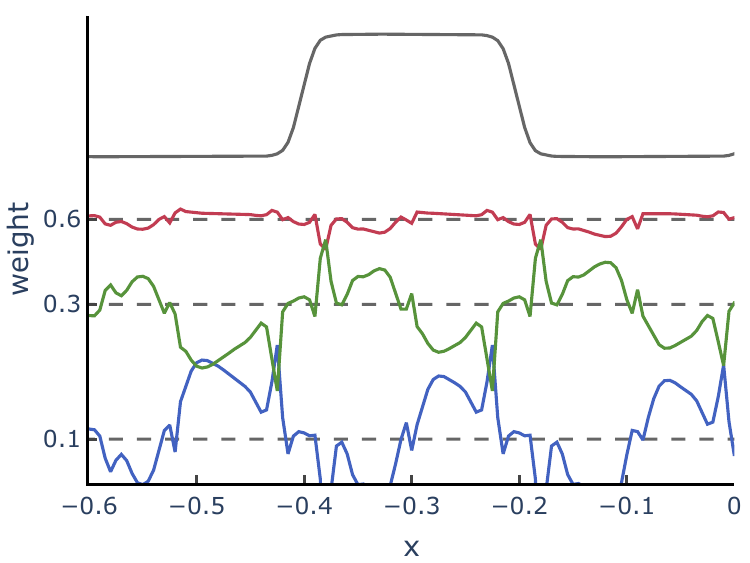}
         \caption{WENO-ZC+, zoom}
     \end{subfigure}
    \caption{Weights \(\omega_k\) of WENO-Z+ and WENO-ZC+ schemes for the GSTE advection problem at \(t=2\) with $N=400$ and CFL = $0.45$.}
    \label{fig:4d}
\end{figure}

The weights distributions for WENO-Z+ and WENO-ZC+ is shown in Figure \ref{fig:4d}, where we can observe substantial differences, since the WENO-ZC+ weights have a much more centered configuration than WENO-Z+. This indicates that the dispersion improvement indeed comes from the centered weights.

In the absence of critical points, WENO-ZC+ achieves the optimal order 5. However, at critical points of order \(\ncp=1\) and \(\ncp=2\), the optimal order is not achieved, being respectively $4$ and $3$, as shown in Table \ref{tab:weno-zcplus-ncp0}. The theoretical results are somewhat long; for this reason, they are included in 
\ref{sec:appendix-order-convergence-weno-zcplus}.
It is not the $c_k$ modification that decreases the order of WENO-ZC+ at critical points, rather the extraction of the $\Delta x$ factor at the numerator of the anti-dissipative term $\xi^{ZC+}$, for the same happens to WENO-Z+ if the grid factor is removed from $\xi^{Z+}$. This leads us to the discussion of trading a suboptimal convergence for better numerical characteristics, such as dissipation and dispersion, particularly in problems containing shocks, as has been hinted before in \cite{henrick05}, when justifying this same characteristic of WENO-M. Furthermore, WENO-D \cite{wang19} is a scheme showing optimal convergence for critical points of any order and yet its results are very similar to the WENO-Z ones, still suffering from the long-term dispersion problem.

\begin{table}[t]
\centering
\begin{tabular}{lllllll}
\toprule
 & \multicolumn{2}{l}{\(f_0(x)\) (no critical points)} & \multicolumn{2}{l}{\(f_1(x)\) (\(\ncp = 1\))} & \multicolumn{2}{l}{\(f_2(x)\) (\(\ncp = 1,\, 2\))} \\       
\cmidrule(lr){2-3} \cmidrule(lr){4-5} \cmidrule(lr){6-7}
\(1/\Delta x\) & \(L^1\) error & \(L^1\) order & \(L^1\) error & \(L^1\) order & \(L^1\) error & \(L^1\) order \\
\midrule
\(25\) & \(2.53798 \times 10^{-5}\) & --- & \(6.53262 \times 10^{-4}\) & --- & \(2.37469 \times 10^{-1}\) & --- \\
\(50\) & \(7.86300 \times 10^{-7}\) & \(5.01246\) & \(2.77585 \times 10^{-5}\) & \(4.55666\) & \(1.03968 \times 10^{-2}\) & \(4.51352\) \\
\(100\) & \(2.45530 \times 10^{-8}\) & \(5.00111\) & \(1.32262 \times 10^{-6}\) & \(4.39146\) & \(9.52538 \times 10^{-4}\) & \(3.44822\) \\
\(200\) & \(7.65285 \times 10^{-10}\) & \(5.00376\) & \(7.53484 \times 10^{-8}\) & \(4.13368\) & \(9.03779 \times 10^{-5}\) & \(3.39774\) \\
\(400\) & \(2.38826 \times 10^{-11}\) & \(5.00196\) & \(4.72904 \times 10^{-9}\) & \(3.99396\) & \(1.02251 \times 10^{-5}\) & \(3.14386\) \\
\(800\) & \(7.56719 \times 10^{-13}\) & \(4.98006\) & \(2.85514 \times 10^{-10}\) & \(4.04992\) & \(1.22910 \times 10^{-6}\) & \(3.05644\) \\
\bottomrule
\end{tabular}
\caption{Numerical results of the accuracy test for the WENO-ZC+ scheme. The expressions for the test functions \(f_0(x)\), \(f_1(x)\), and \(f_2(x)\) are given in Eqs. \eqref{eq:f0}--\eqref{eq:f2}.}
\label{tab:weno-zcplus-ncp0}
\end{table}

%% file: 4.2-ADR-Analysis.tex
\subsection{ADR Analysis}

In this section, we present a quasilinear wavenumber space analysis \cite{lele92,pirozzoli06,jia15} to get a more precise idea on the qualitative behavior of the numerical solution generated by the WENO schemes here discussed. The original von Neumann theory \cite{vichnevetsky82} considers the one-dimensional and linear advection of sinusoidal disturbances in an unbounded domain, where semidiscretization turns the PDE
\[
   u_{t}+cu_{x} = 0, \quad -\infty<x<\infty, \quad u\left(x,0\right) = u_{0}e^{i\omega x}
\]
into the set of ODEs
\begin{equation}
   \frac{du_{j}}{dt}+\frac{c}{\Delta x}\sum_{k=l}^{r}a_{j}u_{j+k}=0,
   \label{eq:ADRode}
\end{equation}
where $u_{j}=u\left(x_{j},t\right)$; $l$ and $r$ are the left and right extents of a particular explicit Finite Difference formula; and $\Delta x$ is the spatial discretization parameter. The exact solution of equation \eqref{eq:ADRode} is given by
\begin{equation}
   u_{j}\left(t\right)=\hat{u_{j}}\left(t\right)e^{i\omega j\Delta x},
   \quad\text{with}\quad
   \hat{u_{j}}\left(t\right)=e^{-i\frac{ct}{\Delta x}\phi_{j}}\hat{u_{0}},
   \label{eq:wavedynamics}
\end{equation}
where $\phi_{j}=\frac{1}{i}\sum_{k=-l}^{r}a_{l}e^{ik\omega h}$ is the modified wavenumber \cite{pirozzoli06}. It is straightforward to see from \eqref{eq:wavedynamics} that the real and imaginary parts of $\phi_{j}$ account, respectively, for the phase and amplitude dynamics of the specific discrete wave being analyzed.

In the WENO case, the conservative Finite Difference is of the type $\frac{1}{\Delta x}\left(\hat{u}_{j+\frac{1}{2}}-\hat{u}_{j-\frac{1}{2}}\right)$, where $\hat{u}_{j+\frac{1}{2}}=\hat{u}\left(u_{j-l+1},\ldots,u_{j+r}\right)$ is a nonlinear combination of the stencil values. Pirozzoli \cite{pirozzoli06} has developed the Approximate Dispersion Relation (ADR) for the nonlinear case, by observing the real and imaginary parts of the modified wavenumber for a very short integration step represented by $\Delta t$:
\[
   \phi_{j}=\frac{i\Delta x}{c\Delta t}\log\left(\frac{\hat{u_{j}}\left(\Delta t\right)}{\hat{u}_{0}}\right).
\]

\begin{figure}
     \centering
     \begin{subfigure}[b]{0.49\textwidth}
         \centering
         \includegraphics[width=\textwidth]{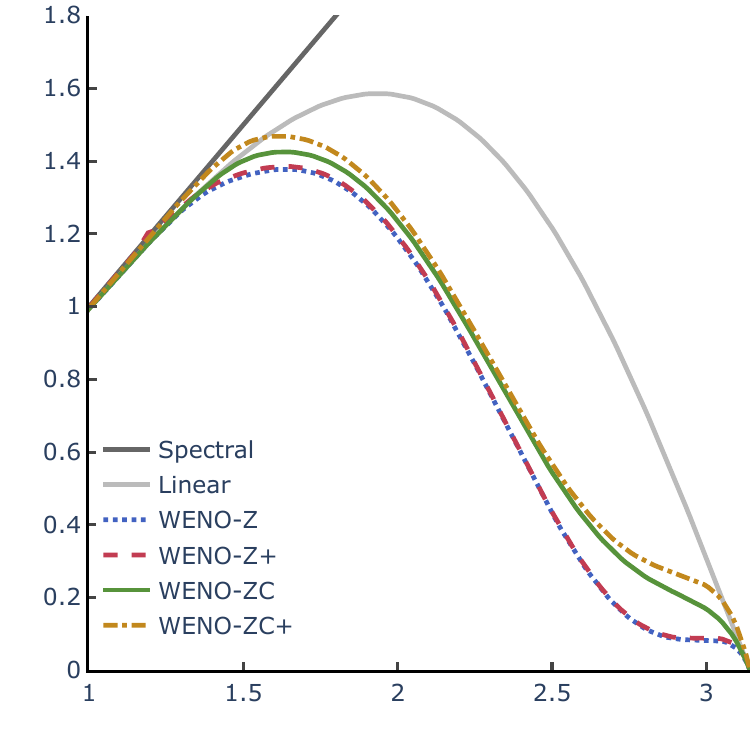}
         \caption{Real part}
     \end{subfigure} 
     \begin{subfigure}[b]{0.49\textwidth}
         \centering
         \includegraphics[width=\textwidth]{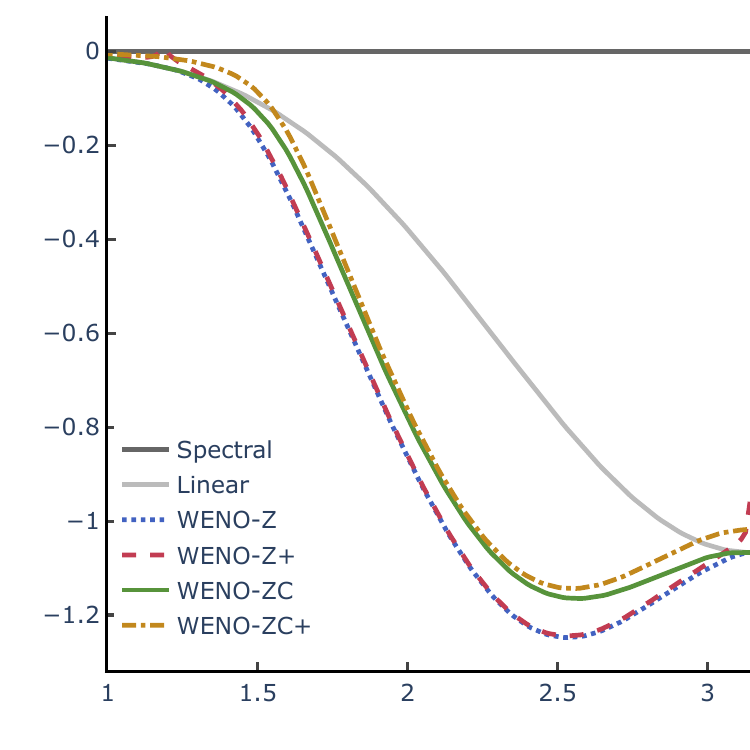}
         \caption{Imaginary part}
     \end{subfigure} \\
     \centering
     \begin{subfigure}[b]{0.49\textwidth}
         \centering
         \includegraphics[width=\textwidth]{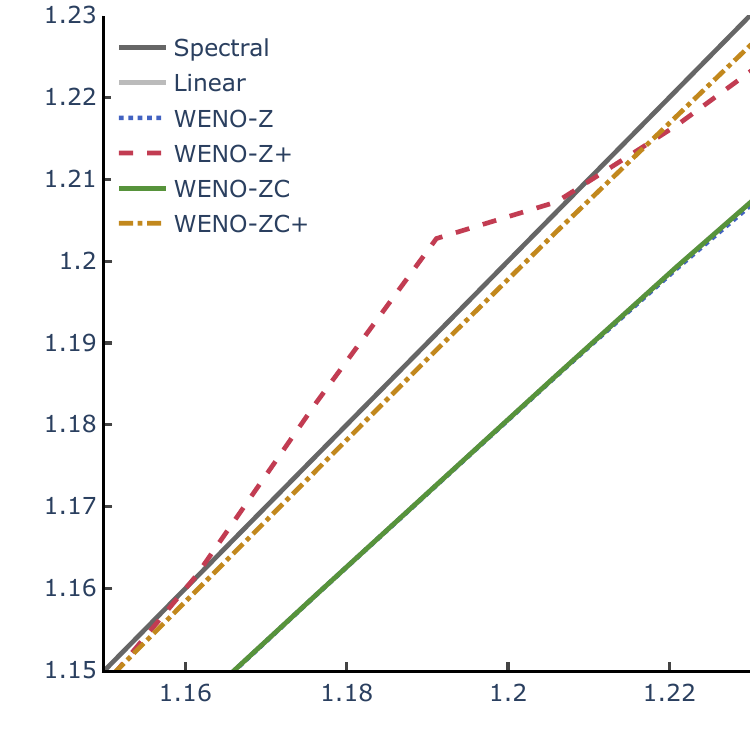}
         \caption{Real part, zoom}
     \end{subfigure}
     \begin{subfigure}[b]{0.49\textwidth}
         \centering
         \includegraphics[width=\textwidth]{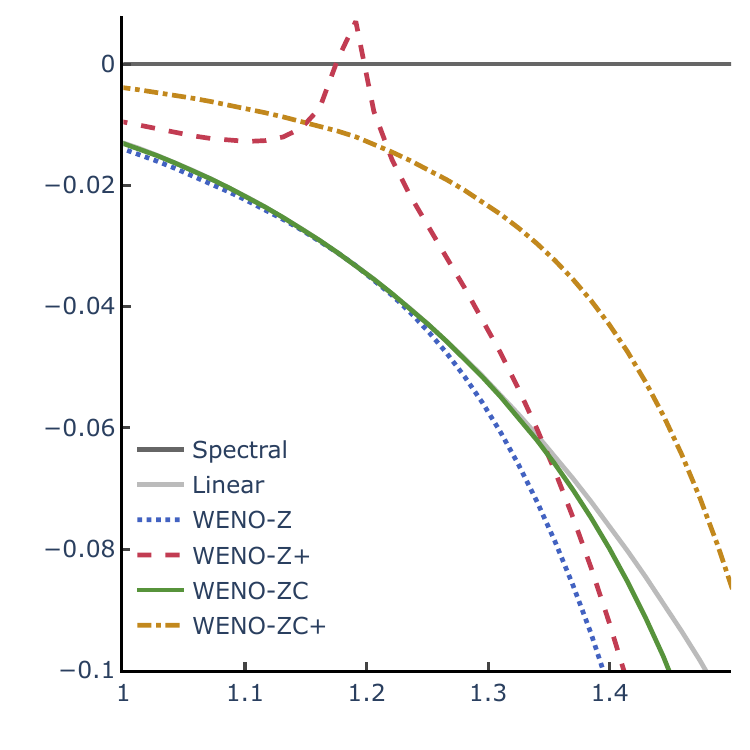}
         \caption{Imaginary part, zoom}
     \end{subfigure}
     \caption{Approximate Dispersion Relation for various WENO schemes with $N=422$, CFL = $0.5$ and $\Delta t=10^{-10}$.}
     \label{fig:4e}
\end{figure}

Although the final behavior of a particular wavenumber is affected by the nonlinear mechanisms underlying shock-capturing schemes, their leading-order effects can be observed in the ideal situation described above, with the guidance of the qualitative results of Figure \ref{fig:4e}. We see that the dispersive behavior of the centered WENO schemes is significantly improved with respect to WENO-Z, and although this quasilinear analysis cannot account for long-time integration effects, we ought to, at least partially, attribute the better numerical advection results for large final times to the superior dispersive properties of the centered schemes. In particular, in the case of WENO-ZC+, it is smaller than the central linear scheme for wavenumbers $\left|\omega\right|\in\left[1,1.5\right].$ Analogously, dissipation is greatly reduced, and marginally resolved waves also have phase and amplitude improved by the use of the centered WENO schemes.
It is also worth noting in Figure \ref{fig:4e} that WENO-ZC+ does not present the numerical artifact of WENO-Z+ at $\omega\approx 1.18$, where it shows negative dissipation, as noted in \cite{luo21b},

%% file: 4.3-Weight-Distribution.tex
\subsection{Distribution of the WENO Weigths}

As an alternative way to visualize how different WENO schemes behave throughout the numerical integration, we wish to use the relation $\sum_{k=0}^{2} \omega_{k} = 1$ and plot $\omega_{0}$ and $\omega_{2}$ in the $x$ and $y$ axis respectively, leaving $\omega_{1}$ to be implicitly represented by the aforementioned constraint.

Nevertheless, the intrinsic asymmetry of the ideal weights, $(d_{0}, d_{1}, d_{2}) = \left(1/10,\, 6/10,\, 3/10\right)$, generates distortions on the graph, making interpretation less intuitive, and to remove this bias, we divide each unnormalized weight $\alpha_{k}$ by $d_{k}$, normalizing their sum to unit as follows:
\begin{equation}
    \begin{aligned}
        \lambda_{k}^{*} &= \frac{\alpha_{k}}{d_{k}} & 
        \text{and} && 
        \lambda_{k} &= \frac{\lambda_{k}^{*}}{\sum_{j=0}^2\lambda_{j}^{*}}, &&
        k = 0,\, 1,\, 2,
    \end{aligned}
\end{equation}
This allows us to generate an indirect graphical representation of the weights, which we call the distribution map of the weights.

For instance, by using a specific WENO scheme, we obtain the numerical solution for the Tirarev--Toro shock-density wave problem at $t_{i}=i/2$, for $i \in \{0, 1, \dots, 10\}$, with $N=1000$ and CFL = $0.5$, and compute the $\lambda_{k}$ for all grid points. Figure \ref{fig:4f} shows the results for some of the WENO schemes considered in this work. The horizontal axis corresponds to the values of \(\lambda_0\); the vertical, to \(\lambda_2\); and \(\lambda_1\), as stated before, is given by \(1 - \lambda_0 - \lambda_2\). It is worth noting that, by design of the WENO weights, \((\lambda_0, \lambda_2)\) should be close to $(1,0)$, $(0,1)$ or $(0,0)$ in the presence of shocks or discontinuities, and be 
close to $(1/3,1/3)$ at the smooth parts of the numerical solution. 
The Titarev--Toro problem is an interesting case since, due to its high-frequency oscillations, \((\lambda_0, \lambda_2)\) display a wide range of values in between these extreme cases above. This allows us to see the behavior of WENO schemes at smoother, yet not completely smooth, parts of the solution.

\begin{figure}[ht]
    \centering
    \begin{subfigure}[b]{0.5\textwidth}
         \centering
         \graficopesado{width=\textwidth}{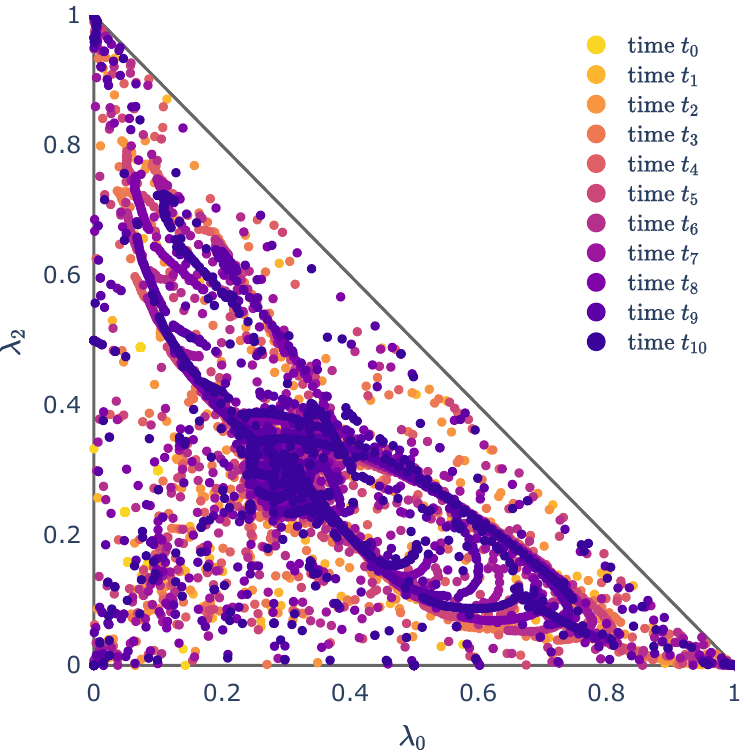}
         \caption{WENO-JS}
         \label{fig:4f-WENO-JS}
     \end{subfigure}\hfill
    \begin{subfigure}[b]{0.5\textwidth}
         \centering
         \graficopesado{width=\textwidth}{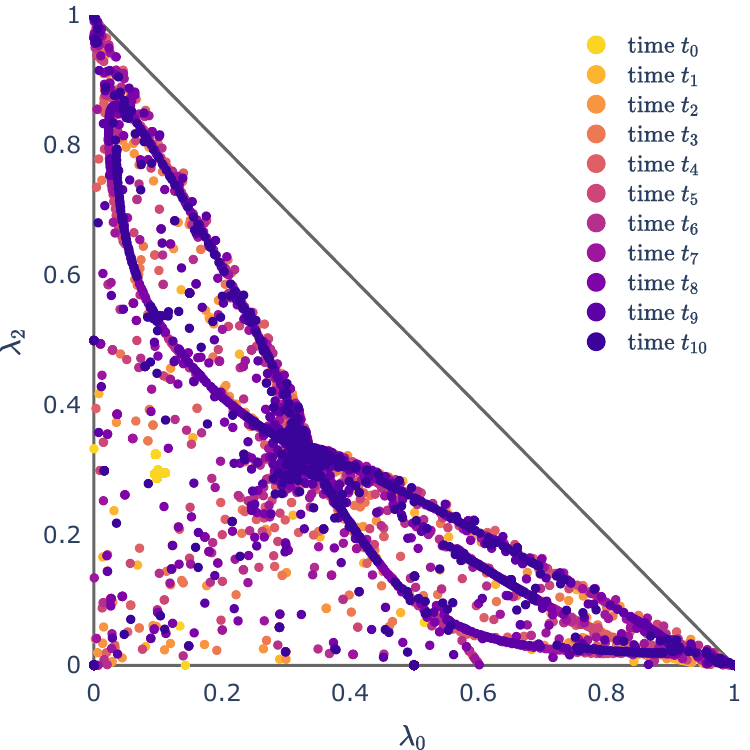}
         \caption{WENO-Z}
         \label{fig:4f-WENO-Z}
     \end{subfigure}\hfill
\end{figure}
\begin{figure}\ContinuedFloat
    \centering
    \begin{subfigure}[b]{0.5\textwidth}
         \centering
         \graficopesado{width=\textwidth}{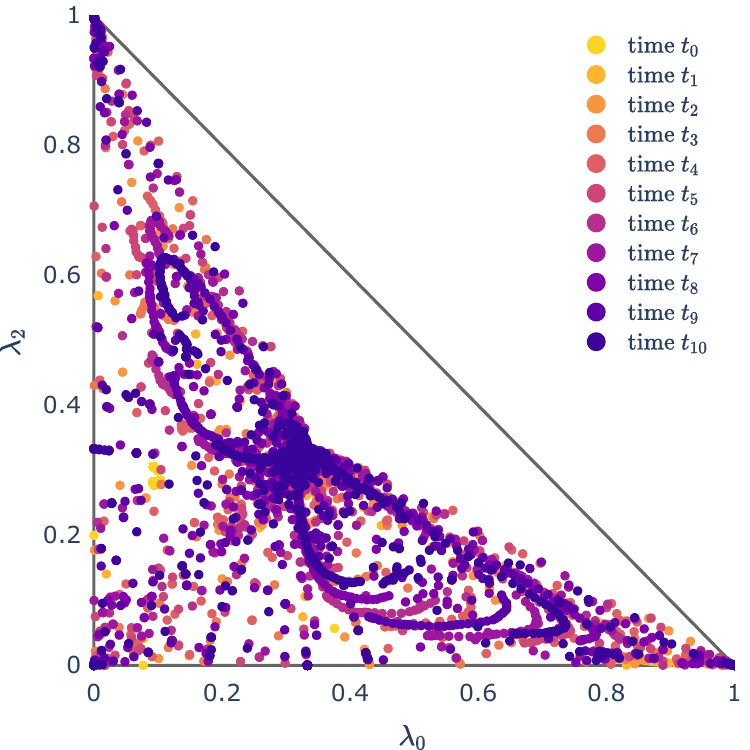}
         \caption{WENO-C}
         \label{fig:4f-WENO-C}
     \end{subfigure}\hfill
    \begin{subfigure}[b]{0.5\textwidth}
         \centering
         \graficopesado{width=\textwidth}{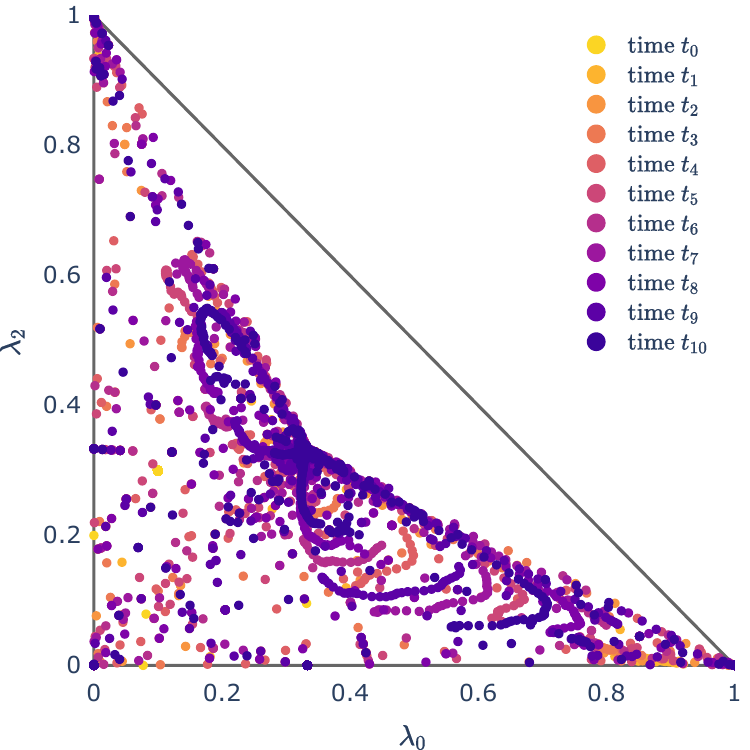}
         \caption{WENO-ZC}
         \label{fig:4f-WENO-ZC}
     \end{subfigure}\hfill
\end{figure}
\begin{figure}\ContinuedFloat
    \centering
    \begin{subfigure}[b]{0.5\textwidth}
         \centering
         \graficopesado{width=\textwidth}{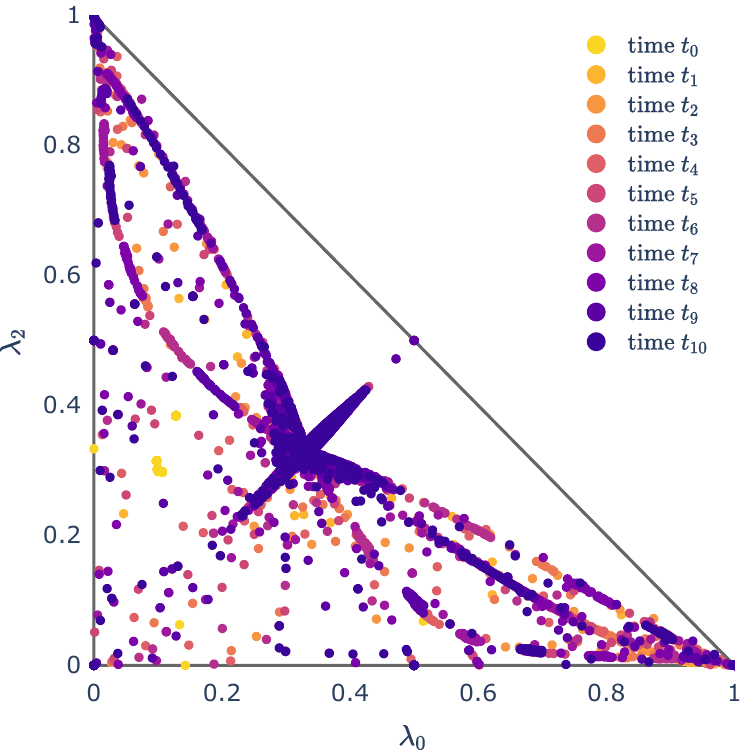}
         \caption{WENO-Z+}
         \label{fig:4f-WENO-Z+}
     \end{subfigure}\hfill
     \begin{subfigure}[b]{0.5\textwidth}
         \centering
         \graficopesado{width=\textwidth}{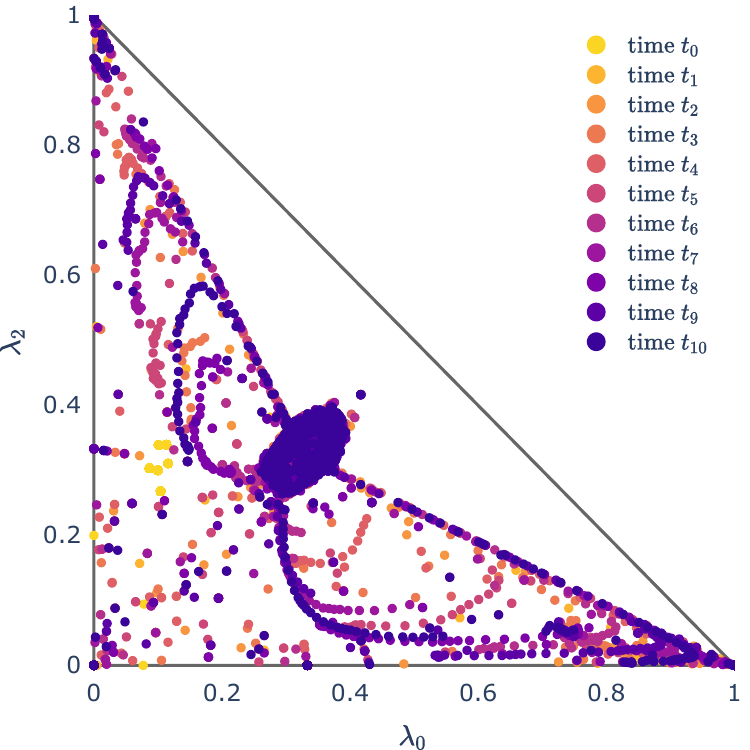}
         \caption{WENO-ZC+}
         \label{fig:4f-WENO-ZC+}
     \end{subfigure} 
     
     \caption{Distribution map of the weights for the Titarev--Toro problem at different times \(t_i = i/2\), with $N=1000$ and CFL = $0.5$.}
     \label{fig:4f}
\end{figure}

For some insight, let us compare the distribution map for some of the methods represented in Figure \ref{fig:4f}. When comparing Figs. \ref{fig:4f-WENO-JS} and \ref{fig:4f-WENO-Z}, the most visually striking difference between WENO-JS and WENO-Z is the empty region near the diagonal line. For points near the diagonal, we have \(\lambda_0 + \lambda_2 \approx 1\), which implies \(\lambda_1 \approx 0\). In other words, for the points in this region, we would have most of the nonlinear weights allocated to the lateral substencils. However, this would paradoxically imply that the lateral substencils are smooth and the central substencil contains some discontinuity, even though $S_{1} \subset (S_{0} \cup S_{2})$. For the WENO-JS scheme, this strange allocation of weights happens because each $\alpha_{k}^{JS}$ only carries information related to $S_{k}$, but the WENO-Z scheme uses the global smoothness indicator $\tau$ to transmit information from one substencil to another, completely avoiding this case. And considering the schemes derived from the WENO-Z, e.g. the WENO-ZC scheme in Figure \ref{fig:4f-WENO-ZC}, we can similarly see a mostly blank region near the diagonal.

Next, when comparing Figures \ref{fig:4f-WENO-Z} and \ref{fig:4f-WENO-Z+}, the most noticeable change from the WENO-Z to the WENO-Z+ is the emergence of a concentration of points satisfying $\lambda_{0}^{Z+} \approx \lambda_{2}^{Z+}$. To better understand this pattern, consider the additional anti-dissipative term $\eta\frac{\beta_{k}}{\tau + \varepsilon}$ for the WENO-Z+, that has a significant impact on the final distribution of the weights for $\tau \approx 0$. It is easy to check that
\begin{equation}
    \tau = |\beta_{2}-\beta_{0}| \approx 0 \Longrightarrow \beta_{0} \approx \beta_{2} \Longrightarrow (\lambda_{0}^{*})^{Z+} \approx (\lambda_{2}^{*})^{Z+} \Longrightarrow \lambda_{0}^{Z+} \approx \lambda_{2}^{Z+},
\end{equation}
thus explaining the resulting pattern. If we then compare Figures \ref{fig:4f-WENO-Z+} and \ref{fig:4f-WENO-ZC+}, we see that WENO-ZC+ presents a similar behavior, but with the modified points more concentrated near $\left(1/3, 1/3\right)$. We attribute this change to the upper limit on the size of its anti-dissipative term, which restrains his overall impact on the scheme's weights.

Table \ref{tab:Relative-error-of all} completes Table \ref{tab:GSTE, weights Relative-error} from Section \ref{sec:dispersion-error}. It shows the relative error of the weights of all schemes with respect to the ideal weights $d_{k}$. Notice that the weights $\omega_{k}$ are getting, in average, closer and closer to the ideal weights as the scheme progress from WENO-JS, WENO-Z, WENO-ZC and WENO-ZC+. The schemes seem to run from the more to the less dissipative as well. This indicates the existence of much room for improvement, but surprisingly, the results with WENO-D shows that making the nonlinear weights as closer as possible to the ideal ones not necessarily implies in a better scheme.
\begin{table}
   \centering
   \begin{tabular}{lcccc}
   \toprule
   & $e_{0}$ & $e_{1}$ & $e_{2}$ & $e_{0} + e_{1} + e_{2}$ \\
   \midrule
   WENO-JS  & 2.29721 & 0.38900 & 1.25938 & 3.94561 \\
   WENO-JSC & 1.67972 & 0.39635 & 1.17703 & 3.25311 \\
   WENO-Z+  & 1.51620 & 0.26511 & 0.85672 & 2.63804 \\
   WENO-Z   & 1.52174 & 0.25985 & 0.85565 & 2.63724 \\
   WENO-C   & 1.03149 & 0.17170 & 0.67824 & 1.88145 \\
   WENO-ZC  & 0.68140 & 0.11323 & 0.45194 & 1.24658 \\
   WENO-ZC+ & 0.60492 & 0.10601 & 0.40117 & 1.11211 \\
   \bottomrule
   \end{tabular}
   \caption{Relative error of the weights considering different WENO schemes for the solution of the GSTE at time $t=2$ with $N=400$ and CFL = $0.45$ obtained with the WENO-Z scheme.}
   \label{tab:Relative-error-of all}
\end{table}

%% file: 5-Numerical-Experiments.tex
\section{Numerical Experiments} \label{sec:numerical-experiments}
In this section, we aim to show that the new centered schemes behave on a stable and nonoscillatory manner when facing the standard set of numerical tests of the current literature on WENO schemes, i.e., WENO-ZC and WENO-ZC+ were able to complete all the one-dimensional Euler 1D experiments, along with the classic two-dimensional Rayleigh--Taylor Instability and Double Mach Reflection tests, showing improved results over the WENO-Z and WENO-Z+ ones.

\subsection{Shock Tube Tests}
We start this section showing the numerical results of the new methods when applied to the standard set of Euler 1D shock tube tests. These experiments test the ability of a compressible code to represent shocks, contact discontinuities and rarefaction waves. The Lax and Sod cases, described below, have exact dynamics, obtained through the solutions of Riemann problems. With them, it is possible to compare the distinct WENO schemes with respect to shock-capturing and oscillatory behavior.

In the Riemann problem of Lax \cite{lax54}, the initial conditions are given by:
\[
   (\rho_{0}, u_{0}, p_{0}) =
      \begin{cases}
         (0.445,\: 0.698,\: 3.528), & x \leq 0, \\
         (0.5,\: 0,\: 0.571), & x > 0,
      \end{cases}
\]
where $x \in [-0.5,\: 0.5]$, with free boundary conditions and final time $T=0.13$. The reference solution at the final time consists, from left to right, of a rarefaction wave, a contact discontinuity and a shock.

\begin{figure}[ht]
    \centering
    \includegraphics[width=\textwidth]{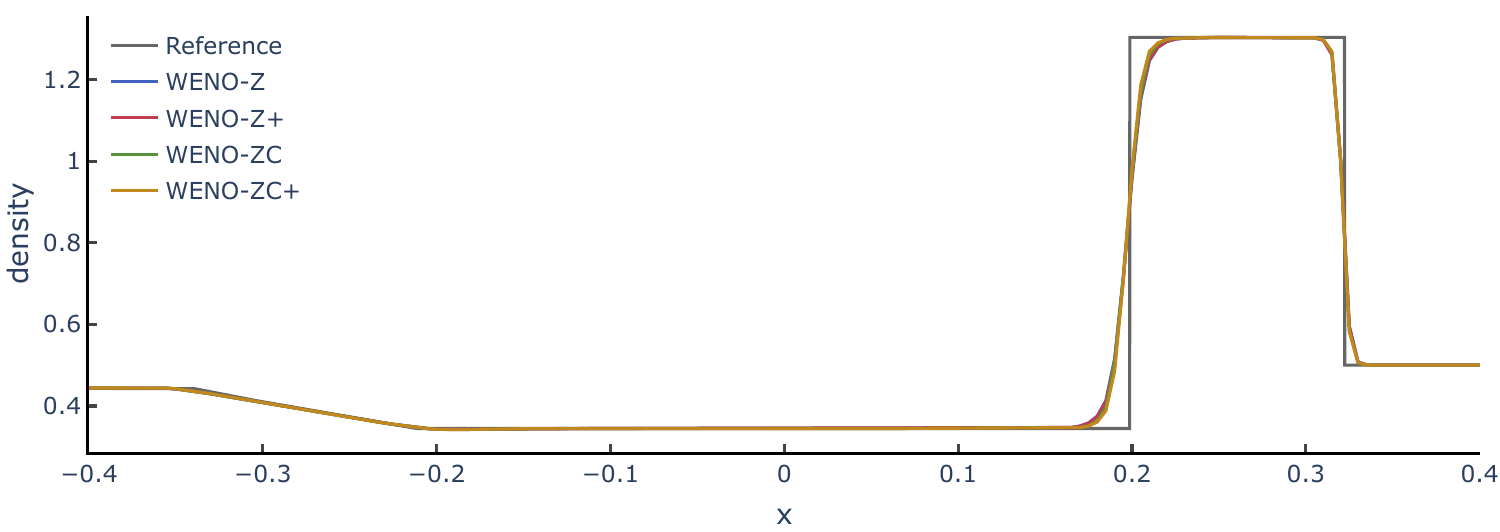} \\
    \includegraphics[width=\textwidth]{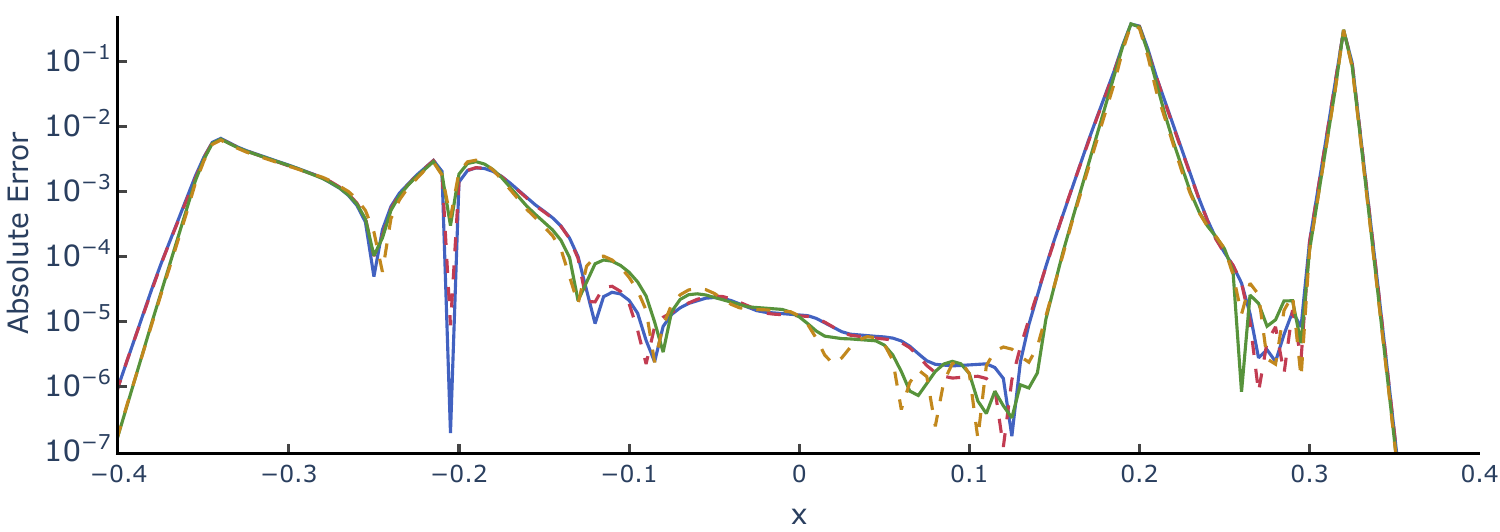}
    \caption{Numerical solutions of the Riemann problem of Lax by the WENO-Z, WENO-ZC and WENO-ZC+ schemes with $N=200$ and CFL = $0.5$. Top: density \(\rho\). Bottom: absolute error in comparison with the solution obtained by an exact Riemann solver.}
    \label{fig:Lax}
\end{figure}

Figure \ref{fig:Lax} shows that WENO-ZC and WENO-ZC+ capture the contact discontinuity with higher precision than WENO-Z, which is also confirmed by the graph of their absolute errors, where the slighter oscillatory behavior of the new schemes occurs at a very small scale and only depicts their lesser dissipation.

\begin{figure}[ht]
    \centering
    \includegraphics[width=\textwidth]{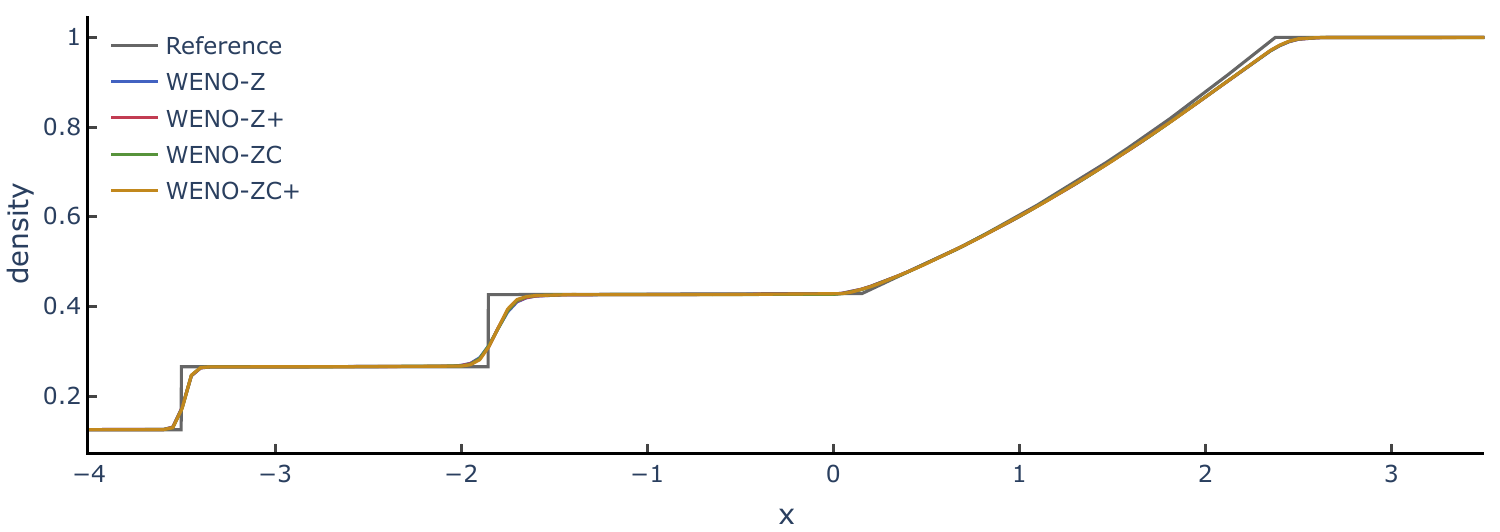} \\
    \includegraphics[width=\textwidth]{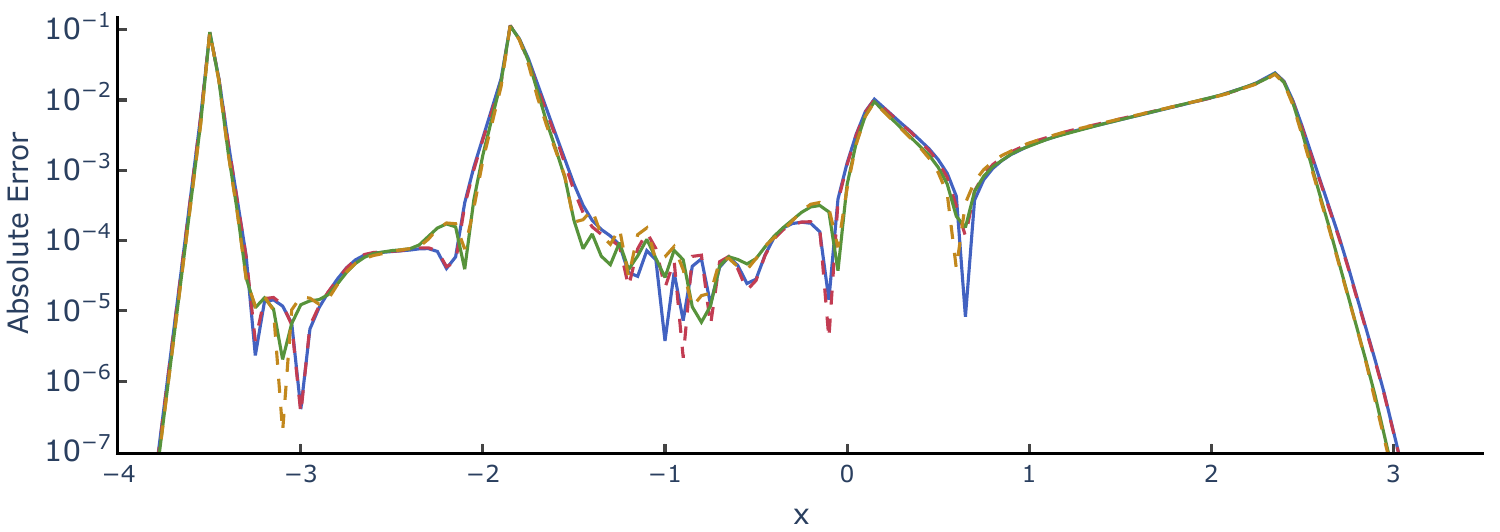}
    \caption{Numerical solutions of the Riemann problem of Sod by the WENO-Z, WENO-ZC and WENO-ZC+ schemes with $N=200$ and CFL = $0.5$. Top: density \(\rho\). Bottom: absolute error in comparison with the solution obtained by an exact Riemann solver.}
    \label{fig:Sod}
\end{figure}

Analogously, in the Riemann problem of Sod \cite{sod78}, the initial conditions are given by:
\[
   (\rho_{0}, u_{0}, p_{0}) =
      \begin{cases}
         ({1}/{8},\: 0,\: {1}/{10}), & x \leq 0, \\
         (1,\: 0,\: 1), & x > 0,
      \end{cases}
\]
where $x \in [-0.5,\: 0.5]$, with free boundary conditions and final time $T=2$. Its exact solution is computed in the same way as in the Lax problem, having the same structure of shock, contact discontinuity, and rarefaction wave. The numerical results are depicted in Figure \ref{fig:Sod}, showing once again the smaller dissipation of the centered schemes.

The Shu--Osher shock-density wave problem \cite{shu89} simulates a Mach 3 shock wave passing through a sinusoidal wave distribution of the density. As a result, a region with shocklets and another with higher-frequency waves are formed behind the shock. The initial conditions are
\[
   (\rho_{0}, u_{0}, p_{0}) =
      \begin{cases}
         ({27}/{7},\: 4\sqrt{35}/9,\: {31}/{3}), & x < -4, \\
         (1+\sin(5x)/5,\: 0,\: 1), & x \geq -4,
      \end{cases}
\]
and the domain of this test is \(x \in [-5, 5]\). 

\begin{figure}[ht]
    \centering
    \includegraphics[width=\textwidth]{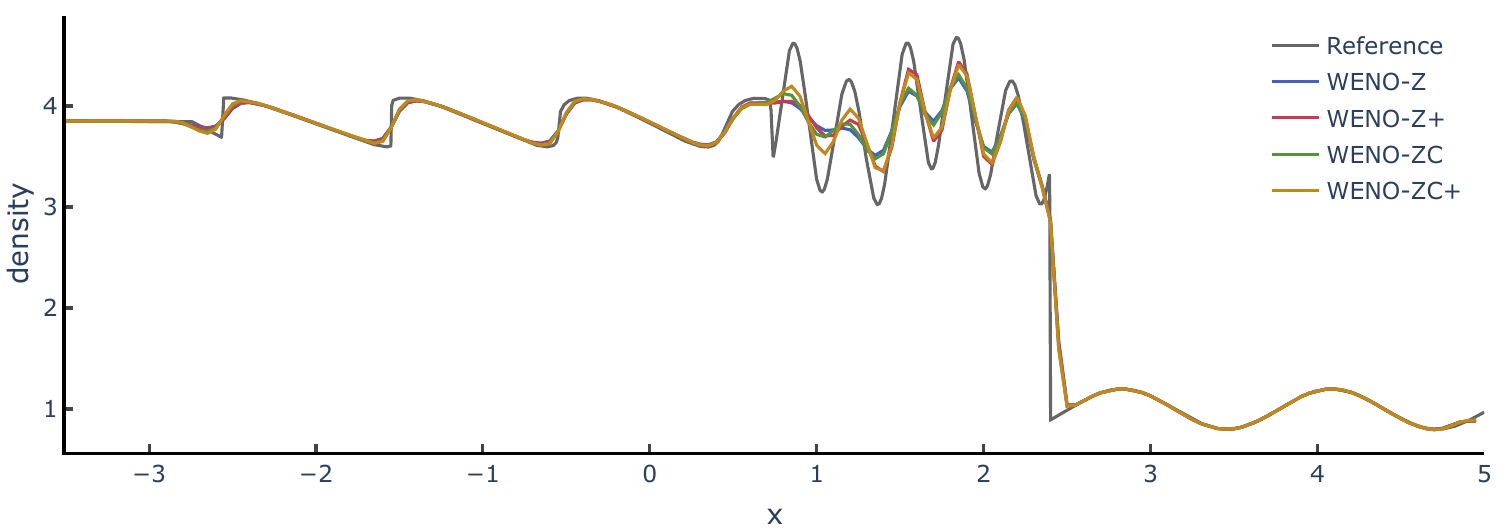} \\
    \includegraphics[width=\textwidth]{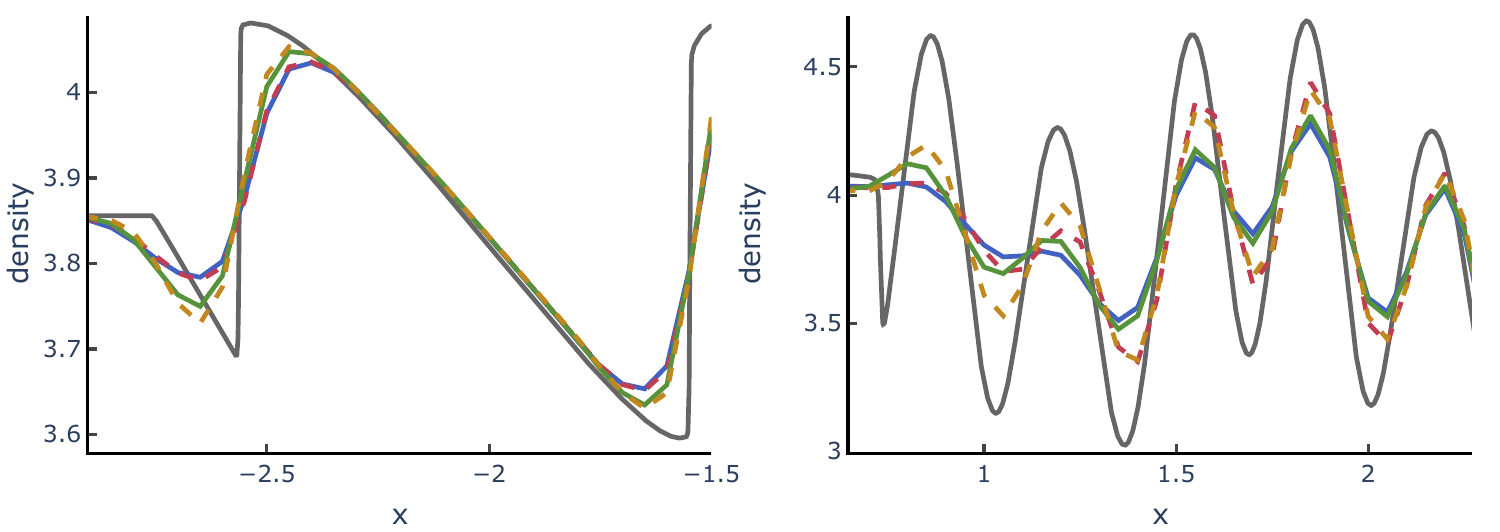}
    \caption{Numerical solutions of the Shu--Osher shock-density wave problem \cite{shu89} by the WENO-Z, WENO-ZC and WENO-ZC+ schemes at \(t=1.8\) with $N=200$ and CFL = $0.5$.}
    \label{fig:Shu-Osher}
\end{figure}

Figure \ref{fig:Shu-Osher} shows the numerical results for the Shu--Osher problem. It can be seen that WENO-ZC+ achieves the best results, followed by WENO-Z+---even though these two schemes are the ones which theoretically lose accuracy at critical points. Also, notice that WENO-ZC is comparable to---but less dissipative than---WENO-Z.

\subsection{Stability Tests}

In the following tests, schemes with unstable tendencies often fail to converge. The first one consists of two interacting blast waves with strong shocks in the solution that are computationally hard to solve \cite{woodward84}:
\[
   \left(\rho_{0},u_{0},p_{0}\right) =
      \begin{cases}
         \left(1,\: 0,\: 1000\right), & x < 0.1,\\
         \left(1,\: 0,\: 100\right), & x > 0.9,\\
         \left(1,\: 0,\: 0.01\right), & \text{otherwise},
      \end{cases}
\]
where $x\in\left[0,1\right]$, the final time is $T=0.038$ and reflexive boundary conditions are applied.

\begin{figure}[ht]
    \centering
    \includegraphics[width=\textwidth]{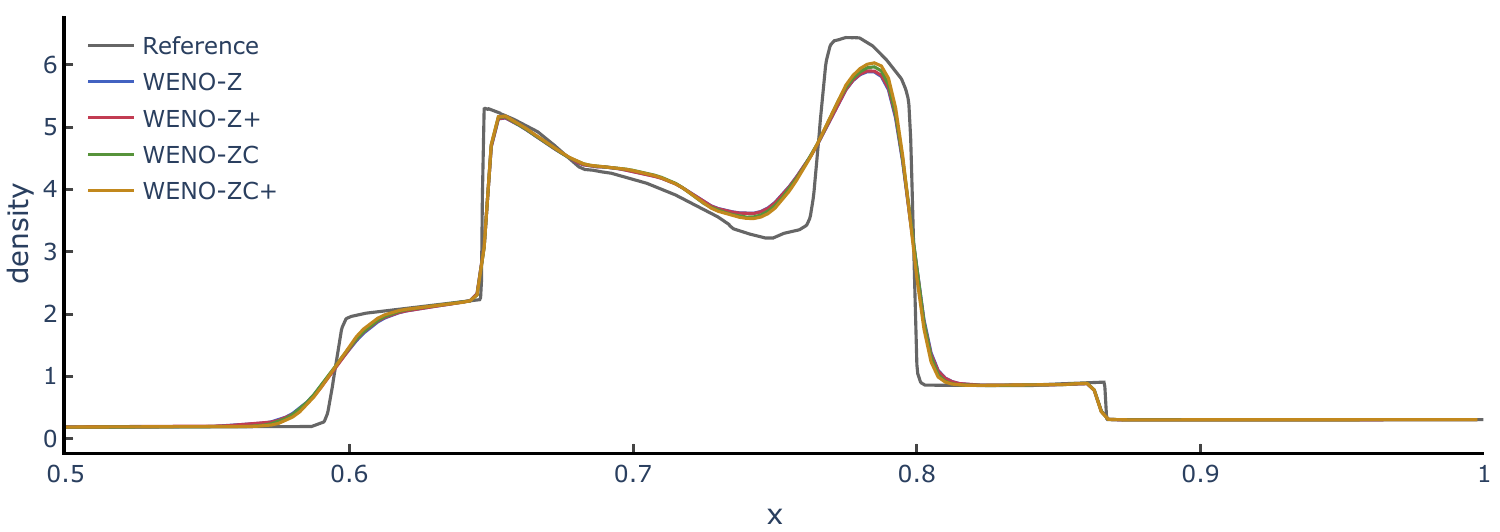} \\
    \includegraphics[width=\textwidth]{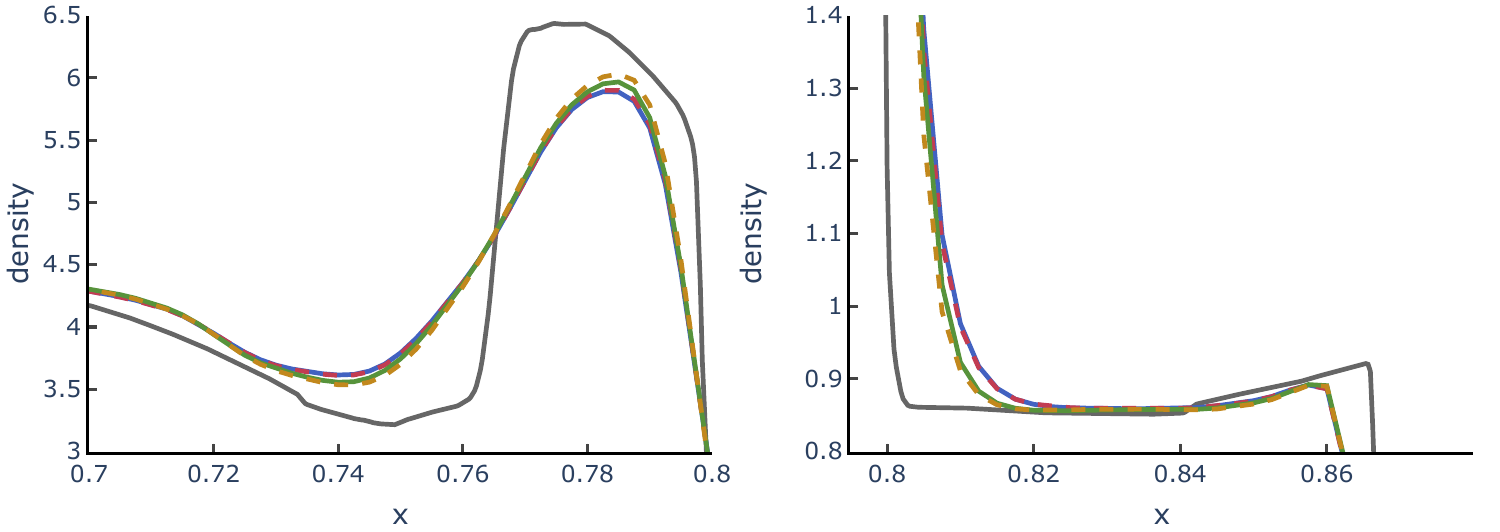}
    \caption{Numerical solutions of the two interacting blast waves problem at \(t = 0.038\) with $N=400$ and CFL = $0.5$.}
    \label{fig:Blast-Waves}
\end{figure}

WENO-ZC and WENO-ZC+ showed numerical solutions which are superior to WENO-Z, as can be seen in Figure \ref{fig:Blast-Waves}. It is worth noting that the WENO-ZC solution improves over WENO-Z and WENO-Z+, which are practically the same, and WENO-ZC+ improves it even further, a behavior that will be seen again in the experiments below.

The second one is the 1D Sedov blast wave problem \cite{sedov93}. Its initial condition simulates a delta distribution of a very strong pressure gradient at the origin. The solution contains a quasi-vacuum state around \(x=0\), and even very small oscillations make the numerical method diverge. The initial condition is
\[
   \left(\rho_{0},u_{0},p_{0}\right) =
      \begin{cases}
         \left(1,\: 0,\: 4\times10^{-13}\right), & \delta<|x|\leq 2,\\
         \left(1,\: 0,\: 2.56\times10^{8}\right), & |x|\leq\delta,
      \end{cases}
\]
with \(\delta = \dx/2\), \(\gamma = 7/5\) and final time \(T = 10^{-3}\).

\begin{figure}[ht]
    \centering
    \includegraphics[width=\textwidth]{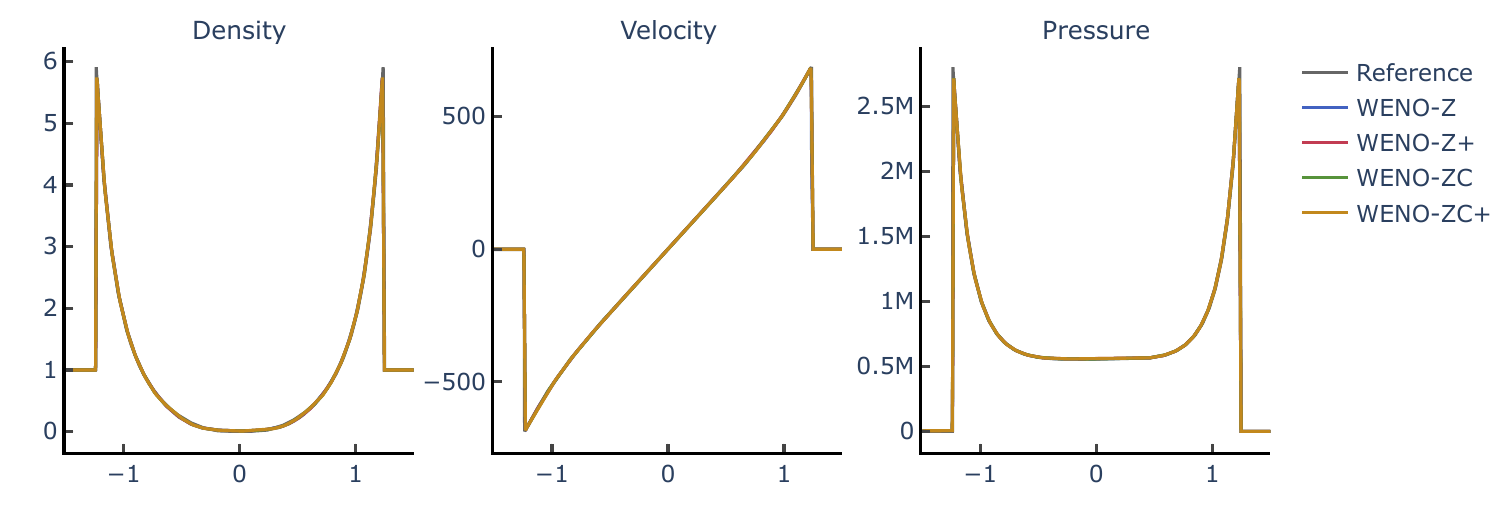} \\
    \includegraphics[width=\textwidth]{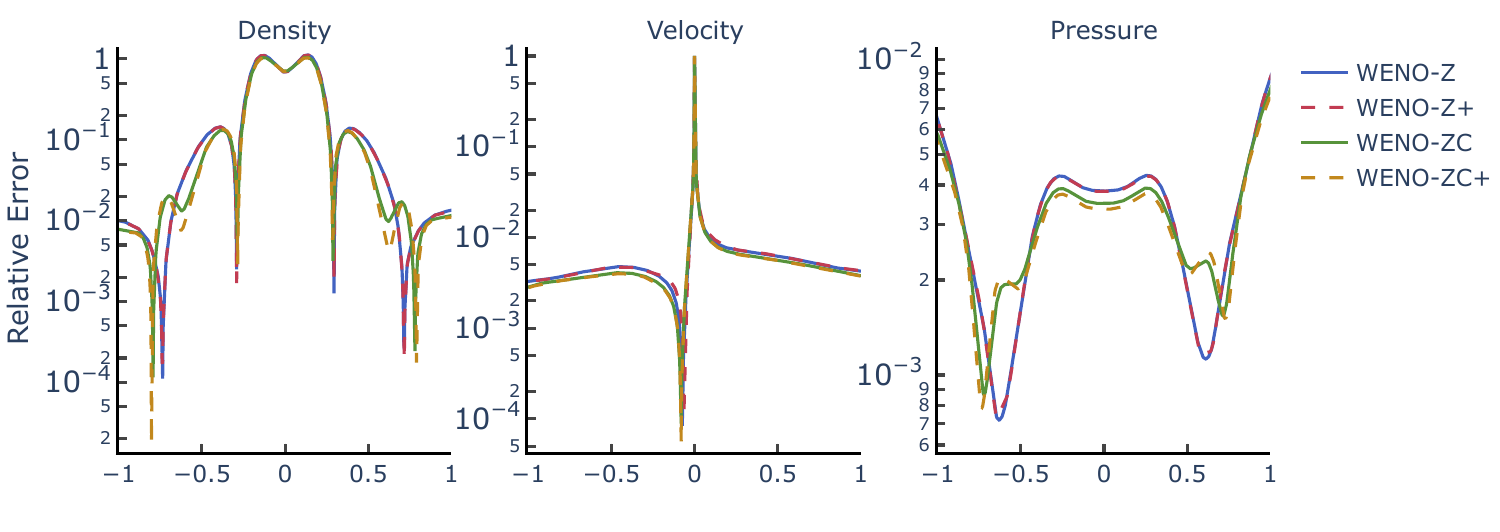}
    \caption{Numerical solutions of the Sedov problem \cite{sedov93} with $N=1250$ and CFL = $0.5$.}
    \label{fig:Sedov}
\end{figure}

Both WENO-ZC and WENO-ZC+ pass this stability test. As Figure \ref{fig:Sedov} shows, the difference in the solutions are very small, but WENO-ZC+ is marginally better than WENO-ZC, which is also slightly better than WENO-Z.

\begin{remark}
   All tests were also run with several other grid sizes, namely: $N = 100$, $200$, $400$, $800$, $1200$, and $1600$. The results were qualitatively the same; e.g., in the Shu--Osher test, WENO-ZC+ showed better results than the other schemes for all tested grid sizes, not only for the showcased \(N=200\); all schemes passed the stability tests for this whole set of grid sizes; etc. For brevity, these other results will not be shown here.
\end{remark}

\subsection{2D Euler Tests}
For the Rayleigh--Taylor instability test, the initial conditions are
\[
   \left(\rho_{0},u_{0},v_{0},p_{0}\right) =
      \begin{cases}
         \left(2,\: 0,\: -0.025\,a\cos(8\pi x),\: 2y+1\right), & y < 1/2,\\
         \left(1,\: 0,\: -0.025\,a\cos(8\pi x),\: y+3/2\right), & y \geq 1/2,
      \end{cases}
\]
where \(a = \sqrt{\gamma p/\rho}\) and \(\gamma = 5/3\). The computational domain is \((x, y) \in [0, 0.25] \times [0, 1]\), with reflective boundaries on the left and right, and fixed values on the top and bottom. The source term \(s(x, y, t) = (0, 0, \rho, \rho v)\)
is added to the right-hand side of the equation. The final time is \(T = 1.95\).

\begin{figure}[p]
   \centering
   \begin{subfigure}[b]{4.5cm}
      \centering
      \includegraphics[width=\textwidth]{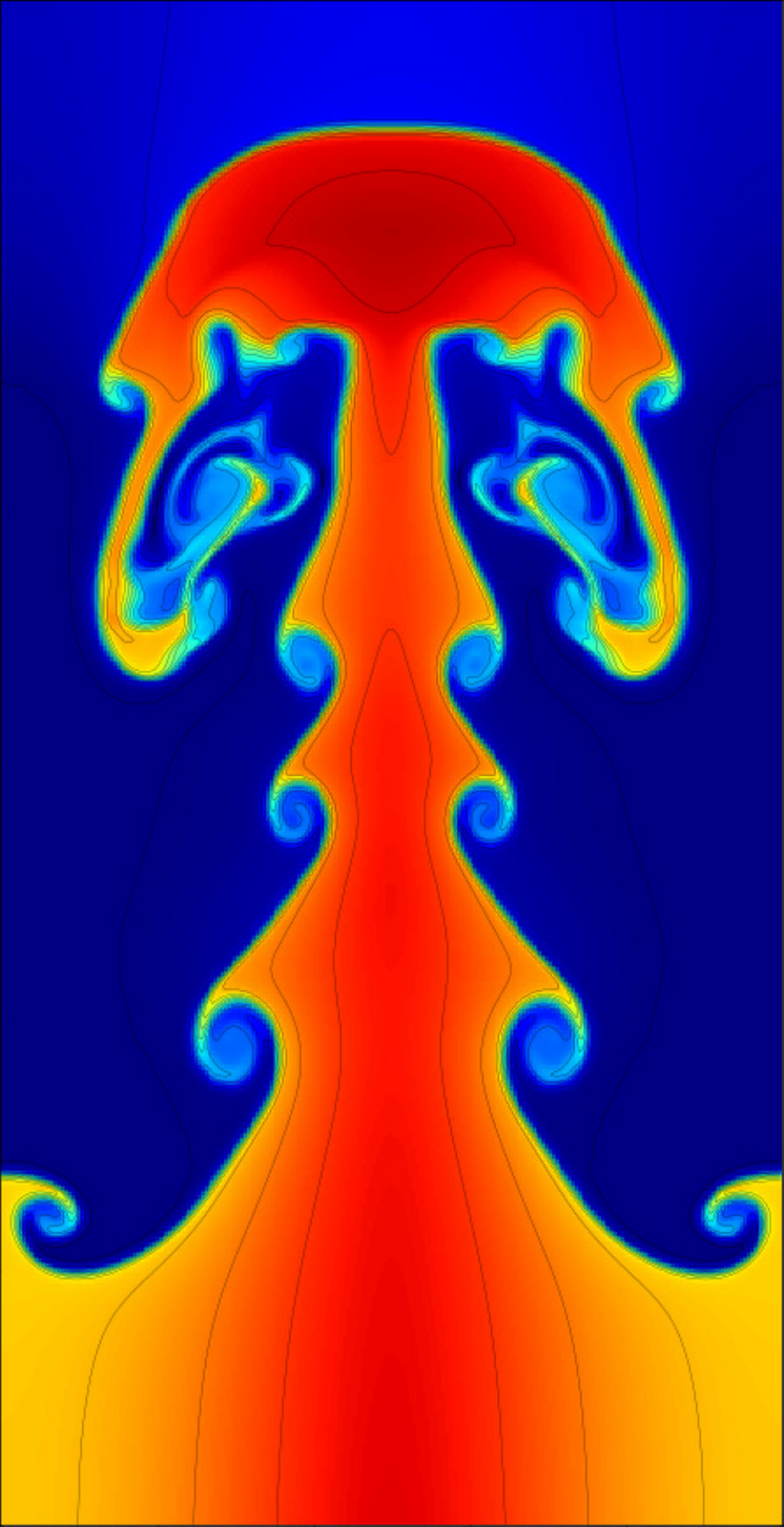}
      \caption{WENO-Z}
   \end{subfigure}
   \hspace{0.5cm}
   \begin{subfigure}[b]{4.5cm}
      \centering
      \includegraphics[width=\textwidth]{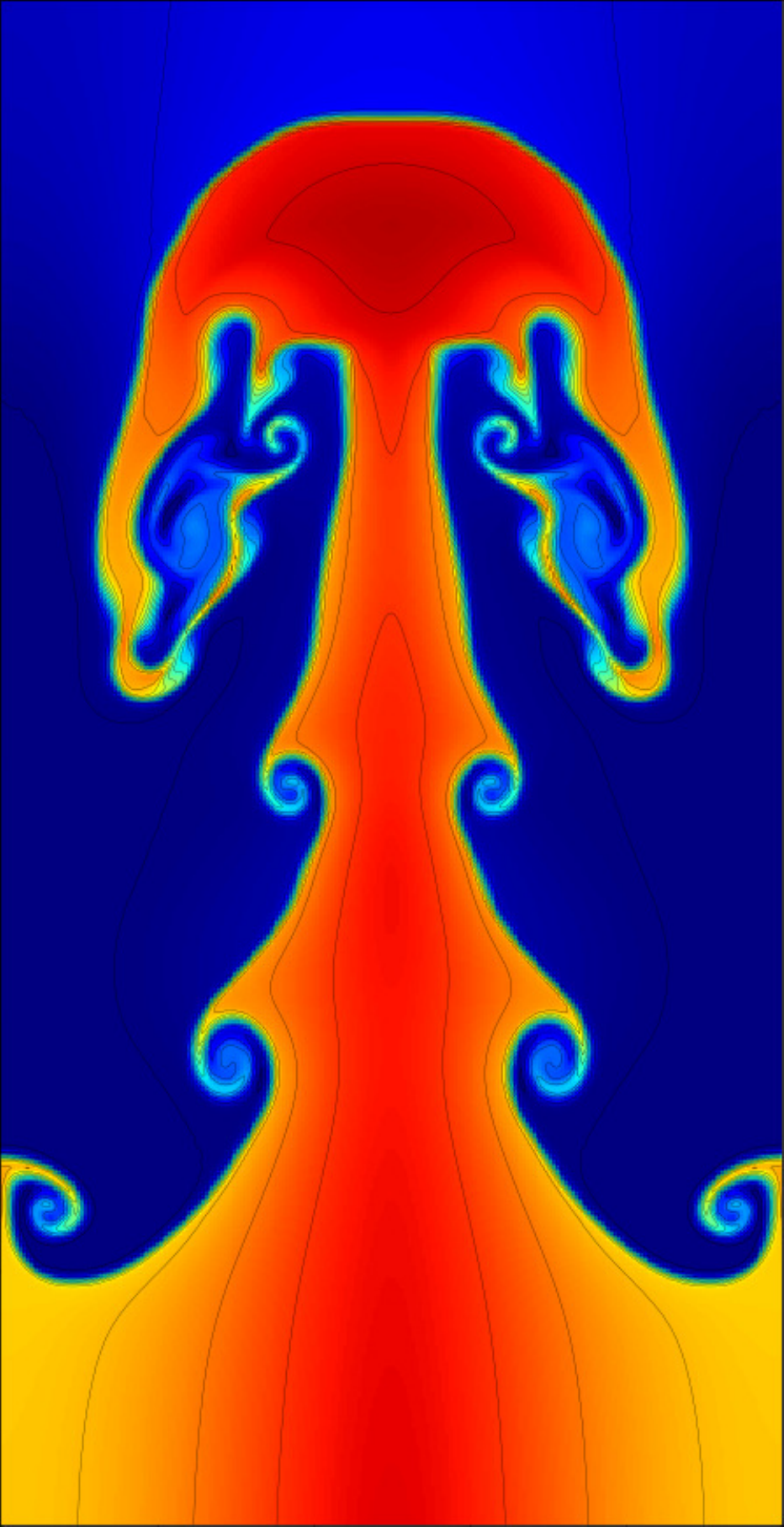}
      \caption{WENO-ZC}
   \end{subfigure}
   \\
   \begin{subfigure}[b]{4.5cm}
      \centering
      \includegraphics[width=\textwidth]{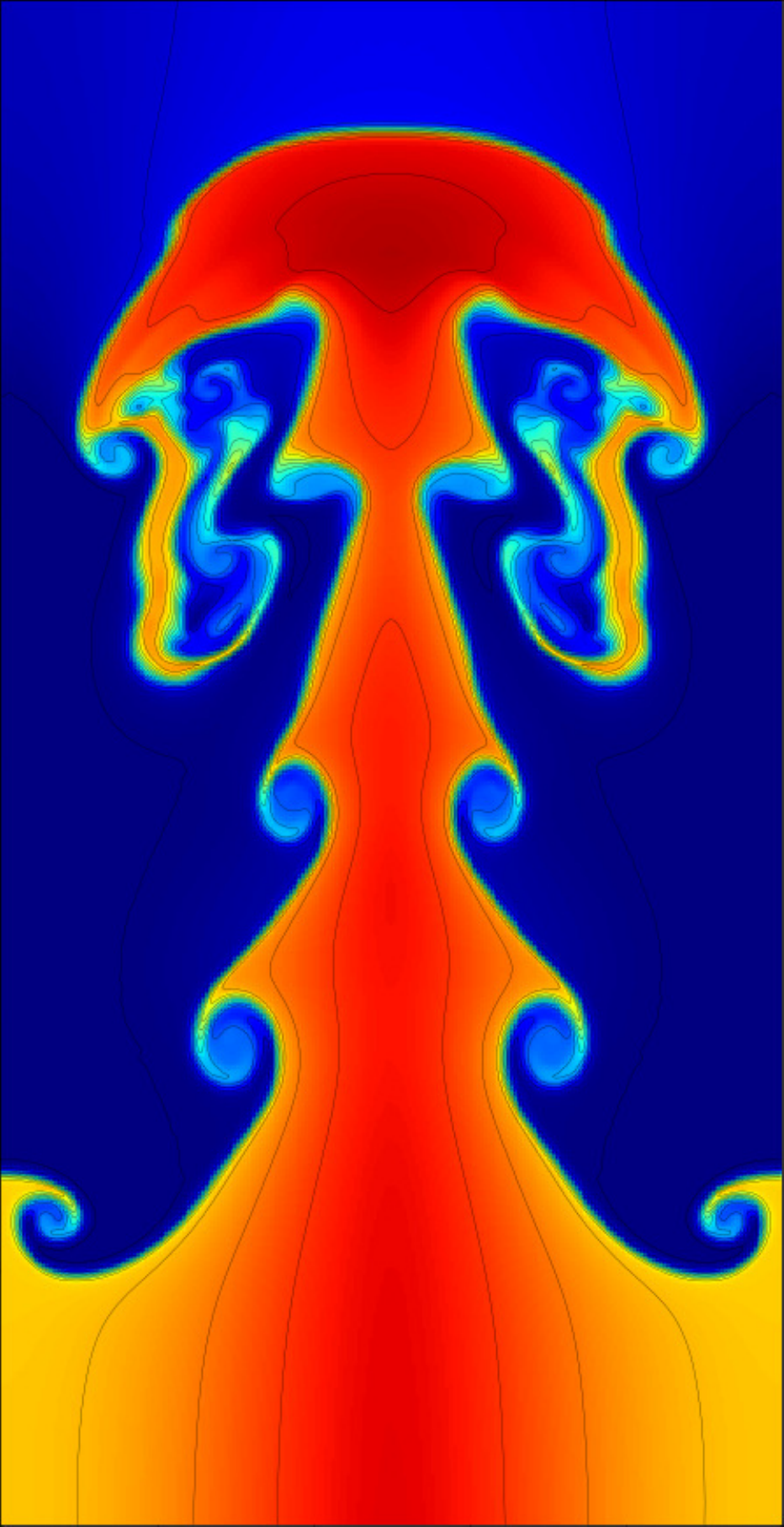}
      \caption{WENO-Z+}
   \end{subfigure}
   \hspace{0.5cm}
   \begin{subfigure}[b]{4.5cm}
      \centering
      \includegraphics[width=\textwidth]{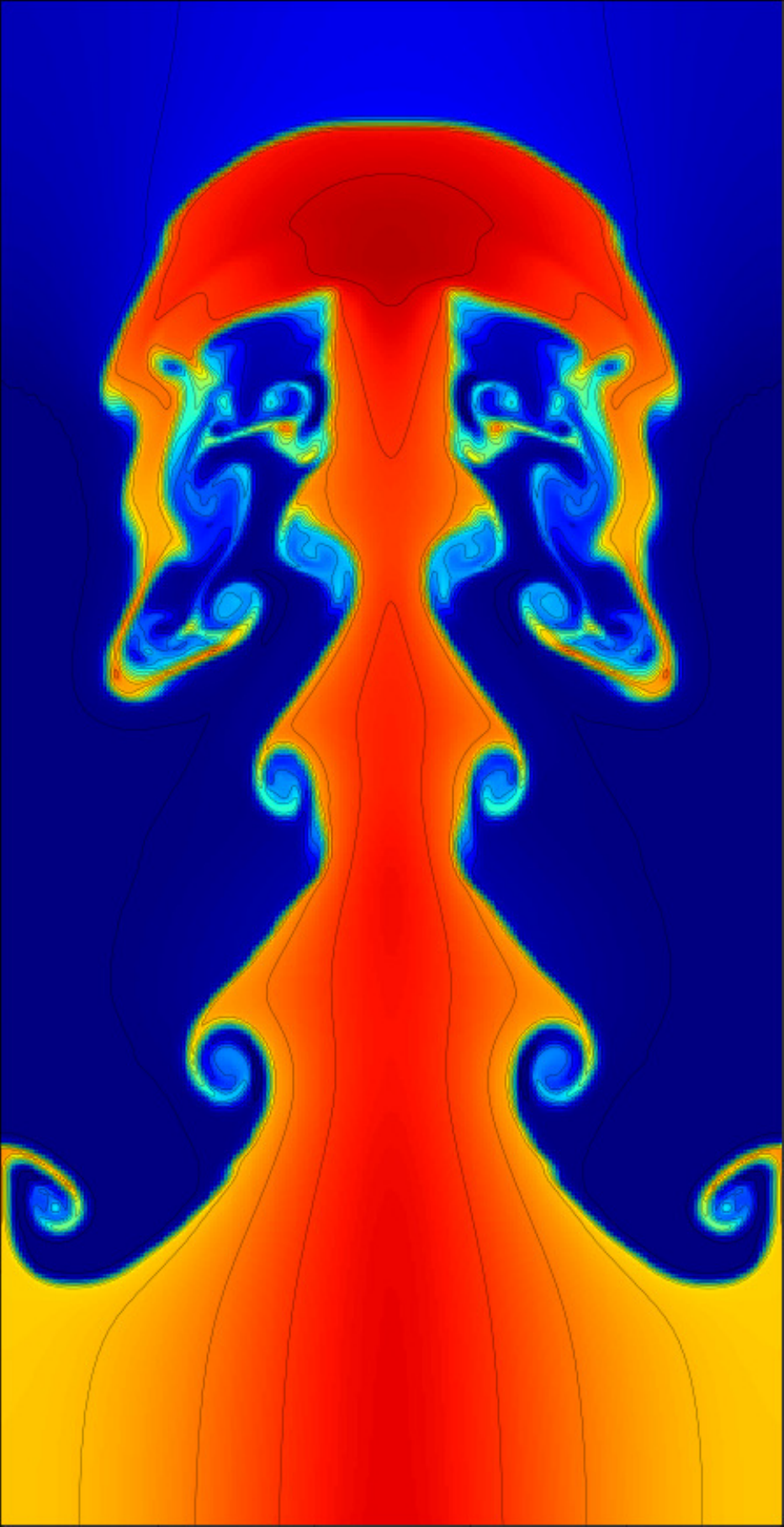}
      \caption{WENO-ZC+}
   \end{subfigure}
   \caption{Numerical solutions of the Rayleigh--Taylor instability test at \(t=1.95\) with $N_{x} \times N_{y} = 3840 \times 950$ and CFL = $0.3$.}
   \label{fig:RTI}
\end{figure}

Figure \ref{fig:RTI} shows the numerical results of the Rayleigh--Taylor instability test for a grid of $3840 \times 950$ points. It can be seen that all schemes present fairly symmetrical solutions with WENO-ZC+ showing more fine structures than the other schemes.

For the Double Mach reflection test \cite{woodward84}, the initial conditions are
\[
   \left(\rho_{0},u_{0},v_{0},p_{0}\right) =
      \begin{cases}
         \left(8,\: 4.125\sqrt{3},\: -4.125,\: 116.5\right), & x<1/6+y/\sqrt{3},\\
         \left(1.4,\: 0,\: 0,\: 1\right), & x\geq1/6+y/\sqrt{3}.
      \end{cases}
\]
The computational domain is \((x, y) \in [0, 4] \times [0, 1]\). In the bottom, the boundary is fixed for \(x \leq 1/6\) and reflexive otherwise. In the top, the boundary is set to follow the Mach 10 oblique shock. In the left and right, inflow and outflow boundaries are respectively used. The final time is \(T = 0.2\).

\begin{figure}
   \centering
   \begin{subfigure}[b]{8cm}
      \centering
      \includegraphics[width=\textwidth]{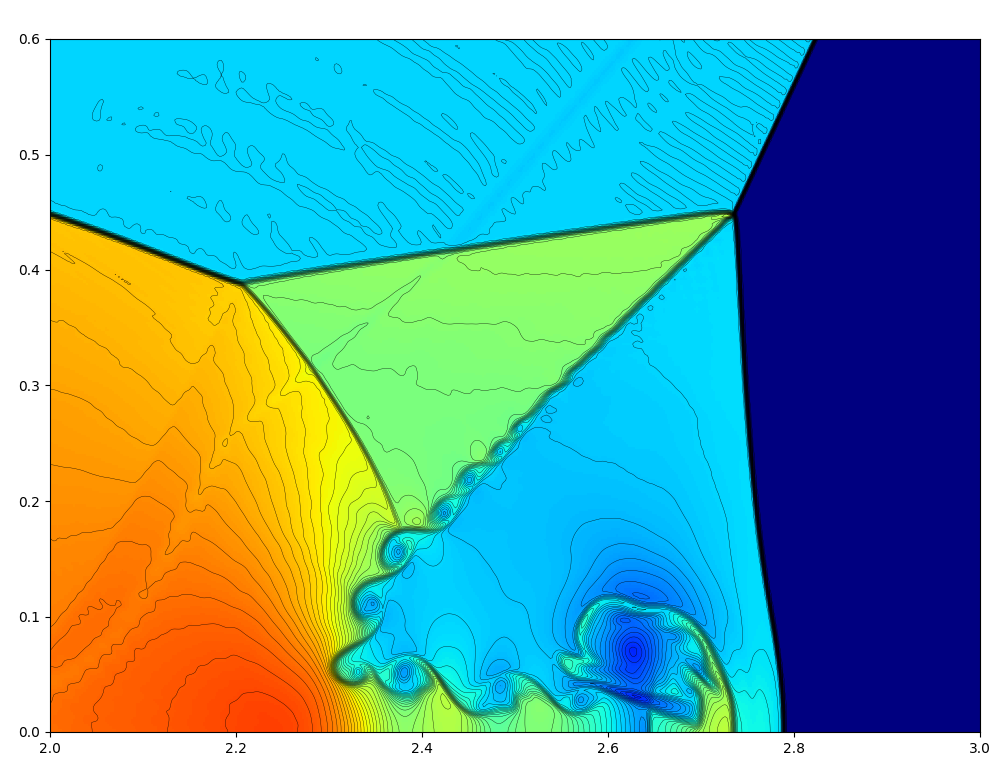}
      \caption{WENO-Z}
   \end{subfigure}
   \begin{subfigure}[b]{8cm}
      \centering
      \includegraphics[width=\textwidth]{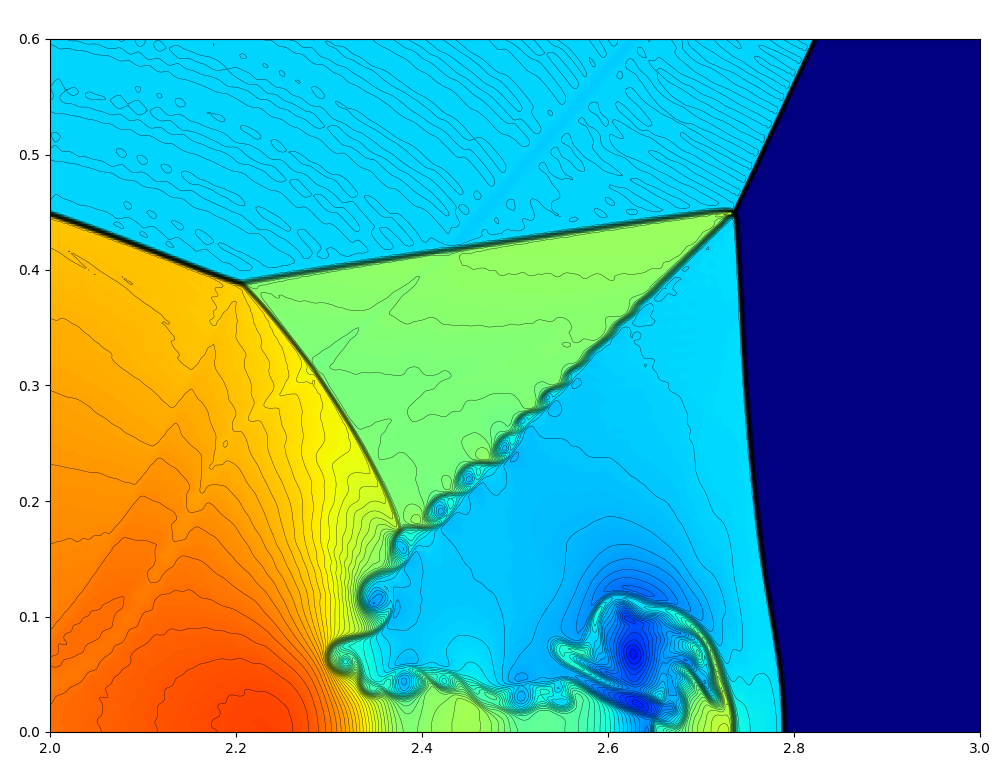}
      \caption{WENO-ZC}
   \end{subfigure}
   \\
   \begin{subfigure}[b]{8cm}
      \centering
      \includegraphics[width=\textwidth]{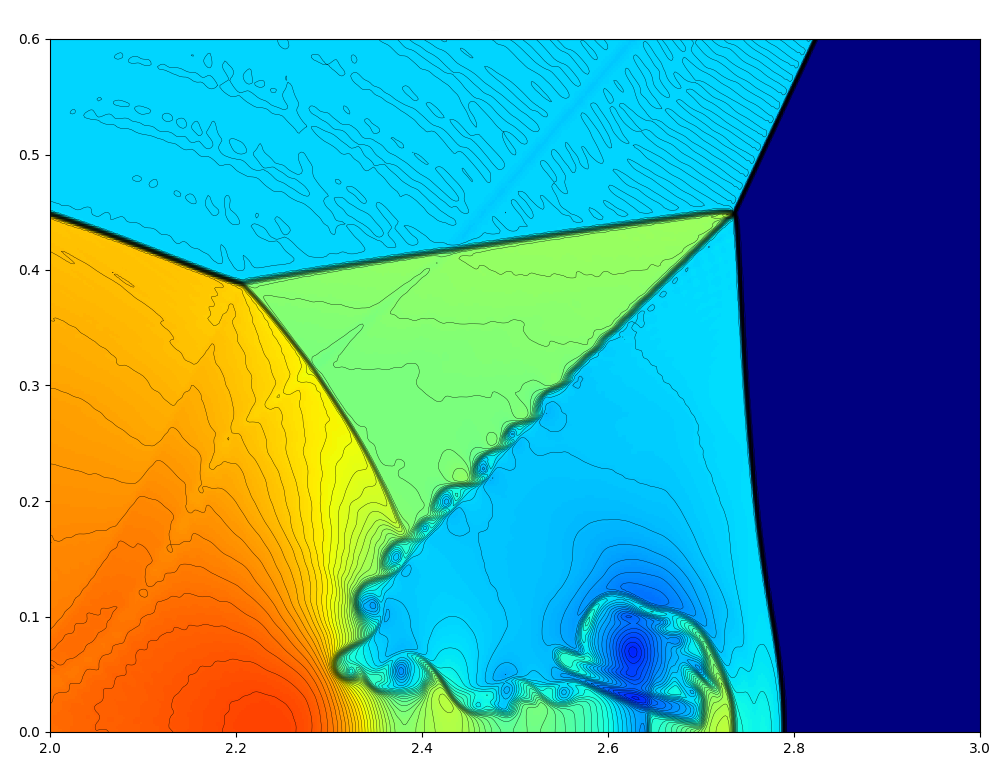}
      \caption{WENO-Z+}
   \end{subfigure}
   \begin{subfigure}[b]{8cm}
      \centering
      \includegraphics[width=\textwidth]{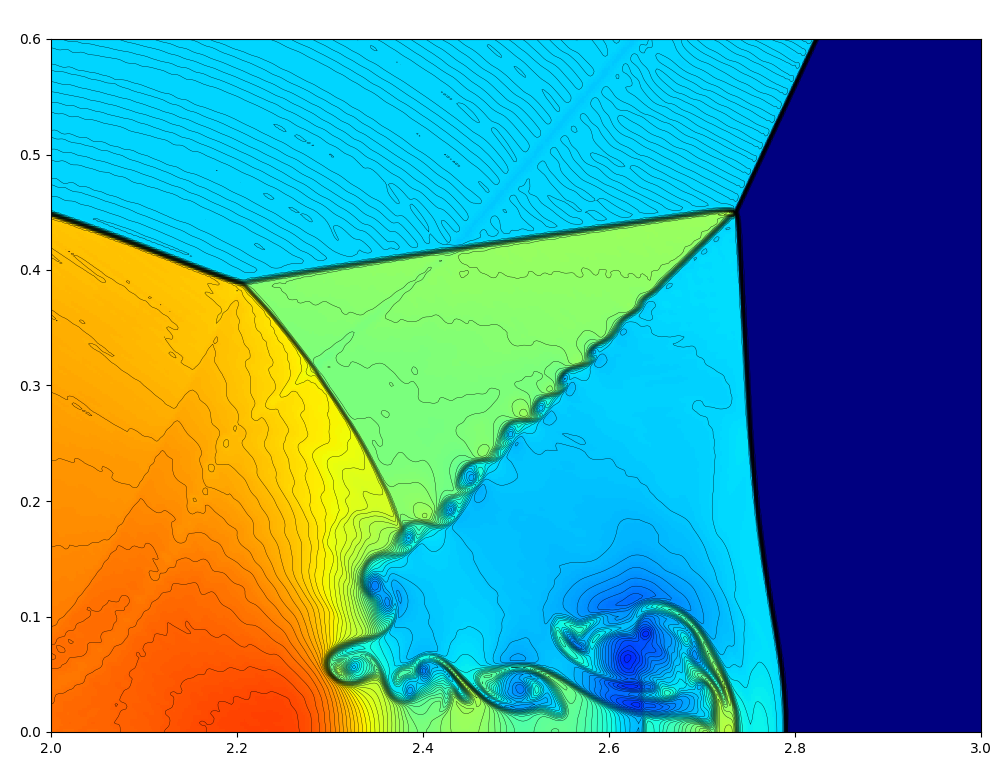}
      \caption{WENO-ZC+}
   \end{subfigure}
   \caption{Numerical solutions of the Double Mach Reflection problem \cite{woodward84} with $N_{x} \times N_{y} = 2000 \times 500$ and CFL = $0.45$.}
   \label{fig:DMR}
\end{figure}

The results of the Double Mach Reflection test for a grid of $2000 \times 500$ points are shown in Figure \ref{fig:DMR} with WENO-ZC+ and WENO-ZC achieving the less dissipative solutions.

%% file: 6-Conclusions.tex
\section{Conclusions} \label{sec:conclusions}

We have proposed a more centered composition  for the nonlinear convex combination of  fifth-order WENO schemes through an overvaluation of the nonlinear component of the weight of the central substencil. The aim was to recover the good dispersive and dissipative properties of the original $5$th order upwind central scheme, which are naturally weakened due to the strengthening of the lateral substencils when high gradients or discontinuities are present in the central substencil. The new scheme, WENO-ZC, fixes a well-known long-term dispersion error, which is typical to most of the standard WENO schemes, also showing less numerical dissipation and a consequential improved shock-capturing ability.

A less dissipative extension to WENO-ZC, similar to WENO-Z+, has also been proposed, and, this time, no ad hoc parameter nor dependence on any power of the grid size was necessary. The extended scheme, dubbed WENO-ZC+, has an equivalent anti-dissipative term, allowing a better representation of the amplitude of curvature features of the solution without the overamplification of WENO-Z+, as shown by the improved results attained with the classical shock-sine wave interaction experiments.  Numerical experiments with the 2D set of the Euler Equations were also carried out. They confirmed the robustness and the improvements of the new weighting strategies compared to the standard WENO-Z and WENO-Z+ schemes. An ADR analysis is also presented to corroborate all these numerical findings.

A thorough analysis of the numerical convergence of the new schemes at critical points in the solution shows that WENO-ZC has the same convergence properties as WENO-Z, i.e., fifth-order is attained at smooth parts of the solution, as well as at first-order critical points. Nevertheless, despite showing better numerical results than WENO-Z+, WENO-ZC+ possesses weaker convergence at critical points of the solution, bringing more to the discussion of the practical importance of this convergence for problems with shocks and discontinuities when compared to other issues, such as keeping the centered structure of the linear scheme when facing an unavoidable suboptimal convergence, as it is the case of the typical problems WENO schemes deal with.

%% file: A-Taylor-Series.tex
\section{Taylor Series Expansions} \label{sec:appendix-taylor}

\subsection[The Numerical Flux Function h(x)]{The Numerical Flux Function \(h(x)\)}

Assume that \(h(x)\) can be written as
\[
   h(x) = a_0 f(x) + a_1 f^{'}(x) \Delta x + a_2 f^{''}(x) \Delta x^2 + a_3 f^{'''}(x) \Delta x^3 + a_4 f^{(4)}(x) \Delta x^4 + 
   \ldots
\]
The Fundamental Theorem of Calculus applied to Eq. \eqref{eq:numerical-flux} gives
\begin{equation}
   f'(x) = \dfrac{h(x + \dx/2) - h(x + \dx/2)}{\dx}.
   \label{eq:h-finite-difference}
\end{equation}
Taking the Taylor series of \(h(x + \dx/2)\) and \(h(x - \dx/2)\) in \eqref{eq:h-finite-difference} allows us to find the coefficients \(a_0\), \(a_1\), \(a_2\), \(\ldots\):
\begin{equation}
   h(x) = f(x) - \frac{1}{24} f^{''}(x) \Delta x^2 + \frac{7}{5760} f^{(4)}(x) \Delta x^4 + \Ord(\Delta x^6).
   \label{eq:h-taylor-x}
\end{equation}

Applying \eqref{eq:h-taylor-x} at \(x_{i\pm\half} \equiv x_i \pm \dx/2\) and taking the Taylor series expansion at \(x = x_i\) gives
\begin{equation}
   h_{i\pm\frac{1}{2}}
      = f_i \pm \frac{1}{2} f_i' \dx + \frac{1}{12} f_i'' \dx^2 - \frac{1}{720} f_i^{(4)} \dx^4 + \Ord(\dx^6).
   \label{eq:taylor-h-fhatk}
\end{equation}
Notice that, with the exception of the coefficient of \(f_i'\), all other coefficients are the same for both \(h_{i-\frac{1}{2}}\) and \(h_{i+\frac{1}{2}}\). This is to be expected since, by \eqref{eq:h-finite-difference}, \(h_{i+\frac{1}{2}} - h_{i-\frac{1}{2}} = f_i' \dx\).

\subsection[The Numerical Flux Approximations fk(x)]{The Numerical Flux Approximations \(\fhat^k(x)\)}

The Taylor series expansion of \(f^{k}(x_{i\pm\frac{1}{2}})\) (Eq. \eqref{eq:fhatk}) at \(x = x_i\) give
   \begin{align}
      f^0_{i-\frac{1}{2}}
         &= f_i - \dfrac{1}{2} f_i' \dx + \dfrac{1}{12} f_i'' \dx^2 - \dfrac{1}{4} f_i''' \dx^3 + \dfrac{61}{144} f_i^{(4)} \dx^{4} - \dfrac{91}{240} f_i^{(5)} \dx^{5} + \Ord(\dx^6), \nonumber \\
      f^1_{i-\frac{1}{2}}
         &= f_i - \dfrac{1}{2} f_i' \dx + \dfrac{1}{12} f_i'' \dx^2 + \dfrac{1}{12} f_i''' \dx^3 - \dfrac{11}{144} f_i^{(4)} \dx^{4} + \dfrac{3}{80} f_i^{(5)} \dx^{5} + \Ord(\dx^6), \nonumber \\
      f^2_{i-\frac{1}{2}}
         &= f_i - \dfrac{1}{2} f_i' \dx + \dfrac{1}{12} f_i'' \dx^2 - \dfrac{1}{12} f_i''' \dx^3 + \dfrac{1}{144} f_i^{(4)} \dx^{4} - \dfrac{1}{240} f_i^{(5)} \dx^{5} + \Ord(\dx^6), \nonumber \\
      f^0_{i+\frac{1}{2}}
         &= f_i + \dfrac{1}{2} f_i' \dx + \dfrac{1}{12} f_i'' \dx^2 - \dfrac{1}{4} f_i''' \dx^3 + \dfrac{25}{144} f_i^{(4)} \dx^{4} - \dfrac{19}{240} f_i^{(5)} \dx^{5} + \Ord(\dx^6), \nonumber \\
      f^1_{i+\frac{1}{2}}
         &= f_i + \dfrac{1}{2} f_i' \dx + \dfrac{1}{12} f_i'' \dx^2 + \dfrac{1}{12} f_i''' \dx^3 + \dfrac{1}{144} f_i^{(4)} \dx^{4} + \dfrac{1}{240} f_i^{(5)} \dx^{5} + \Ord(\dx^6), \nonumber \\
      f^2_{i+\frac{1}{2}}
         &= f_i + \dfrac{1}{2} f_i' \dx + \dfrac{1}{12} f_i'' \dx^2 - \dfrac{1}{12} f_i''' \dx^3 - \dfrac{11}{144} f_i^{(4)} \dx^{4} - \dfrac{3}{80} f_i^{(5)} \dx^{5} + \Ord(\dx^6).
         \label{eq:fhatk-taylor}
   \end{align}

\subsection[The Smoothness Indicators betak]{The Smoothness Indicators \(\beta_k\)}

To find the Taylor series of \(\beta_k(x_{i\pm\half}) \equiv \beta_k^{\pm}\), we first need to write the smoothness indicators \eqref{eq:betaksintegral} in a more convenient form. It has been shown that \cite{jiang96}:
\begin{align}
   \beta_{0}^+ &= \frac{1}{4}(f_{i-2} - 4f_{i-1} + 3f_{i})^2 + \frac{13}{12}(f_{i-2} - 2f_{i-1} + f_{i})^2, \nonumber \\
   \beta_{1}^+ &= \frac{1}{4}(-f_{i-1} + f_{i+1})^2 + \frac{13}{12}(f_{i-1} - 2f_{i} + f_{i+1})^2, \nonumber \\
   \beta_{2}^+ &= \frac{1}{4}(-3f_{i} + 4f_{i+1} - f_{i+2})^2 + \frac{13}{12}(f_{i} - 2f_{i+1} + f_{i+2})^2. \label{eq:betak-quadratic-form}
\end{align}
The expressions for \(\beta_{k}^-\) can be obtained from those of \(\beta_{k}^+\) by shifting the stencils 1 point to the left; for instance, \(\beta_{1}^- = \dfrac{1}{4}(-f_{i-2} + f_{i})^2 + \dfrac{13}{12}(f_{i-2} - 2f_{i-1} + f_{i})^2\).

Taylor series expansions of  \eqref{eq:betak-quadratic-form} at \(x=x_i\) give
\begin{align}
   \beta_{0}^+ &= \big(f_{i}'\big)^{2}\Delta x^{2} + \left(\dfrac{13}{12}\big(f_{i}''\big)^{2}-\dfrac{2}{3}f_{i}'f_{i}'''\right)\Delta x^{4} - \left(\dfrac{13}{6}f_{i}''f_{i}'''-\dfrac{1}{2}f_{i}'f_{i}^{(4)}\right)\Delta x^{5} \nonumber \\
   &\qquad {} + \left(\frac{43}{36}\big(f_{i}'''\big)^{2} + \dfrac{91}{72}f_{i}''f_{i}^{(4)}-\dfrac{7}{30}f_{i}'f_{i}^{(5)}\right)\Delta x^{6} + \Ord(\Delta x^{7}), \nonumber \\
   \beta_{1}^+ &= \big(f_{i}'\big)^{2}\Delta x^{2} + \left(\dfrac{13}{12}\big(f_{i}''\big)^{2}+\dfrac{1}{3}f_{i}'f_{i}'''\right)\Delta x^{4} \nonumber \\
   &\qquad {} + \left(\frac{1}{36}\big(f_{i}'''\big)^{2} + \dfrac{13}{72}f_{i}''f_{i}^{(4)}+\dfrac{1}{60}f_{i}'f_{i}^{(5)}\right)\Delta x^{6} + \Ord(\Delta x^{8}), \nonumber \\
   \beta_{2}^+ &= \big(f_{i}'\big)^{2}\Delta x^{2} + \left(\dfrac{13}{12}\big(f_{i}''\big)^{2}-\dfrac{2}{3}f_{i}'f_{i}'''\right)\Delta x^{4} + \left(\dfrac{13}{6}f_{i}''f_{i}'''-\dfrac{1}{2}f_{i}'f_{i}^{(4)}\right)\Delta x^{5} \nonumber \\
   &\qquad {} + \left(\frac{43}{36}\big(f_{i}'''\big)^{2} + \dfrac{91}{72}f_{i}''f_{i}^{(4)}-\dfrac{7}{30}f_{i}'f_{i}^{(5)}\right)\Delta x^{6} + \Ord(\Delta x^{7}), \nonumber \\
   \beta_{0}^- &= \big(f_{i}'\big)^2 \Delta x^{2} - 2 f_i' f_i'' \Delta x^{3} + \left(\dfrac{25}{12} \big(f_i''\big)^2 + \dfrac{1}{3} f_i' f_i'''\right) \Delta x^{4} - \left(\dfrac{14}{3} f_i'' f_i''' - \dfrac{5}{6} f_i' f_i^{(4)}\right) \Delta x^{5} \nonumber \\
   &\qquad {} + \left(\frac{157}{36}\big(f_{i}'''\big)^{2} + \dfrac{265}{72}f_{i}''f_{i}^{(4)}-\dfrac{59}{60}f_{i}'f_{i}^{(5)}\right)\Delta x^{6} + \Ord(\Delta x^{7}), \nonumber \\
   \beta_{1}^- &= \big(f_{i}'\big)^2 \Delta x^{2} - 2 f_i' f_i'' \Delta x^{3} + \left(\dfrac{25}{12} \big(f_i''\big)^2 + \dfrac{4}{3} f_i' f_i'''\right) \Delta x^{4} - \left(\dfrac{7}{2} f_i'' f_i''' + \dfrac{2}{3} f_i' f_i^{(4)}\right) \Delta x^{5} \nonumber \\
   &\qquad {} + \left(\frac{55}{36}\big(f_{i}'''\big)^{2} + \dfrac{139}{72}f_{i}''f_{i}^{(4)}+\dfrac{4}{15}f_{i}'f_{i}^{(5)}\right)\Delta x^{6} + \Ord(\Delta x^{7}), \nonumber \\
   \beta_{2}^- &= \big(f_{i}'\big)^2 \Delta x^{2} - 2 f_i' f_i'' \Delta x^{3} + \left(\dfrac{25}{12} \big(f_i''\big)^2 + \dfrac{1}{3} f_i' f_i'''\right) \Delta x^{4} - \left(\dfrac{1}{3} f_i'' f_i''' + \dfrac{1}{6} f_i' f_i^{(4)}\right) \Delta x^{5} \nonumber \\
   &\qquad {} + \left(\frac{1}{36}\big(f_{i}'''\big)^{2} + \dfrac{25}{72}f_{i}''f_{i}^{(4)}+\dfrac{1}{60}f_{i}'f_{i}^{(5)}\right)\Delta x^{6} + \Ord(\Delta x^{7}).
   \label{eq:betataylor}
\end{align}

\subsection[The Global Smoothness Indicators tau]{The Global Smoothness Indicators \(\tau\)}
From \eqref{eq:betak-quadratic-form}, the Taylor expansion of the global smoothness indicator $\tau^{\pm}=|\beta_{2}^{\pm}-\beta_{0}^{\pm}|$ at \(x = x_i\) is given by:
\begin{align}
   \tau^+ &= \left|\left(\dfrac{13}{3}f''_{i}f'''_{i} - f_{i}'f_{i}^{(4)}\right)\Delta x^{5} + \left(\dfrac{103}{36}f'''_{i}f^{(4)}_{i} + \dfrac{13}{12}f''_{i}f^{(5)}_{i} - \dfrac{1}{6}f_{i}'f_{i}^{(6)}\right)\Delta x^{7} + \Ord(\Delta x^{9})\right|, \nonumber \\
   \tau^- &= \left|\left(\dfrac{13}{3}f''_{i}f'''_{i} - f_{i}'f_{i}^{(4)}\right)\Delta x^{5} - \left(\dfrac{13}{3}\big(f'''_{i}\big)^2 + \dfrac{10}{3}f''_{i}f^{(4)}_{i} - f_{i}'f_{i}^{(5)}\right)\Delta x^{6}\right. \nonumber \\
      & \qquad\left.+ \left(\dfrac{319}{36}f'''_{i}f^{(4)}_{i} + \dfrac{9}{4}f''_{i}f^{(5)}_{i} - \dfrac{2}{3}f_{i}'f_{i}^{(6)}\right)\Delta x^{7} + \Ord(\Delta x^{8})\right|.
   \label{eq:tautaylor}
\end{align}

%% file: B-Proof-Conditions.tex
\section{Proof of Conditions \ref{cond:sufficient} and \ref{cond:necessary}} \label{sec:proof-conditions}

\begin{proof}
   From \eqref{eq:h-taylor-x} and \eqref{eq:fhatk-taylor}, we can express \(f^k_{i\pm\frac{1}{2}}\) as
   \begin{equation}
      f^k_{i\pm\frac{1}{2}} = h_{i\pm\frac{1}{2}} + A_k f_i''' \dx^3 + \Ord(\dx^4),
      \quad \text{with} \quad
      A_0 = -\dfrac{1}{4}, \:\:
      A_1 = \dfrac{1}{12}, \:\:
      A_2 = -\dfrac{1}{12}.
      \label{eq:fhatk-Ak}
   \end{equation}
   By \eqref{eq:convex combination} and \eqref{eq:taylor-h-fhatk}, we have
   \begin{align*}
      \hat{f}_{i\pm\frac{1}{2}}
         &= \sum_{k=0}^{2}d_{k}\hat{f}_{i\pm\frac{1}{2}}^{k}+\sum_{k=0}^{2}(\omega_{k}^{\pm}-d_{k})\hat{f}_{i\pm\frac{1}{2}}^{k} \\
         &= \left[h_{i\pm\frac{1}{2}} - \frac{1}{60} f_i^{(5)} \Delta x^{5} + \Ord(\Delta x^{6})\right]+\sum_{k=0}^{2}(\omega_{k}^{\pm}-d_{k})\hat{f}_{i\pm\frac{1}{2}}^{k}.
   \end{align*}
   The second term above may be expanded by using \eqref{eq:fhatk-Ak}:
   \begin{gather*}
      \sum_{k=0}^{2} (\omega_{k}^{\pm}-d_{k}) \hat{f}_{i\pm\frac{1}{2}}^{k} 
         = \sum_{k=0}^{2} (\omega_{k}^{\pm}-d_{k}) \left[h_{i\pm\frac{1}{2}} + A_{k}f_i'''\Delta x^{3} + \Ord(\Delta x^{4})\right] \\
         \qquad = h_{i\pm\frac{1}{2}} \sum_{k=0}^{2} (\omega_{k}^{\pm}-d_{k}) + f_i'''\Delta x^{3} \sum_{k=0}^{2}(\omega_{k}^{\pm}-d_{k}) A_{k} + \sum_{k=0}^{2} (\omega_{k}^{\pm}-d_{k}) \Ord(\Delta x^{4}).
   \end{gather*}
   From here, we see that if \(\omega_k^{\pm}\) satisfy the Condition \ref{cond:sufficient}, then we have
   \begin{gather*}
      \sum_{k=0}^{2} (\omega_{k}^{\pm}-d_{k}) \hat{f}_{i\pm\frac{1}{2}}^{k} 
         = \Ord(\dx^6) \quad \therefore \\
      \quad\therefore\quad \dfrac{\fhat_{i+\frac{1}{2}} - \fhat_{i-\frac{1}{2}}}{\dx}
         = \dfrac{h_{i+\frac{1}{2}} - h_{i-\frac{1}{2}} + \Ord(\dx^6)}{\dx}
         = f_i' + \Ord(\dx^5).
   \end{gather*}
   However, Condition \ref{cond:sufficient} is not strictly necessary: Eq. \eqref{eq:condition-dx3} can be relaxed to Eqs. \eqref{eq:idealweightscond}--\eqref{eq:Aks}, yielding the same result. Assuming Condition \ref{cond:necessary} is satisfied, we have
   \begin{gather*}
      \sum_{k=0}^{2} (\omega_{k}^{\pm}-d_{k}) \hat{f}_{i\pm\frac{1}{2}}^{k} 
         = \Ord(\dx^6) + f_i'''\Delta x^{3} \sum_{k=0}^{2}(\omega_{k}^{\pm}-d_{k}) A_{k} + \Ord(\dx^6) \quad \therefore \\
      \quad\therefore\quad \dfrac{\fhat_{i+\frac{1}{2}} - \fhat_{i-\frac{1}{2}}}{\dx}
         = \dfrac{h_{i+\frac{1}{2}} - h_{i-\frac{1}{2}} + f_i'''\Delta x^{3} \sum_{k=0}^{2}(\omega_{k}^{+} - \omega_{k}^{-}) A_{k} + \Ord(\dx^6)}{\dx} = \\
         \qquad\quad = f_i' + \Ord(\dx^5) - \dfrac{f_i'''\Delta x^{2}}{12} \left[3(\omega_{0}^{+}-3\omega_{0}^{-}) - (\omega_{1}^{+}-\omega_{1}^{-}) + (\omega_{2}^{+}-\omega_{2}^{-})\right] \\
         \qquad\quad = f_i' + \Ord(\dx^5).
   \end{gather*}   
\end{proof}

%% file: C-Convergence-WENO-ZDplus.tex
\section{Theoretical Results on the Order of Convergence of WENO-ZC+} \label{sec:appendix-order-convergence-weno-zcplus}

\subsection{In the Absence of Critical Points}

In what follows, consider that \(f\) is a smooth function with no critical points. For simplicity, let us neglect \(\eps\). From Eq. \eqref{eq:weno-zcplus-weights}, we get the following relation for the left \((-)\) and right \((+)\) stencils:
\begin{equation} \label{eq:weno-zcplus-weights-convenient}
   \big(\tau^{\pm} + \overline{\beta}^{\pm}\big) \alpha_{k}^{ZD+(\pm)}
      = d_{k} \left[\tau^{\pm} + \overline{\beta}^{\pm} + \beta_{k}^{\pm} + c_{k}\big(\tau^{\pm} + \overline{\beta}^{\pm}\big)\bigg(\frac{\tau^{\pm}}{\beta_{k}^{\pm}}\bigg)^{p} \bigg(\frac{\tau^{\pm}}{\tau^{\pm} + \overline{\beta}^{\pm}}\bigg)^{p}\right].    
\end{equation}
Eqs. \eqref{eq:weno-zcplus-weights}, \eqref{eq:betataylor}, \eqref{eq:tautaylor} and the results of Section \ref{sec:convergence-analysis-weno-zc} give
\begin{align*}
   \big(\tau^{+} + \overline{\beta}^{+}\big) \alpha_{k}^{ZD+(+)}
      &= d_k \bigg[2 (f_{i}')^{2} \Delta x^{2} + \dfrac{13}{6} (f_i'')^2 \Delta x^{4} + A_k f_i' f_i''' \Delta x^{4} + \Ord(\dx^5) + {} \\
      & \quad {} + c_k \Ord(\dx^2) \Ord(\dx^{3p}) \Ord(\dx^{3p})\bigg] \\
      &= d_k \left[2(f_{i}')^{2} \Delta x^{2} + \dfrac{13}{6} (f_i'')^2 \Delta x^{4} + A_k f_i' f_i''' \Delta x^{4} + \Ord(\dx^5)\right], \\
   \big(\tau^{-} + \overline{\beta}^{-}\big) \alpha_{k}^{ZD+(-)}
      &= d_k \bigg[2 (f_{i}')^{2} \Delta x^{2} - 4 f_i' f_i'' \Delta x^3 + \dfrac{25}{6} (f_i'')^2 \Delta x^4 + B_k f_i' f_i''' \Delta x^4 + {} \\
      & \quad {} + \Ord(\dx^5) + c_k \Ord(\dx^2) \Ord(\dx^{3p}) \Ord(\dx^{3p})\bigg] \\
      &= d_k \bigg[2 (f_{i}')^{2} \Delta x^{2} - 4 f_i' f_i'' \Delta x^3 + \dfrac{25}{6} (f_i'')^2 \Delta x^4 + {} \\
      & \quad {} + B_k f_i' f_i''' \Delta x^4 + \Ord(\dx^5)\bigg],
\end{align*}
where the coefficients \(A_k\) and \(B_k\) are
\begin{align}
   (A_0, A_1, A_2) &= (-1, 0, -1), \label{eq:weno-zcpplus-Ak} \\
   (B_0, B_1, B_2) &= (1, 2, 1). \label{eq:weno-zcpplus-Bk}
\end{align}
This results in
\begin{align}
    \omega_{k}^{ZD+(+)}
       &= \dfrac{\alpha_{k}^{ZD+(+)}}{\sum_{j=0}^2 \alpha_{j}^{ZD+(+)}}
       = \dfrac{\big(\tau^{+} + \overline{\beta}^{+}\big)\alpha_{k}^{ZD+(+)}}{\sum_{j=0}^2 \big(\tau^{+} + \overline{\beta}^{+}\big) \alpha_{j}^{ZD+(+)}} \nonumber \\
       &= \dfrac{d_k \left[2(f_{i}')^{2} \Delta x^{2} + \dfrac{13}{6} (f_i'')^2 \Delta x^{4} + A_k f_i' f_i''' \Delta x^{4} + \Ord(\dx^5)\right]}{\displaystyle\sum_{j=0}^2 d_j \left[2(f_{i}')^{2} \Delta x^{2} + \dfrac{13}{6} (f_i'')^2 \Delta x^{4} + A_j f_i' f_i''' \Delta x^{4} + \Ord(\dx^5)\right]} \nonumber \\
       &= \dfrac{d_k \left[2(f_{i}')^{2} \Delta x^{2} + \dfrac{13}{6} (f_i'')^2 \Delta x^{4} + A_k f_i' f_i''' \Delta x^{4} + \Ord(\dx^5)\right]}{2(f_{i}')^{2} \Delta x^{2} + \dfrac{13}{6} (f_i'')^2 \Delta x^{4} - \dfrac{2}{5} f_i' f_i''' \Delta x^{4} + \Ord(\dx^5)} \nonumber \\
       &= \dfrac{d_k \left[1 + \dfrac{A_k f_i' f_i''' \Delta x^{4}}{2(f_{i}')^{2} \Delta x^{2} + \dfrac{13}{6} (f_i'')^2 \Delta x^{4}} + \Ord(\dx^3)\right]}{1 - \dfrac{\dfrac{2}{5} f_i' f_i''' \Delta x^{4}}{2(f_{i}')^{2} \Delta x^{2} + \dfrac{13}{6} (f_i'')^2 \Delta x^{4}} + \Ord(\dx^3)} \nonumber \\
       &= d_k \left[1 + \left(A_k + \dfrac{2}{5}\right) \dfrac{f_i' f_i'''}{2(f_{i}')^{2} + \dfrac{13}{6} (f_i'')^2 \Delta x^{2}} \Delta x^2 + \Ord(\dx^3)\right] \nonumber \\
       &= d_k \left[1 + \left(A_k + \dfrac{2}{5}\right) \dfrac{f_i'''}{2f_{i}'} \Delta x^2 \left(1 - \dfrac{13}{12} \dfrac{(f_i'')^2}{(f_{i}')^{2}} \Delta x^{2} + \Ord(\dx^4)\right) + \Ord(\dx^3)\right] \nonumber \\
       &= d_k \left[1 + \left(A_k + \dfrac{2}{5}\right) \dfrac{f_i'''}{2 f_{i}'} \Delta x^2 + \Ord(\dx^3)\right]. \label{eq:accuracy-weno-zcplus-plus}
\end{align}
Analogously,
\begin{align}
    \omega_{k}^{ZD+(-)}
       &= \dfrac{d_k \left[2 (f_{i}')^{2} \Delta x^{2} - 4 f_i' f_i'' \Delta x^3 + \dfrac{25}{6} (f_i'')^2 \Delta x^{4} + B_k f_i' f_i''' \Delta x^{4} + \Ord(\dx^5)\right]}{2 (f_{i}')^{2} \Delta x^{2} - 4 f_i' f_i'' \Delta x^3 + \dfrac{25}{6} (f_i'')^2 \Delta x^{4} + \dfrac{8}{5} f_i' f_i''' \Delta x^{4} + \Ord(\dx^5)} \nonumber \\
       &= d_k \left[1 + \left(B_k - \dfrac{8}{5}\right) \dfrac{f_i'''}{2 f_{i}'} \Delta x^2 + \Ord(\dx^3)\right]. \label{eq:accuracy-weno-zcplus-minus}
\end{align}
Eqs. \eqref{eq:weno-zcpplus-Ak}--\eqref{eq:accuracy-weno-zcplus-minus} give
\begin{align*}
   \omega_{0}^{ZD+(\pm)} &= d_k - \dfrac{3}{100} \dfrac{f_i' f_i'''}{(f_{i}')^{2}} \Delta x^2 + \Ord(\dx^3), \\
   \omega_{1}^{ZD+(\pm)} &= d_k + \dfrac{6}{50} \dfrac{f_i' f_i'''}{(f_{i}')^{2}} \Delta x^2 + \Ord(\dx^3), \\
   \omega_{2}^{ZD+(\pm)} &= d_k - \dfrac{9}{100} \dfrac{f_i' f_i'''}{(f_{i}')^{2}} \Delta x^2 + \Ord(\dx^3),
\end{align*}
which implies that WENO-ZC+ satisfies Condition \ref{cond:necessary}.

\subsection[At a Critical Point of Order ncp = 1]{At a Critical Point of Order \(\ncp = 1\)}

Now, suppose that \(x_i\) is a critical point of order \(\ncp = 1\) of \(f\). Eq. \eqref{eq:weno-zcplus-weights-convenient}, together with Eqs. \eqref{eq:weno-zcplus-weights}, \eqref{eq:betataylor}, \eqref{eq:tautaylor} and the results of Section \ref{sec:convergence-analysis-weno-zc}, give
\begin{align*}
   \big(\tau^{+} + \overline{\beta}^{+}\big) \alpha_{k}^{ZD+(+)}
      &= d_k \bigg[\dfrac{13}{6} (f_i'')^2 \dx^4 + \frac{13}{3}\big|f_i'' f_i'''\big| \dx^5 + C_k f_i'' f_i''' \dx^5 + \Ord(\dx^6) + {} \\
      & \quad {} + c_k \Ord(\dx^4) \Ord(\dx^{p}) \Ord(\dx^{p})\bigg] \\
      &= d_k \left[\dfrac{13}{6} (f_i'')^2 \dx^4 + \frac{13}{3}\big|f_i'' f_i'''\big| \dx^5 + C_k f_i'' f_i''' \dx^5 + \Ord(\dx^6)\right], \\
   \big(\tau^{-} + \overline{\beta}^{-}\big) \alpha_{k}^{ZD+(-)}
      &= d_k \bigg[\dfrac{25}{6} (f_i'')^2 \dx^4 + \frac{13}{3}\big|f_i'' f_i'''\big| \dx^5 + D_k f_i'' f_i''' \dx^5 + \Ord(\dx^6) + {} \\
      & \quad {} + c_k \Ord(\dx^4) \Ord(\dx^{p}) \Ord(\dx^{p})\bigg] \\      
      &= d_k \left[\dfrac{25}{6} (f_i'')^2 \dx^4 + \frac{13}{3}\big|f_i'' f_i'''\big| \dx^5 + D_k f_i'' f_i''' \dx^5 + \Ord(\dx^6)\right],
\end{align*}
where the coefficients \(C_k\) and \(D_k\) are
\begin{align}
   (C_0, C_1, C_2) &= \left(-\dfrac{13}{6},\: 0,\: \dfrac{13}{6}\right), \label{eq:weno-zcpplus-Ck} \\
   (D_0, D_1, D_2) &= \left(-\dfrac{15}{2},\: -\dfrac{19}{3},\: -\dfrac{19}{6}\right). \label{eq:weno-zcpplus-Dk}
\end{align}
Analogously to \eqref{eq:accuracy-weno-zcplus-plus}, we have
\begin{align*}
    \omega_{k}^{ZD+(+)}
       &= \dfrac{d_k \left[\dfrac{13}{6} (f_i'')^2 \dx^4 + \dfrac{13}{3}\big|f_i'' f_i'''\big| \dx^5 + C_k f_i'' f_i''' \dx^5 + \Ord(\dx^6)\right]}{\dfrac{13}{6} (f_i'')^2 \dx^4 + \dfrac{13}{3}\big|f_i'' f_i'''\big| \dx^5 + \dfrac{13}{30} f_i'' f_i''' \dx^5 + \Ord(\dx^6)} \nonumber \\
       &= d_k \left[1 + \left(C_k - \dfrac{13}{30}\right) \dfrac{f_i'''}{f_i''} \dx + \Ord(\dx^2)\right], \\ 
    \omega_{k}^{ZD+(-)}
       &= \dfrac{d_k \left[\dfrac{25}{6} (f_i'')^2 \dx^4 + \dfrac{13}{3}\big|f_i'' f_i'''\big| \dx^5 + D_k f_i'' f_i''' \dx^5 + \Ord(\dx^6)\right]}{\dfrac{25}{6} (f_i'')^2 \dx^4 + \dfrac{13}{3}\big|f_i'' f_i'''\big| \dx^5 - \dfrac{11}{2} f_i'' f_i''' \dx^5 + \Ord(\dx^6)} \nonumber \\
       &= d_k \left[1 + \left(D_k + \dfrac{11}{2}\right) \dfrac{f_i'''}{f_i''} \dx + \Ord(\dx^2)\right].
\end{align*}
As such, the WENO-ZC+ scheme does not achieve the optimal order 5 in the presence of critical points.